\documentstyle[11pt]{article}
\newtheorem{tm}{Theorem}[subsection]
\newtheorem{lm}[tm]{Lemma}
\newtheorem{pr}[tm]{Proposition}
\newtheorem{rmk}[tm]{Remark}
\newtheorem{cor}[tm]{Corollary}
\newtheorem{ex}[tm]{Example}
\newtheorem{exe}[tm]{Exercise}

\newtheorem{??}[tm]{Question}
\newtheorem{defi}[tm]{Definition}

\font\tenmsb=msbm10
\font\sevenmsb=msbm7
\font\fivemsb=msbm5

\newfam\msbfam
\textfont\msbfam=\tenmsb
\scriptfont\msbfam=\sevenmsb
\scriptscriptfont\msbfam=\fivemsb
\def\Bbb#1{{\fam\msbfam #1}}

\font\teneufm=eufm10
\font\seveneufm=eufm7
\font\fiveeufm=eufm5
\newfam\eufmfam
\textfont\eufmfam=\teneufm
\scriptfont\eufmfam=\seveneufm
\scriptscriptfont\eufmfam=\fiveeufm

\def\lorw{\longrightarrow}
\newcommand\n{\noindent}

\newcommand\ci{\cite}
\newcommand\s{\sigma}
\newcommand\rat{{\Bbb Q}}
\newcommand\comp{{\Bbb C}}
\newcommand\real{{\Bbb R}}
\newcommand\zed{{\Bbb Z}}

\newcommand\pn[1]{{\Bbb P}^{#1}}

\newcommand\blacksquare{{\hspace*{\fill} $\Box$}}

\newcommand\e{\epsilon}

\newcommand{\im}{ \hbox{\rm Im} }

\newcommand{\ke}{ \hbox{\rm Ker} }

\newcommand{\ov}[1]{ \overline{#1}}

\newcommand{\surj}{ -\!\!\!-\!\!\!\gg}
\newcommand{\inj}{\hookrightarrow}

\title{Lectures on the Hodge theory of projective manifolds}
\author{
Mark Andrea A.  de Cataldo\thanks{
Partially supported by N.S.F. Grant DMS 0202321, N.S.A. Grant 
MDA904-02-1-0100.}
}

\date{}

\begin{document}\maketitle

    \vspace{3cm}
   
 \begin{verse}
   {\em    \hspace{ 7cm}
 After the mountain, \\
\hspace{ 7cm} another mountain} \\

\hspace{ 7cm} (traditional Korean saying)
\end{verse}



\newpage
As my friend and colleague Luca Migliorini
once wrote to me, the topology of algebraic 
varieties is a mystery and a miracle.

\medskip
These  lectures are an attempt to introduce the reader
to the Hodge theory of algebraic varieties.

\n
The geometric implications of Hodge theory
for a compact oriented  manifold
become progressively richer and more beautiful
as one specializes from Riemannian, to  complex,  to K\"ahler
and finally to projective manifolds.

\n
The structure of these lectures tries to reflect this
fact.

\medskip
I delivered these eight  one-hour lectures at the
July 22 -- July 27, 2003
Summer School on Hodge Theory at the Byeonsan Peninsula 
in South Korea.

\n
I would like to thank  Professor JongHae Keum
and Dr. Byungheup Jun for organizing the event,
Professor 
Jun-Muk Hwang and Professor Yongnam Lee for supporting it.

\n
I would like to thank all of those who have attended
the lectures for the warm atmosphere I have found
that has made my stay a wonderful experience.
In particular, I would like to thank Professor Dong-Kwan Shin
for showing me with great humor some aspects of the 
Korean culture.

\medskip
Beyond the choice of topics, exercises and exposition style, nothing
in these written-up version of the lectures is original.
The reader is assumed to have some familiarity with smooth and complex
manifolds. 
These lectures are not self-contained and at times
a 
remark or an exercise require knowledge
of notions and facts which are
not covered here. I see no harm in that, but 
rather as an encouragement to the reader to explore the subjects
involved.
The location of the
exercises in the lectures  is a suggestion to the reader 
to solve them before proceeding to the  theoretical
facts that follow.

\medskip
The table of contents should be self-explanatory.
The only exception is $\S$\ref{apv} where
I discuss, in a simple example, a  technique
for studying the class map
for homology classes on the fibers of a map
and one 
for  approximating
a certain kind of primitive vectors. These  techniques
have been
introduced in 
\ci{demig1} and \ci{demig2}.

\bigskip
I would like to thank  Fiammetta Battaglia,
Gabriele LaNave,
Jungool Lee and Luca Migliorini for suggesting improvements
and correcting some mistakes. The remaining mistakes
and the shortcomings in the exposition are mine.

\newpage

$\;$

\vspace{7cm}

\centerline{Alla memoria di  Meeyoung }

\newpage

\tableofcontents

\newpage

\section{Lecture 1: Calculus on smooth manifolds}
\label{l1}
We introduce some  basic structures on a finite dimensional
real vector space with a metric: the metric on the exterior algebra,
orientations, the volume element and 
the  star isomorphism.
We introduce smooth differential forms, de Rham cohomology 
groups, state the de Rham Theorem and discuss Weil's sheaf-theoretic 
approach
to this theorem. We briefly discuss  Riemannian metrics, 
orientations on manifolds and integration
of top forms on oriented manifolds.

\bigskip
References for this lecture are
 \ci{wa}, 
 \ci{b-t} and \ci{dem1}.

 \bigskip
Manifolds are assumed to be {\em connected} 
and to  satisfy the second axiom of 
countability. One of the advantages of the second assumption is that it 
implies the existence of partitions of unity, an important
tool for the study of smooth manifolds.

\subsection{The Euclidean structure on the exterior algebra}
\label{tim}
Let $V$  be a $m-$dimensional real  
vector  space with inner product  
$\langle\, , \,\rangle ,$
i.e. $\langle \, , \, \rangle: V \times V \lorw \real$ 
is a symmetric, positive definite bilinear form on $V.$

\bigskip
Let
$$
\Lambda(V) =\bigoplus_{0 \leq p
\leq m}{ \Lambda^{p}(V) }
$$
be the  exterior algebra associated with $V.$

\medskip
\n
If $ \{ e_{1}, \ldots , e_{m} \}$ is a basis for $V,$
then  the elements $e_{I}=e_{i_1} \wedge \ldots  \wedge  e_{i_p} ,$
where $I =(i_{1}, \ldots,  i_{p} )$ 
ranges in the corresponding set of multi-indices 
with $ 1 \leq  i_{1} < \ldots < i_{p} \leq m ,$
form a basis for $\Lambda^{p}(V).$

\medskip
\n
The elements of
$\Lambda^{p}(V)$ can be seen
as the alternating $p-$linear form
on $V^{*}$ as follows:
$$
v_{1} \wedge \ldots \wedge v_{p} 
\, = \, 
 \, \sum_{\nu \in S_{p}}{  
\, \e(\nu) \, v_{\nu_1} \otimes \ldots
\otimes v_{\nu_{p}}},
$$ 
where $S_{p}$ is the symmetric group
in $p$ elements and $\e(\nu)$  is the sign
 of the  permutation $\nu.$

\bigskip
When dealing with exterior algebras
it is costumary to write summations such as
$\sum_{I}{ u_{I} \,  e_{I}} $ or $\sum_{|I|=p}{ u_{I}\, e_{I} }$
meaning that the summation is over the set of ordered multi-indices
$I$
as above and $u_I \in \real.$

\bigskip
There is a natural inner product on 
$\Lambda(V)$ defined by declaring any two distinct  spaces
$\Lambda^{p}(V)$ and $\Lambda^{p'}(V)$  mutually orthogonal
and setting
\begin{equation}
    \label{eq1}
\langle \, v_{1} \wedge \ldots 
\wedge v_{p} \, , \,  w_{1} \wedge \ldots \wedge
w_{p} \, \rangle : = \det{  || \, \langle v_{j}, w_{k}\rangle \, ||}, 
\qquad v_j, w_k \in V.
\end{equation}
If the basis $\{e_{1}, \ldots , e_{m} \}$ for $V$  is orthonormal,
then so is the corresponding one for $\Lambda (V).$

\bigskip
\begin{exe}
    \label{ittwo}
    {\rm
Let $V = V' \oplus V''$ be an orthogonal direct sum decomposition,
$v', \, w' \, \in \Lambda^{p'}(V')$ and $v'', \, w'' \, \in \Lambda^{p''}(V'').$
Show that
$$
\langle \langle v' \wedge v''  \, , \, w' \wedge w'' \rangle \rangle
\, = \, \langle \langle v', w' \rangle \rangle \, \langle
\langle v'', w'' \rangle \rangle.
$$
}
\end{exe}
\subsection{The star isomorphism on $\Lambda(V)$}
\label{tsi}
By definition, $\Lambda^{0}(V)=\real.$
The positive half-line $\real^{+} \subseteq \Lambda^{0}(V)$
is  defined without ambiguity.

\bigskip
The real vector space $\Lambda^{m}(V)$  is one-dimensional and 
$\Lambda^{m}(V) \setminus \{0 \}$ has two connected components,
i.e. two half-lines. 

\medskip
\n
However, unlike the case
of $\Lambda^{0}(V),$ there is no canonical way to distinguish
either of them. 

\n
\medskip
A choice is required. This choice
gives rise to an isometry $\Lambda^{0}(V) \simeq
\Lambda^{m}(V).$ 

\bigskip
The $\star$ operator is the operator that naturally 
arises when we want to complete the picture
with linear isometries $\Lambda^{p}(V) \simeq \Lambda^{m-p}(V).$

\begin{defi}
    \label{or}
    {\rm 
    ({\bf Orientations of $V$})
 The choice of 
a connected component  $\Lambda^{m}(V)^{+}$
of $\Lambda^{m}(V) \setminus \{0 \}$ is called 
an {\em orientation} on $V.$
}
\end{defi}

\bigskip

Let $V$ be oriented and $\{ e_i \}$
be an ordered, orthonormal basis such that $e_{1}\wedge \ldots
\wedge e_{m} \in \Lambda^{m}(V)^{+}.$
This element is uniquely defined since any two  bases
as above are related by a orthogonal matrix with determinant $+1.$

\begin{defi}
    \label{volel}
    {\rm ({\bf The volume element})
    The  vector 
\begin{equation}
    \label{vlmt}
    \fbox{$
dV\, : = \,  e_1 \wedge \ldots \wedge e_m \, \in \,  \Lambda^{m}(V)^{+}  
$}
\end{equation}
 is called  {\em the volume element} associated with the oriented
$(V, \langle \, , \,\rangle ,\Lambda^{m}(V)^{+}) .$
}
\end{defi}

\bigskip
\begin{defi}
    \label{doso}
{\rm ({\bf The $\star$ operator}) 
The {\em $\star$ operator}  is the unique linear isomorphism
$$
 \star \,  :\,  \Lambda (V)\,  \simeq \,   \Lambda (V)   
$$
defined by the properties
$$
\star \,  :\,  \Lambda^{p}(V) \, \simeq \, \Lambda^{m-p}(V),
$$
$$
\fbox{$ \,  u\wedge \star \,v
\, =\,  \langle \,  u,v  \, \rangle  \, dV,
\quad \quad \forall \; u,\, v\,  \in \, \Lambda^{p}(V), \; 
\forall p.$}
$$
}
\end{defi}

\bigskip
The $\star$ operator depends on the inner product {\em and}
on the chosen orientation.

\bigskip
Let us check that the operator $\star$ exists and is unique.

\smallskip
\n
Consider the non-degenerate pairing
$$
\Lambda^{p}(V) \times \Lambda^{m-p}(V) \lorw \real,\qquad
(u,w) \lorw  (u\wedge w)/dV. 
$$
More explicitly, let $\{e_{i} \} $ be an oriented  orthonormal
basis,
 $dV = e_{1} \wedge \ldots  \wedge e_{m}$ be the volume element,
 $u= \sum_{I}{ u_{I}\,e_{I}},$
$w= \sum_{J}{ w_{J}\,e_{J}}.$ Then
$$
u \wedge w \, =\,  \e(I,CI) \, \sum_{I}{ u_{I}\, w_{CI} \, dV},
$$
where, if
$I$ is an ordered set of indices, then $CI$ is the
complementary ordered set of indices and $\e(I,CI)$ is the sign of the 
permutation $(I,CI)$ of the ordered  set $\{1, \ldots, m  \}.$

\n
Let $v= \sum_{I}{ v_{I}\, e_{I}}$ and define
\begin{equation}
    \label{defstarok}
    \fbox{$
 \star \, v\,  := \, \sum_{I}{ \e (I,CI) \, v_{I} \, e_{CI}}.$}
\end{equation}
We have that 
$$
u\wedge \star\,  v \, = \, \sum_{I}{u_{I}\, e_{I}}
\wedge \sum_{I'}{ \e (I',CI') \, v_{I'} \, e_{CI'} } \, = \,
\sum_{I}{ u_{I} \, v_{I} \,\e(I,CI) \, e_{I}\wedge e_{CI}}\, = \,
\langle u, v \rangle \, dV.
$$
This shows the existence of $\star$, which seems to depend on the choice
of the orthonormal basis.  

\n
Assume $\star'$ is another such operator.
Then
$ u \wedge (\star - \star')(v) =
(\langle u,v \rangle - \langle u,v \rangle  )  \, dV = 0$ 
for every $u$ and $v$
which implies 
that $\star =\star'.$

\begin{exe}
\label{exonst}
{\rm  Verify the following statements.

\n
The $\star$ operator is a linear isometry.

\n
Let 
$dV= e_1 \wedge \ldots \wedge e_m$ be the  volume
element for the fixed orientation.
$$
 \star \,  (1)\, =\,  dV, \quad 
 \star\, (dV)\, =\,   1, 
\quad  \star\,  (e_{1}\wedge \ldots \wedge e_{p}) =  e_{p+1}\wedge
\ldots \wedge e_{m}, 
$$
\begin{equation}
\label{stst}
\star \, \star_{| \Lambda^{p}(V)} \,  =\, 
 (-1)^{p(m-p)} \, Id_{\Lambda^{p}(V)} ,
\end{equation}
$$
 \langle u, v\rangle  \,  =\,  \star ( v \wedge \star u) \,  = \,  \star (u
\wedge  \star \, v), 
\quad \forall \, u,\, v \,  \in\,  \Lambda^{p}(V), \quad
\forall \; 0 \leq p \leq m.
$$
}
\end{exe}

\subsection{The tangent  and cotangent bundles of a smooth manifold}
\label{tbosm}
A smooth manifold $M$ of dimension 
$m$ comes equipped with natural smooth
vector bundles. 

\bigskip
Let $(U; x_{1}, \ldots, x_{m})$ be a chart
centered at a point  $q \in M.$

\begin{itemize}

\item $T_M$ the {\em tangent bundle} of $M.$ The fiber
$T_{M,q}$ can be identified with the linear span
$\real \langle \partial_{x_{1}}, \ldots, \partial_{x_{m}} \rangle.$

\item
$T^*_M$ the {\em cotangent bundle} of $M.$ 
Let $\{ dx_{i} \}$ be the dual basis of
the basis
$\{ \partial_{x_{i}} \}.$ 
The fiber $T^{*}_{M,q}$ can be indetified with the span
$\real \langle d x_{1}, \ldots, d x_{m} \rangle.$

\item $\Lambda^p(T^*_M)$ the {\em $p-$th exterior
bundle} of $T^*_M$. The fiber
$\Lambda^p(T^*_M)_q = \Lambda^p(T^*_{M,q})$ can be identified
with the linear span
$\real \langle \{ d x_I \}_{|I|=p}  \rangle.$

\item $\Lambda (T^*_M) \,:= \, \oplus_{p=0}^m{ \Lambda^p(T^*_M)   }$ 
the {\em exterior algebra bundle} of $M.$

\end{itemize}

\bigskip
As it is costumary and by 
slight abuse of notation, one can use the symbols
 $\partial_{x_{i}}$ and
$dx_{j}$ to denote  the corresponding sections
over $U$ of the tangent and cotangent
bundles.

\medskip
\n
In this case, the collections 
$\{ \partial_{x_{i}} \}$ and
$\{ dx_{i} \}$  form {\em local frames} for the bundles in question.

\bigskip

\begin{exe}
\label{trftc}
{\rm
Let $M$ be a smooth manifold of dimension $m.$
Take any definition in the literature
for the tangent and cotangent bundles and verify the
assertions that follow. 
Let $(U; x_{1}, \ldots, x_{m})$ and
$(U'; x'_{1}, \ldots, x'_{m})$
be two local charts
centered at $q \in M.$
Show that the transition functions
$\tau_{U'U}(x)  : U \times \real^{m} \simeq
U' \times \real^{m}$ ($\gamma_{U'U}(x),$ resp.) for the tangent bundle
$T_{M}$ (cotangent bundle $T^{*}_{M},$ resp.)
are given, using the column notation for vectors
in $\real^{m}$, by
\[
\tau_{U'U}(x) =  \left[ J \left( \begin{array}{c} x' \\ x \end{array} 
\right) (x) \right] ^{t},
\qquad
\gamma_{U'U}(x) = \left[ J \left( \begin{array}{c} x' \\ x \end{array} 
\right) (x) \right]^{-1},
\qquad
\]
where
\[
J \left( \begin{array}{c} x' \\ x \end{array} 
\right) (x) =  
\left(
\begin{array}{ccc}
  \frac{\partial x'_{1} }{\partial x_{1}}(x) & \ldots & 
    \frac{\partial x'_{m} }{\partial x_{1}}(x) \\
    \vdots & \ldots & \vdots \\
     \frac{\partial x'_{1} }{\partial x_{m}}(x) & \ldots & 
    \frac{\partial x'_{m} }{\partial x_{m}}(x)
    \end{array}
    \right).
    \]
 }
\end{exe}

\begin{exe}
\label{tgsp}
{\rm
Determine the tangent bundle of the sphere $S^{n} \subseteq 
\real^{n+1}$  given by the equation
$\sum_{j=1}^{n+1}{ x_{j}^{2}} =1.$
Let $S^{n-1} \subseteq S^{n}$ be an  ``equatorial'' embedding,
i.e. obtained by intersecting $S^{n}$ with the hyperplane
$x_{n+1}=0.$ Study the ``normal bundle'' exact sequence
$$
0 \lorw T_{S^{n-1}} \lorw (T_{S^{n}})_{|S^{n-1}}  \lorw N_{S^{n-1}, 
S^{n}} \lorw 0.
$$
}
\end{exe}

\bigskip

\subsection{The de Rham cohomology groups}
\label{tdrcg}

\begin{defi}
    \label{defofep}
 {\rm ({\bf $p-$forms}) 
The  elements of the real vector space
$$
\fbox{$
E^{p}(M) \, : = \; C^{\infty}( M, \Lambda^{p}( T^{*}_{M} ) ) $}
$$
of smooth real-valued  
sections of the vector bundle $\Lambda^{p}( T^{*}_{M} )$
 are called (smooth) {\em $p-$forms}
on $M.$
}
\end{defi}

\bigskip

Let $d: E^{p}(M) \lorw E^{p+1}(M)$ denote the exterior derivation
of differential forms.

\bigskip

\begin{exe}
    \label{disd}
{\rm
A $p-$form  $u$ on $M$ can be written
on $U$ as $u = \sum_{|I| =p}{ u_{I} dx_{I}}.$

\n
Show that exterior derivation of forms is well-defined.
More precisely, show that if one defines, locally on the chart 
$(U;x_{1}, \ldots, x_{m}),$
$$
du\, : =\,   \sum_{j}{ \frac{\partial u_{I}}{\partial x_{j}} \, dx_{j}
\wedge dx_{I}},
$$
then $du$ is independent of the choice of coordinates.
}
\end{exe}

\bigskip
\begin{defi}
    \label{defco}
    {\rm ({\bf Complexes, cohomology  of a complex})
A {\em complex} is a sequence of maps of vector spaces
$$  
\ldots \lorw  V^{i-1} \stackrel{\delta^{i-1}}\lorw V^i 
\stackrel{\delta^i}\lorw V^{i+1} \lorw \ldots, \quad i \in \zed,
$$
also denoted by $(V^{\bullet}, \delta ),$
such that $\delta^{i} \circ \delta^{i-1} =0$ for every 
index $i,$ i.e. such that $\im{ \,\delta^{i-1}} \subseteq
\ke{\, \delta^i}.$ 

\n
The vector spaces $H^i(V^{\bullet}, \delta )=:
\ke{\,\delta^i}/ \im{ \,\delta^{i-1}}$ are called the {\em
cohomology groups} of the complex. 

\n
A complex is said to be {\em exact
at $i$} if  $\im{ \,\delta^{i-1}}=\ke{ \,\delta^i},$ i.e.
if $H^i( V^{\bullet}, \delta   )=0$, and {\em exact} if it 
is exact for every $i \in \zed.$

}
\end{defi}

\bigskip
\begin{exe}
\label{drcdef}
{\rm ({\bf The de Rham complex})
Show that $d^{2}=0$ so  that we get the so-called
{\em de Rham}  complex of  vector spaces,
\begin{equation}
    \label{drcplx}
    \fbox{$
0 \lorw E^{0}(M) \stackrel{d}\lorw E^{1}(M) \stackrel{d}\lorw
\ldots \stackrel{d}\lorw E^{m-1}(M) \stackrel{d}\lorw
E^{m}(M) \lorw 0. $}
\end{equation}
}
\end{exe}

\bigskip
\begin{defi}
 \label{cve}
 {\rm ({\bf Closed/exact})
 A $p-$form $u$ is said to be {\em closed} if $du=0$ and is said to be 
{\em exact} if there exists $v \in E^{p-1}(M)$ such that
$dv =u.$
}
\end{defi}

\bigskip
\begin{exe}
\label{execlo}
{\rm
Let $u$ and $v$ be closed forms. Show that $u \wedge v$ is closed.
Assume, in addition,  that $v$ is exact and  show that $u\wedge v$ is exact.
}
\end{exe}

\bigskip
\begin{exe}
\label{414}
{\rm
Let
$$
u \, = \, (2x+ y \,\cos{xy})\,dx \, + \, (x\,\cos{xy}) \, dy 
$$
on $\real^{2}.$
Show  that $u$ is exact.
What is the integral of $u$ along any closed curve in $\real^{2}$?
}
\end{exe}

\bigskip
\begin{exe}
\label{415}
{\rm
Let
$$
u \, = \, \frac{1}{2\pi} \, \frac{ x\,dy \, - \, y\, dx}{ x^{2} \, + 
\, y^{2}}
$$
on $\real^{2} \setminus \{0 \}.$
Show that $u$ is closed.
Compute the integral of $u$ over the unit circle $S^{1}.$
Is $u$ exact?
Is $u_{|S^{1}}$ exact?
}
\end{exe}

\bigskip
\begin{exe}
\label{416}
{\rm
(a) Prove that every closed $1-$form on $S^{2}$ is exact.

\smallskip
\n
(b) Let
$$
u \, = \frac{ x \, dy\wedge dz - y \,  dx\wedge dz + 
z \, dx \wedge dy}{(x^{2} +y^{2} +z^{2})^{3/2}}
$$
on $\real^{3} \setminus \{0 \}.$
Show  that $u$ is closed.

\smallskip
\n
(c) Evaluate $\int_{S^{2}}{ u}.$
Conclude that $u$ is not exact.

\smallskip
\n
(d) 
Let
$$
u \, = \frac{ x_{1} \, dx_{1} \ldots + x_{n} \, dx_{n} }{(x_{1}^{2} + 
\ldots  +x_{n}^{2})^{n/2}}
$$
on $R^{n} \setminus \{0 \}.$
Show  that   $\star \, u$ is closed.

\smallskip
\n
(e) Evaluate $\int_{S^{n-1}}{ \star \, u}.$
Is  $\star \, u$  exact?
}
\end{exe}

The de Rham complex (\ref{drcplx})  is {\em not} exact and
its  deviation from exactness  is an important
invariant of $M$ and is measured by
the so-called de Rham cohomology groups of $M;$
see Theorem \ref{drit}.

\bigskip
\begin{defi}
\label{defofdrcg}
{\rm ({\bf The de Rham cohomology groups})
The real {\em de Rham cohomology groups} 
$H_{dR}^{\bullet}(M, \real)$  of $M$ are the 
cohomology groups of the complex
(\ref{drcplx}), i.e. 
$$
\fbox{$
H^{p}_{dR}(M, \real) \,  \simeq \,  \frac{\mbox{closed $p-$forms on M}}{
\mbox{exact $p-$forms 
on M.}}
$}
$$
}
\end{defi}

\bigskip
The de Rham complex 
is {\em locally exact} on $M$
by virtue of the important

\bigskip
\begin{tm}
    \label{pl}
    ({\bf Poincar\'e Lemma})
    Let $p>0.$
A  closed $p-$form $u$   on $M$   is locally exact, i.e. 
for every $q \in M$ there exists
an open neighborhood $U$ of $q$ and $v \in E^{p-1}(U)$ such that
$$
\fbox{
$u_{|U} \, = \, dv.$}
$$
    \end{tm}
    {\rm Proof.} See \ci{b-t}, $\S 4$, \ci{wa}, $\S 4.18.$ \blacksquare

\bigskip
The following result of de Rham's is fundamental.
In what follows the algebra structures are given by the wedge and cup 
products, respectively.

\bigskip
\begin{tm}
\label{drit}
({\bf The de Rham Theorem})
Let $M$ be a not necessarily orientable smooth manifold.
There is a  canonical isomorphisms of $\,\real-$algebras 
$$
\fbox{$
H^{\bullet}_{dR}(M,\real) \, \simeq \, H^{\bullet}(M, \real).
$}
$$ 
\end{tm}

\bigskip
\begin{rmk}
\label{alliso}
{\rm 
Theorem \ref{drit} is obtained using integration
over differentiable simplices and then using
the various canonical identifications
between real singular cohomology and real differentiable singular cohomology.
Of course, one also has canonical isomorphisms with
real Alexander-Spanier cohomology, sheaf cohomology with coefficients
in the locally constant sheaf $\real_M$ and  
\v{C}ech cohomology
with coefficients in $\real_M.$
See \ci{wa}, $\S$ $4.7, \, 4.17,\, 5.34-5.38,$ $5.43-5.45,$
and \ci{mu} 5.23--5.28. 
}
\end{rmk}

\bigskip
\begin{rmk}
    \label{weil}
    {\rm  ({\bf Sheaf-theoretic de Rham Theorem})
It is important to know that what above, and more, admits a 
sheaf-theoretic re-formulation, due to A. Weil. Here is a sketch of this 
re-formulation.
Let ${\cal E}^{p}_{M}$ be the sheaf of germs
of smooth $p-$forms on $M,$ ${\cal E}_{M}:= {\cal E}^{0}_{M}$
be the sheaf of germs of smooth functions on $M,$ $\real_{M}$
be the sheaf of  germs of locally constant functions
on $M.$ 
Due to the existence
of partitions of unity,
the sheaf ${\cal E}_{M}$ is what one calls a {\em fine} sheaf
(see \ci{gh}, page 42). Any 
sheaf of ${\cal E}_{M}-$modules is fine.
In particular, the sheaves 
${\cal E}^{p}(M)$ are fine.
Fine sheaves have  trivial higher sheaf cohomology groups. 
The operator $d$ is defined locally
and gives rise to maps of sheaves
$d: {\cal E}^p_M \lorw {\cal E}^{p+1}_M.$
In this context, the Poincar\'e Lemma \ref{pl} implies that
the complex of sheaves
\begin{equation}
    \label{cosss}
0 \lorw \real_{M} \hookrightarrow {\cal E}^{0}_{M}\stackrel{d}\lorw
{\cal E}^{1}_{M} \stackrel{d}\lorw  \ldots \stackrel{d}\lorw
{\cal E}^{m-1} \stackrel{d}\lorw {\cal E}^{m}_{M} \lorw 0
\end{equation}
is exact.
In short, ${\real}_{M} \lorw ({\cal E}^{\bullet}, d^{\bullet})$ is a 
resolution of $\real_M$ by fine sheaves.
The sheaf cohomology of a sheaf $F$ on $X$ is defined
by considering a resolution
$F \lorw I^{\bullet}$ of $F$ by injective sheaves
and setting $H^{p}_{Sheaf}(X, F) := H^p ( H^{0}(X, I^{\bullet} )).$
Choosing another injective  resolution gives
canonically isomorphic
sheaf cohomology groups. 
One can take a resolution of $F$ by fine sheaves as well.
Taking global  sections in (\ref{cosss})
and ignoring $\real_M,$ we obtain the 
de Rham complex (\ref{drcplx}).
It follows that
there is a
natural isomorphism between
the sheaf cohomology of $\real_M$ and de Rham cohomology:
$$
H^p_{Sheaf}(X, \real_M) \simeq H^p (E^{\bullet},d)=: H^p_{dR}(X, \real).
$$
The de Rham Theorem \ref{drit} follows
from the natural  identification of the standard
 singular cohomology with real coefficients
with the sheaf cohomology of ${\real}_{M}.$ See \ci{wa}. 

\n
Surprisingly, while the left-hand-side is a topological invariant,
the right-hand-side depends on the smooth structure. In other
words: {\em the ``number'' of linearly independent
smooth closed $p-$forms  which are not exact is a topological
invariant, independent of the smooth structure on $M.$}
The advantage of this sheaf-theoretic approach
to the de Rham Isomorphism Theorem 
is that it gives rise to a variety of isomorphisms
in different contexts, as soon as one has Poincar\'e  Lemma-type
results.
One other example is the sheaf-theoretic
proof of the Dolbeault Isomorphism Theorem 
based on the Grothendieck-Dolbeault Lemma 
\ref{gdl}. See Remark \ref{weildol}.
}
\end{rmk}
 
\bigskip

\begin{rmk}
    \label{curr}
    {\rm 
It is important to know that  
the complexes (\ref{cosss}) and (\ref{drcplx}) have counterparts in
the theory of currents on a manifold:
\begin{equation}
    \label{cosssc}
0 \lorw \real_{M} \hookrightarrow {\cal D}^{0}_{M}\stackrel{d}\lorw
{\cal D}^{1}_{M} \stackrel{d}\lorw  \ldots \stackrel{d}\lorw
{\cal D}^{m-1} \stackrel{d}\lorw {\cal D}^{m}_{M} \lorw 0,
\end{equation}
\begin{equation}
    \label{covsooc}
0  \lorw D^{0}(M) \stackrel{d}\lorw D^{1}(M) 
\stackrel{d}\lorw
\ldots \stackrel{d}\lorw D^{m-1}(M) \stackrel{d}\lorw
D^{m}(M) \lorw 0
\end{equation}
yielding  results analogous to the ones of Remark
\ref{weil}. See \ci{gh}, $\S$3 for a quick 
introduction to currents.

\n
The space of currents $D^{p}(X, \real)$ is defined as the topological
dual of the space of compactly supported smooth
$p-$forms
 $E^{m-p}_{c}(M, \real)$ endowed with the 
$C^{\infty}-$topology.

\n
The exterior derivation 
$$
d \, = \, D^{p}(M, \real) \lorw D^{p+1}(M, \real) 
$$
is defined by
$$
d \, (T) (u) : = (-1)^{p+1} T (d u) .
$$

\n
Currents have a local nature which allows to write them as differential
forms with distribution coefficients. 
If $M$ is oriented, then
any differential $p-$form $u$  is also a $p-$current
via the assignment $ v\lorw \int_{M}{u \wedge v},$
$ v \in E^{m-p}_{c}(M, \real).$ 

\n
In this  case, the complexes (\ref{cosss}) and (\ref{cosssc}), 
and (\ref{drcplx})
and (\ref{covsooc}) are quasi-isomorphic, i.e. the natural injection
induces isomoprhisms at the cohomology level.

\n 
By integration, a piecewise smooth oriented $(n-p)-$chain
in $M$ 
gives a $p-$current  in $D^{p}(M).$

\n
Importantly, any closed analytic subvariety
$V$ of complex codimension
$d$ of a complex manifold $X$ 
gives rise, via integration, to a $2d-$current
$\int_{V}$ in $D^{2d}(X, \real).$ 
By Stokes' Theorem for analytic varieties \ci{gh}, page 33, 
such a current is closed and therefore gives rise to
a cohomology class $[\int_{V}]$ which coincides with
the fundamental class $[V]$ of $V \subseteq X.$
}
\end{rmk}

\begin{exe}
\label{dtc}
{\rm 
Let $X \to Y$ be the blowing up of a smooth complex surface
at a point $y \in Y.$
Show that 
$$
Rf_{*} \real_{X} \, \simeq \, Rf_{*} {\cal D}^{\bullet}_{X} 
\, \simeq \, 
f_{*} {\cal D}^{\bullet}_{X}.
$$
Use the currents  of integration
$$
\int_{X} \;  \qquad \mbox{and}
\qquad
 \int_{f^{-1}(y) } 
$$
to construct  an isomorphism
$$
Rf_{*} \real_{X} \, \simeq \, \real_{Y} [0] \, \oplus \, 
\real_{y}[-2].
$$
}
\end{exe}

\subsection{Riemannian metrics}
\label{rm}
A {\em Riemannian metric} $g$ on a smooth manifold $M$
is the datum of 
 a smoothly-varying positive inner product $g(-,-)_{q}$
on the fibers of $T_{M,q}$ of the tangent bundle
of $M.$ This means that, using a chart $(U;x),$ 
the functions 
$$
g_{jk}(q)=  g(\partial_{x_{j}}, 
\partial_{x_{k}} )_{q}
$$ 
are smooth on $U.$  

\medskip
\n
 A  Riemannian metric  induces an isomorphism
 of vector bundles $T_M \simeq_{g} T^*_M$
and we can naturally  define a metric on $T^*_M.$
\begin{exe}
\label{dualmetric}
{\rm ({\bf The dual metric}) 
Verify the following assertions.

\n
Let $E$ be a real, rank $r$ vector bundle on $M,$
$(U;x)$ and $(V;y)$ be two local charts around $q \in M,$
 $\tau_{U}: E_{|U} \lorw U\times \real^{r} \longleftarrow
V \times \real^{r}:\tau_{V}$ be two trivializations of $E$
on the charts and $\tau_{VU}(x)= \tau_{V}\circ \tau_{U}^{-1}:
(U\cap V) \times \real^{r} \lorw (V\cap U) \times \real^{r}$
be the transition functions which we view as an invertible $r\times r$
matrix of real-valued functions of $x \in U \cap V.$

\n
A metric $g$ on $E$ is a smoothly-varying inner product
on the fibers of $E.$ It can be viewed as a 
symmetric bilinear map of vector bundles $E \times E  \lorw \real_{M}.$
Any bilinear map $b: E\times E  \lorw \real_{M}$
can be written, on a chart $(U,x),$ as a $r\times r$ matrix
$b_{U}$ subject to the relation 
$$
b_{V} \, =\,  (\tau_{VU}^{-1})^{t} \, b_{U} \, \tau_{VU}^{-1}.
$$
If $b$ is non-degenerate, then its  local representations
$b_{U}$ are non-degenerate and we can define a non-degenerate
bilinear form on the dual vector bundle $E,$
$b^{*}  : E^{*}\times E^{*} \lorw \real_{M}$
by setting
$$
b^{*}_{U}\, :  = \,  b_{U}^{-1},
$$
where it is understood that we are using
as bases to represent $b^{*}$
the 
dual bases to the ones  employed to represent $b$ on $U.$

\n
If $b$ is symmetric, then so is $b^{*}.$ If $b$ is positive definite,
then
so is $b^{*}.$ If $\{e_{1}, \ldots, \e_{m} \}$
is an orthonormal frame for $E,$ with respect to
$g,$ then the dual frame $ \{ e_{1}^{*}, \ldots , \e_{m}^{*} \}$
is orthonormal for $g^{*}.$

\n
Conclude
that if $(M,g)$ is a Riemannian manifold, then
there is a unique metric $g^{*}$ on $T^{*}_{M}$
such that the  isomorphism induced by the metric (every bilinear form
$V \times V \lorw \real,$ induces a linear map $V \lorw V^{*}$)
$T_M \simeq_{g} T^*_M$ is an isometry on every fiber.
}
\end{exe}

\subsection{Partitions of unity}
\label{pofun}
Recall that a  (smooth) {\em partition of unity} on $M$ is a collection
$\{\rho_{\alpha} \}$ of non-negative
smooth functions on $M$ such that
the sum $\sum_{\alpha}{ \rho_{\alpha} }$ is locally finite
on $M$ and adds-up to the value $1.$ 
This means that for every $q \in M$ there is a neighborhood
$U$ of $q$ in $M$ such that ${\rho_{\alpha}}_{|U} \equiv 0$
for all but finitely many indices $\alpha$ so that
the sum $\sum_{\alpha}{ \rho_{\alpha}(q)}$ is finite and adds to $1.$

\bigskip

\begin{defi}
    \label{pou}
    {\rm ({\bf Partition of unity })
    Let $\{U_{\alpha}\}_{\alpha \in A}$ be an open covering of $M.$
    A {\em partition of unity subordinate to the covering
    $\{U_{\alpha}\}_{\alpha \in A}$}
    is
    a partition of unity $\{ \rho_{\alpha} \} $
    such that the support of each $\rho_{\alpha}$ is contained
    in $U_{\alpha}.$
    }
    \end{defi} 
 
    \bigskip
    Note that, on a non-compact manifold, 
    it is not possible  in general to have a partition 
    of unity subordinate to a given covering and such that
    the functions $\rho_{\alpha}$ have compact support, e.g.
    the covering of $\real$ given by the single open set
    $\real.$

    \begin{tm}
    \label{poucs}
    ({\bf Existence of partitions of unity})
     Let $\{U_{\alpha}\}_{\alpha \in A}$ be an open covering of
    $M.$ Then there are
    
    \smallskip
    \n
   a) a partition of unity subordinate to $\{ U_{\alpha} \}$ and

   \smallskip
   \n
    b)  a partition of unity
$\{ \rho_{j} \}_{j\in J},$ where $J \neq A$ in general, such that
    (i) the support of every  $\rho_{j}$ is compact 
     and (ii)
   for every index  $j$ there is an index
    $\alpha$ such that $supp \, (\rho_{j})
    \subseteq U_{\alpha}. $
    \end{tm} 
    {\em Proof.} \ci{wa}, page 10.

\bigskip
\begin{exe}
    \label{rme}
    {\rm
    Prove that every smooth manifold admits Riemannian 
    metrics on it. (Hint: use partitions of unity. See \ci{b-t}, page 
    42.)
    }
    \end{exe}

\subsection{Orientation and integration}
\label{oai}
Given any smooth manifold $M,$
the space 
$ \Lambda^{m}(T^{*}_M) \setminus M,$ where $M$ is embedded in
the total space of the line bundle
$\Lambda^{m}(T^{*}_M)$ as the zero section, has at most two 
connected components. The reader should verify this.

\bigskip
\begin{defi}
    \label{defori}
    {\rm ({\bf Orientation}) 
    A smooth manifold $M$ is said to be {\em orientable} if
 $ \Lambda^{m}(T^{*}_M) \setminus M$ has two connected components,
{\em non-orientable}
otherwise.

\n
If $M$ is orientable, then
the choice of a connected component of
 $ \Lambda^{m}(T^{*}_M) \setminus M$  is called
 an {\em orientation} of $M$ which is then said to be
 {\em oriented.}
 }
 \end{defi}

 \bigskip
 \begin{ex}
 \label{exoreuc}
 {\rm 
 The {\em standard orientation} on  $\real^{m}= \{ (x_{1}, \ldots, 
 x_{m}) \, | \; x_{j}  \in \real \}$ 
  is the one associated with 
  $$dx_{1}\wedge \ldots \wedge
 dx_{m}.$$
 Similarly, the torus $\real^{m}/\zed^{m} \simeq (S^{1})^{m}$
 is oriented using the  form  above descended to the 
 torus.
 }
 \end{ex}
 
 \bigskip

\begin{exe}
\label{saor}
{\rm Let $M$ be a differentiable manifold of dimension $m.$
Show that the following three statements are equivalent.
See \ci{wa}, $\S 4.2$ and \ci{b-t}, $\S I.3$.

\medskip
\n
(a) $M$ is orientable.

\smallskip
\n
(b) There is a collection of coordinate system
$\{ U^{\alpha}, x^{\alpha} \}$ 
such that 
$$ \det{ || \frac{\partial x^{\beta}_{i}}{ \partial
x^{\alpha}_{j} } (x^{\alpha} ) ||  } \,  >\, 0, \qquad
\mbox{on} \;
U^{\alpha} \cap U^{\beta}.
$$ 

\smallskip
\n
(c) There is a nowhere vanishing $m-$form on $M.$
}
\end{exe}

\bigskip
\begin{rmk}
    \label{popa}
    {\rm 
A collection  of charts with the property
specified  in Exercise \ref{saor}.b is    called
{\em orientation-preserving}.
Such a collection determines uniquely an orientation preserving 
collection of coordinate systems containing it which is maximal with 
respect to inclusion. Such a maximal collection is called
an {\em orientation preserving atlas} for the orientable $M.$
Note that composing with the automorphism $x \to -x$ we obtain a
different orientation preserving atlas.
Given an orientation preserving atlas, we can orient
$M$ by choosing the connected component of $\Lambda^{m}(T_{M}^{*}) 
\setminus M$  containing the vectors $\tau_{\alpha}^{*} \,
(dx_{1}\wedge \ldots \wedge dx_{m}) (q),$ where $\tau_{\alpha} : 
U_{\alpha} \lorw \real^{m}$ is a chart in the atlas, $q \in U_{\alpha}$
and $dx_{1} \wedge \ldots \wedge dx_{n}$ is the canonical orientation
of $\real^{m}.$ Of course, one could choose the opposite one as well.
However, if $M$ is oriented, we consider the orientation preserving
atlas and coordinate systems
that agree with the orientation in the sense mentioned above.
}
\end{rmk}
\bigskip

 \begin{exe}
     \label{spheres}
     {\rm
Let $f: \real^{m+1} \lorw \real$ be a smooth function
such that  $df (q) \neq 0$ for every $q \in \real^{m+1} $
such that $f(q)=0.$

\n
Show that the equation $f=0$ defines a  possibly empty
collection
of smooth and connected 
$m-$dimen\-sio\-nal submanifolds
of $\real^{m+1}.$

\n
Let $M$ be a connected component of the  locus $(f=0).$

\n
Let 
$$
T   \, := \; \{\,  (q,x) \, | \;
x \cdot df (q) =0  \, \} \subseteq M \times R^{m+1} .
$$
Show that 
$$
T_{M} \,  \simeq \, T 
$$
as smooth manifolds and as vector bundles over $M.$

\n
Show that if $M$ is unique, then  $\real^{m+1} \setminus M$ 
has two connected components $A^{+}$ and $A^{-}$ determined by the sign of
the values of $f$ so that $A := A^{-} \coprod M$ is an
$m-$manifold with 
boundary.

\n
Show that $M$ carries a natural orientation associated with $A:$
let $q \in M,$ $\{v_{1}, \ldots, v_{m} \}$
be vectors in  $T_{\real^{m+1},q}^{*}$ such that
 $(df(q), v_{1}, \ldots, v_{m})$ is an 
oriented basis for $\real^{m+1}$ at $q.$ Check that one indeed gets
an orientation for $M$ by taking the vectors $v$ restricted to
$T_{M,q}^{*}.$

\n
Compare what above with the notion of induced orientation on the 
boundary of a manifold with boundary in \ci{b-t}, page 31. They 
coincide.

\n
Compare what above with \ci{wa}, Exercise 4.1
which asks to prove that a codimension one submanifold
of $\real^{m+1}$ is orientable if and only if
there is  a nowhere-vanishing
smooth normal vector field (i.e. a section
$\nu$ of $T_{\real^{m+1}}$
defined over $M \subseteq \real^{m+1}$, smooth over $M$, such that
$\nu (q) \perp T_{M,q}$ for every $q \in M$). 

\n
Compute the volume form associated with the natural orientation on $M$
and the Riemannian metric on $M$ induced by the Euclidean metric on 
$\real^{m+1}.$ The answer is that it is the restriction to $M$
of the contraction of $dx_{1}\wedge \ldots \wedge dx_{m+1}$
with the oriented unit normal vector field along $M.$ See
\ci{wa}, Exercise 20.a. See \ci{wa}, page 61 and \ci{dem1}, page 22
for the definition, properties and explicit form of the contraction
operation.

\n
Make all the above explicit in the case when  
$M= S^{m}$ is  defined by the equation $f:= \sum_{j=1}^{m+1}{x_{j}^{2}} 
=1$ in the Euclidean space $\real^{m+1}.$

\n
Compute everything using the usual spherical coordinates of calculus
books.

\n Solve Exercise 4.20.b of \ci{wa} which gives the 
volume form for surfaces $(x,y, \varphi(x,y))$ in $\real^{3}.$
}
 \end{exe}

\bigskip
\begin{exe}
     \label{exeonor}
     {\rm 
Prove that the M\"obius strip, the Klein
bottle and $\real \pn{2}$ are non-orientable.

\n
Prove
that $\real \pn{m}$ is orientable iff
$m$ is odd. (Hint: the antipodal map $S^{m} \lorw S^{m}$ 
is orientation preserving
iff $m$ is odd.)
}
\end{exe}

\bigskip
Let $M$ be an oriented manifold of dimension $m.$
In particular, $M$ admits an orientation preserving atlas.
It is using this atlas and partitions of unity
that we can define the operation of integrating
$m-$forms with compact support on $M$
and verify Stokes' Theorem \ref{st}.

\bigskip

Let us discuss how integration is defined.
See \ci{b-t}, $\S I.3$ (complemented by Remark \ref{popa} above)
or \ci{wa}, $\S$4.8.

\medskip

\n
Let $\s$ be a  $m-$form on $M$  with compact support  contained 
in a chart $(U;x)$ with trivialization $\tau_{U} : U \lorw \real^{m}.$

\n
The $m-$form $(\tau_{U}^{-1})^{*} \s = f(x) \, dx_{1} \wedge \ldots
\wedge dx_{m},$ for a unique compactly supported  real-valued function $f(x)$
on $\real^{m}.$

\n
Define
$$
\int_{M}^{U}{ \s} \, := \, \int_{\real^{m}}{ f(x) \, dx} ,
$$ 
where $dx$ denotes 
the Lebesgue measure on $\real^{m}$ and the left-hand-side is
the Riemann-Stijlties-Lebesgue integral of $f(x).$

\n
The trouble with this definition is that it depends on the chosen chart.
If $(V,y)$ is another such chart and $T: y \to x$ is the patching 
function, then
$$
(\tau_{V}^{-1})^{*} \s \, = \,  g(y) \, dy_{1} \wedge \ldots
\wedge dy_{m}\, = \,  f(T(y)) \, J(T)(y) \,dy_{1} \wedge \ldots
\wedge dy_{m}
$$
so that
$$
\int_{M}^{V}{ \s}\,= \, 
\int_{\real^{m}}{ g(y) \, dy} \, = \,
\int_{\real^{m}}{ f(T(y)) \, J(T)(y) \, dy}.
$$
On the other hand, the change of variable formula for the Riemann 
integral gives
$$
\int_{M}^{U}{ \s} \,:= \, \int_{\real^{m}}{ f(x) \, dx} \, =\, 
\int_{\real^{m}}{ f(T(y)) \, |J(T)(y)| \, dy}.
$$
It follows that the two definitions agree iff
$J(T)(y) >0$ on the support of $\s.$

\bigskip
This suggests that we can define the integral of $m-$forms
only in the presence of an orientation-preserving atlas.

\bigskip
Let $\omega $ be a smooth $m-$form with compact support
on $M.$

\n
Let $\{ U_{\alpha} \}$ be an orientation
preserving collection of coordinate systems (see Remark \ref{popa})
with  orientation preserving trivializations
$\tau_{\alpha}: U_{\alpha} \simeq \real^{m}$ and $\{\rho_{\alpha} \}$ be a 
partition of unity subordinate to the covering $\{ U_{\alpha} \}.$

\n
The forms $\rho_{\alpha} \omega$ have compact support
contained in $U_{\alpha}.$ 

\n
Define
\begin{equation}
    \label{defofint}
    \fbox{$
\int_{M}{ \omega} = \sum_{\alpha}{ \int_{U_{\alpha}}{  \rho_{\alpha}
\omega }}. $}
\end{equation}
One checks that it is well-defined by a simple
partition of unity argument.

\bigskip
Changing the orientation simply changes the signs of the integrals.

\begin{exe}
    \label{welldefint}
{\rm
Check that $\int_{M}{\omega}$, as defined in (\ref{defofint}),
is well-defined, i.e. that the definition does not depend
on the covering and partition of unity chosen.
}
\end{exe}

\bigskip
Note that one can integrate compactly supported
functions $f$  on any Riemannian manifold
so that, if the manifold is also oriented, then the integral coincides
with the integral of the top-form $\star  \,f.$ See \ci{wa}, page 150.

\bigskip
Complex manifolds are always orientable and are usually oriented
using a standard orientation; see Proposition \ref{cmao}.
Using the standard orientation, if $\omega$ is the $(1,1)-$form
associated with the Fubini-Study metric of the complex projective 
space $\pn{n}$
(see Exercise \ref{fsm}), then $\int_{\pn{n}}{ \omega^n } =1.$

\begin{tm}
\label{st}
({\bf Stokes' Theorem} (simple version)). 
Let $M$ be an oriented manifold of dimension $M$
and $u$ be a $(m-1)-$form on $M$ with compact support.
Then
$$
\fbox{$
\int_{M}{du}=0.
$}
$$
\end{tm}
{\em Proof.} See \ci{b-t}, $\S I.3.5$, \ci{wa}, page 148.
\blacksquare

\bigskip
We have stated a weak version of Stokes' Theorem. See the references 
above for the complete version. See also \ci{wa}, page 150-151 and 
Exercise 4.4 for its re-formulation as the Divergence Theorem.

\bigskip
If a Riemannian metric $g$ on the oriented manifold $M$ is given,
then
we have the notion of {\em Riemannian volume element}
associated with  $(M,g)$ and 
with the orientation. It is the unique $m-$form $dV$
on $M$ such that, for every $q \in M,$
$dV_q$ is the 
volume element of  Definition \ref{volel}
for the dual metric $g^{*}_{q}$ on $T^{*}_{M,q}.$ 

\bigskip
If the integral $\int_{M}{dV}$ converges, then  its value is positive
and it is called the {\em volume} of the oriented 
Riemannian manifold.

\bigskip
\begin{exe}
\label{exvol}
{\rm 
Compute the volume of the unit sphere $S^{m}
\subseteq {\real}^{m+1}$ with respect to the metric induced by
the Euclidean metric and the induced orientation.

\n
Do the same  for the tori $\real^{m}/\zed^{m},$ 
where the lattice  $\zed^{m} \subseteq \real^{m}$
is generated by the vectors $(0, \ldots, r_{j}, \ldots, 0)$, $r_{j}
\in \real^{+},$ $ 1 \leq j \leq m.$
}
\end{exe}

\newpage
\section{Lecture 2: The Hodge theory of  a smooth, oriented, compact
Riemannian manifold}
\label{l2}
We define the inner product $\langle \langle \, , \, \rangle \rangle,$
the adjoint $d^{\star},$ the Laplacian, harmonic forms.
We state the Hodge Orthogonal Decomposition Theorem for 
compact oriented Riemannian manifolds, deduce the Hodge 
Isomorphism Theorem and the Poincar\'e Duality Theorem.

\bigskip
Some references for this section are
\ci{dem1},  \ci{wa} and \ci{gh}.

\subsection{The adjoint of $d: \; d^{\star}$}
\label{taod}
Let $(M,g)$ be an oriented  Riemannian manifold of dimension
$m.$

\bigskip

By Exercise \ref{dualmetric},
where we set
$(V, \langle, \rangle )= (T^{*}_{M,q}, g^{*}_{q}),$
the Riemannian metric $g$ on  $M$ 
defines a smoothly  varying 
inner product on the exterior algebra bundle
$\Lambda( T^{*}_{M} ).$

\medskip
\n
By $\S$\ref{tsi}, the orientation on $M$ gives rise
to the  $\star$ operator
on the differential forms on $M:$
\begin{equation}
\label{eq3}
\fbox{$
\star\, :\, E^{p}(M) \lorw E^{m-p}(M).$}
\end{equation}
In fact, the star operator is defined point-wise,
using the metric and the orientation,
on the exterior algebras $\Lambda( T^{*}_{M,q}  )$ and 
it extends to differential forms.

\bigskip
\begin{rmk}
    {\rm 
    The example  $M = \real$ with the standard orientation and
    the Euclidean metric shows that
    $\star$ and $d$ do not commute. In particular, $\star$
    does not preserve closed forms. See Exercise \ref{lsc}
    and Lemma \ref{dd*closed}.
    } 
    \end{rmk}

\bigskip
\begin{defi}
    \label{defofrlangle}
    {\rm
Define  an inner product 
on the space of compactly
supported $p-$forms  on $M$ by setting  
\begin{equation}
\label{ipl2}
\fbox{$
\langle \langle \, u,v \,\rangle \rangle  : = \int_{M}{ \langle u,v\rangle \,
dV }
$} \, = \, \int_{M}{ u \wedge \star \, v}.
\end{equation}
For the last equality, see Definition \ref{doso}.
}
\end{defi}

\bigskip
Once we have a metric, we may start talking about formal adjoints.
For a more thorough discussion of this notion, 
see \ci{dem1}, $\S 2.6.$
\begin{defi}
    \label{defofadj}
    {\rm ({\bf Formal adjoints})
Let $T: E^{p}(M) \lorw E^{p'}(M)$ be a linear map.
We say that a linear map 
$$
T^{\star}: E^{p'}(M)\lorw
E^{p}(M)
$$
is the {\em formal
adjoint} to $T$ (with respect to  the metric) if,
for  every  compactly supported $u \in E^{p}(M)$ and
$v \in E^{p'}(M):$
$$
\fbox{$
\langle \langle \,  Tu,v \, \rangle \rangle \, = \,  \langle \langle \, u , 
T^{\star} v \, \rangle \rangle .
$}$$

    }
    \end{defi}

\medskip
\begin{defi}
\label{d*}
{\rm
({\bf Definition of  $d^{\star}$})
Define
$ d^{\star} : E^{p}(M)  \lorw E^{p-1}(M)$
as
$$
\fbox{$ 
d^{\star} \, := \;(-1)^{m(p+1)+1}\, \star \,  d\, \star.
$}
$$
}
\end{defi}

\bigskip
\begin{pr}
\label{dadjd*}
The operator  $d^{\star}$ is the formal adjoint of $d.$
\end{pr}
{\em Proof.}
We need to show that if $u_{p-1} \in E^{p-1}(M)$ and $v_{p} \in E^{p}(M)$
have compact support, then 
$$
\langle \langle  \,du_{p-1}, v_{p}\,\rangle \rangle = \langle \langle \, 
u_{p-1}, 
d^{\star}v_{p}\, 
\rangle \rangle .
$$
We have
$$
\langle \langle \, du_{p-1},v_{p}\rangle \rangle =
\int_{M}{ \langle \, du_{p-1},v_{p} \, \rangle \, dV } = 
\int_{M}{ du_{p-1} \wedge \star\, v_{p} }=
$$
$$
=\int_{M}{ d( u_{p-1} \wedge \star\, v_{p} ) - (-1)^{p-1} u_{p-1}
\wedge d \star\, v_{p}}
=\int_{M}{  (-1)^{p} u_{p-1} \wedge d \star \, v_{p}  }
$$
$$
=
\int_{M}{  (-1)^{p}  u_{p-1} \wedge (-1)^{(m-p +1)(p-1)}  \star\,  \star
\, d \, \star \, 
v  }=
$$
$$
= \int_{M}{    u \wedge  \star \, ( (-1)^{m(p+1)+1}   \star\,  d\,  
\star\,  v ) }
= \langle \langle \,   u, (-1)^{m(p+1)+1}   \star\, d\, \star\,
v \, \rangle \rangle
\, dV
$$
where  the fourth equality follows from Stokes' Theorem \ref{st},
the fifth from (\ref{stst}) and the last one from 
(\ref{ipl2}) and simple $mod \,2$ congruences.
\blacksquare

\subsection{The Laplace-Beltrami operator of an oriented
Riemannian manifold}
\label{lboorm}
\begin{defi}
\label{deflbo}
{\rm ({\bf The Laplacian})
The {\em Laplace-Beltrami operator,} or Laplacian, is defined 
as 
$\Delta:
E^{p}(M) \lorw E^{p}(M)$
$$
\fbox{ $
 \Delta \, := \;  d^{*}d \, + \, d\, d^{*}.
$}
$$
}
\end{defi}

\bigskip
We have 
$$
\Delta \, = \, (-1)^{m(p+1) +1} d \star d \star \, + \, 
(-1)^{ mp+1    } \star d \star d.
$$

\bigskip
While $\star$ is defined point-wise using the metric,
$d^{\star}$ and $\Delta$ are defined locally (using $d$) and depend
on the metric.

\bigskip
\begin{exe}
    \label{lsc}
    {\rm 
    Show that $\star\, \Delta \,= \,\Delta\, \!  \star.$
    In particular, a form $u$ is harmonic iff $\star \, u$ is harmonic.
    }
    \end{exe}

\bigskip
\n
The Laplacian is  self-adjoint, i.e. it coincides with
its adjoint and is an {\em  elliptic} second order linear partial 
differential operator. See  \ci{dem1}, $\S 4.14$ and \ci{wa}, 6.28, 
6.35, Exercise 9, 16, 18 (the wave equation is not elliptic
and there is no regularity theorem for it), 21 (elliptic
operators for vector bundles).

\begin{exe}
    \label{itslap}
    {\rm
Verify that
if $(M,g)$ is the Euclidean space
$\real^{n}$ with the standard metric and  
$
u = \sum_{|I|=p}{ u_{I} \, dx_{I}},
$
then
$$
\fbox{$
\Delta (u)\,  =\,  - \sum_{I}{ \sum_{j}{  \frac{\partial^{2}u_{I}}{\partial
x_{j}^{2}} } \, dx_{I}}.  $}
$$
See \ci{dem2}, $\S VI.3.12$ and \ci{gh}, page 83.
}
\end{exe}

\subsection{Harmonic forms and the Hodge 
Isomorphism Theorem}
\label{hfhit}
Let $(M,g)$ be a \underline{compact} oriented  Riemannian
manifold.
\bigskip

\begin{defi}
\label{dhf}
{\rm ({\bf Harmonic forms})
Define the space of  {\em real harmonic $p-$forms}
as 
$$
\fbox{$
{\cal H}^{p}(M, \real): = \ke{\, ( \,  \Delta: E^{p}(M) \lorw E^{p}(M)} \, ).
$}
$$
}
\end{defi}

\bigskip
The space of harmonic forms depends on the metric. A form may be 
harmonic with respect to one metric and fail to be harmonic
with respect to another metric.

\bigskip

\begin{lm}
\label{dd*closed}
A $p-$form $u \in E^{p}(M)$ is harmonic
if and only if $du=0$ and $d^{\star}u=0.$
\end{lm}
{\em Proof.}  It follows  from the identity
$$
\langle \langle \Delta u, u \rangle \rangle  = \langle \langle
dd^{\star} u + d^{\star}d u, u\rangle \rangle =
\langle \langle
d^{\star}u ,  d^{\star}u \rangle \rangle +
\langle \langle
du ,  du \rangle \rangle =
||d^{\star}u||^{2} + ||du||^{2}.
$$
\blacksquare

\bigskip
Note that if $M$ is {\em not} compact, then a form
which is $d-$closed  and $d^{\star}-$closed is harmonic by the
definition
of $\Delta.$ However, the converse is not true, e.g.
the function $x$ on $M= \real.$ The argument 
above breaks down in handling not necessarily convergent
integrals.

\bigskip

We now come to the main result concerning the Hodge theory of
a compact  oriented Riemannian compact manifold. The proof
is based on the theory of elliptic differential operators.
See  \ci{wa}, \ci{gh} for self-contained proofs.

\begin{tm}
\label{hdts}
({\bf The Hodge Orthogonal Decomposition Theorem})
Let (M,g) be a \underline{compact} oriented  Riemannian manifold. 

\n
Then 
$$
\fbox{$
\dim_{\real}{ \, {\cal H}^p(M, \real) }  \,  < \, \infty
$}
$$
and we have a direct sum decomposition
into $\langle \langle \, , \, \rangle \rangle-$orthogonal subspaces
$$
\fbox{$
E^{p}(M)\,  =\,  {\cal H}^{p}(M, \real) \, 
\stackrel{\perp}\oplus \,  d \,(E^{p-1}(M))
\, \stackrel{\perp}\oplus \,  d^{\star} ( E^{p+1}(M)).
$}
$$
\end{tm}

\bigskip

\begin{exe}
\label{umv}
{\rm 
Show that if $u$ and $v$ are $m-$forms on a compact oriented 
manifold $M$ of dimension $m$ such that
$\int_{M}{u} = \int_{M}{v},$ then $(u-v)$ is exact.
}
\end{exe}

\bigskip

\begin{exe}
\label{ortdc}
{\rm 
Show that the orthogonality of the direct sum decomposition
is an easy consequence of Lemma \ref{dd*closed}.
}
\end{exe}

\begin{exe}
    \label{adjkerimort}
    {\rm 
    Show that if  $T$ and $T^{*}$ are adjoint, then
    $\ke{\,T} = (\im{\,T^{*}})^{\perp}.$
    }\end{exe}

\begin{cor}
 \label{hit}
 ({\bf The Hodge Isomorphism Theorem})
 Let $(M,g)$ be a compact oriented Riemannian manifold.
There is an isomorphism depending only on the metric $g:$ 
$$
\fbox{$
{\cal H}^{p}(M,\real ) \simeq  H^{p}_{dR}(M,\real).$}
$$
In particular, $\dim_{\real}{H^{p}(M, \real)} < \infty.$
\end{cor}
{\em Proof.} 
By Lemma \ref{dd*closed}, $d$ and $d^{*}$ are adjoint. 
By Exercise \ref{adjkerimort},
we have the equality 
$\ke{\,d} = \im{\,d^{\star} }^{\perp},$ the right-hand-side
of which is ${\cal H}^{p}(M, \real ) \stackrel{\perp}\oplus d E^{p-1}(M)$
by
Theorem  \ref{hdts}.
The corollary follows.
\blacksquare

\bigskip
The following is one way to think about Corollary \ref{hit}.
The inclusion $\im{\,d} \subseteq \ke{\,d}$  and its quotient
$H^{p}_{dR}(M, \real)$ are  independent of the metric and canonically
attached to the differentiable structure of the manifold $M.$
The metric allows to distinguish a vector subspace, ${\cal H}^{p}(M,\real),$ 
complementary to
$\im{\,d}$ in $\ke{\,d}$ and hence isomorphic
to $H^{p}_{dR}(M,\real).$

\bigskip
The following heuristic argument, which is taken from \ci{gh}, page 80, where
the case of Dolbeault cohomology is treated, suggests why
one may think that the results of Hodge theory should 
hold. 

\n
Let $\alpha = [u +dv] \in H^{p}_{dR}(M,\real),$ 
where $u \in E^{p}(M)$ is a closed form representing
$\alpha$ 
and $v\in E^{p-1}(M).$ 

\n
We have that
{\em $||u||$ has minimal norm among all representatives
of $\alpha$  iff $d^{*}u=0$}, i.e., in view of
Lemma \ref{dd*closed}, iff $u$ is harmonic. 

\n
This means 
that the choice of an  orientation and of a metric
allows to distinguish a representative of any cohomology class
in such a way that the norm is minimized and the representatives
form a vector subspace of the closed forms.

\medskip
\n
The condition $d^{\star}u=0$ is sufficient: 
$$
||u + dv||^{2} \; = \;
||u||^{2}+ ||dv||^{2} + 2\langle\langle u, dv \rangle \rangle
\;\geq \; ||u||^{2} + 2\, \langle\langle d^{\star}u, v \rangle \rangle
\;=\; ||u||^{2}.
$$
The condition $d^{\star}u=0$ is necessary:  if $u$ has smallest norm,
then
$$
0= \frac{d}{dt}\,  (\, ||u + t \, dv||^{2}\, ) (0) = 2\, \langle\langle
u, dv \rangle \rangle = 2 \, \langle \langle d^{\star} u, v \rangle
\rangle \quad \forall v \, \in \, E^{p-1}(M),
$$
so that $d^{\star}u=0.$

\bigskip
Let us emphasize again that the space ${\cal H}^{p}(M, \real)$ 
and the identification with de Rham cohomology depend on 
the metric. Once the metric has been fixed, the content of  
Corollary \ref{hit}
is that one has a unique harmonic representative
for every  de Rham  cohomology class.

\bigskip
We are now in the position to give a proof of  the Poincar\'e Duality
isomorphism on an oriented compact manifold. 

\begin{tm}
\label{pd}
({\bf Poincar\'e Duality})
Let $M$ be a compact oriented smooth manifold.
The pairing 
$$
H^{p}_{dR}(M, \real ) \times H^{m-p}_{dR}(M, \real) \lorw \real , 
\qquad (u,v) \lorw 
\int_{M}{u\wedge v}
$$
is non-degenerate, i.e. it 
induces an isomorphism (the Poincar\'e Duality Isomorphism)
$$
\fbox{$
H^{m-p}_{dR}(M, \real ) \simeq H^{p}_{dR}(M, \real)^{\vee}.$}
$$
\end{tm}
{\em Proof.}
Note that the pairing, which is a priori  defined 
at the level of forms, descends to de Rham cohomology
by virtue of Stokes' Theorem (applied to $d(u\wedge v')$,
$d( u' \wedge v)$ and to $du' \wedge dv'$ which is exact (check!)), i.e.
if $u$ and $v$ are closed, then 
$$
\int_{M}{ (u+du')\wedge (v+dv') } = \int_{M}{ u\wedge v}.
$$
Fix a Riemannian metric on $M.$
Using the Hodge Isomorphism Theorem \ref{hit},
we may use harmonic representatives. 
The key point is that $u$ is harmonic if and only if
$\star \,u$ is harmonic; see
Exercise \ref{lsc}.
We have the equality
$$
\int_{M}{u \wedge \star \, u}= \int_{M}{ \langle u,u\rangle \, dV} = ||u||^{2}.
$$
This shows that given any non-zero harmonic form
$u \in H^{p}_{dR}(M, \real),$ there is an harmonic form, $\star \,u
\in H^{m-p}_{dR}(M, \real),$ such that it pairs
non-trivially with
$u.$ This shows that the bilinear pairing
induces  injections
$H^{p}_{dR}(M) \to H^{m-p}_{dR}(M)^{\vee}$ and, symmetrically,
$H^{m-p}_{dR}(M) \to H^{p}(M)^{\vee}.$
The result follows from the fact that injections
of finite dimensional vector spaces of equal dimension are isomorphisms.
\blacksquare

\bigskip

\begin{rmk}
\label{nciso}
{\rm
The proof given above shows also that
the $\star$ operator induces isomorphisms
$$
\star \, : \; {\cal H}^p(M, \real ) \, \simeq \, {\cal H}^{m-p}(M, \real)
$$
for any compact, smooth, oriented Riemannian manifold $M.$
Via the Hodge Isomorphism Theorem \ref{hit},
these induce isomorphisms
$$
H^p_{dR}(M, \real) \simeq H^{m-p}_{dR}(M, \real).
$$
These latter  isomorphisms are not canonical.
On the other hand, 
the Poincar\'e Duality
isomorphisms depend only on the orientation.
Changing the orientation, changes the sign of integrals
and hence of the Poincar\'e Duality  isomorphism.
}
\end{rmk}

\bigskip
\begin{exe}
\label{dualisit}
{\rm
Let $M$ be a smooth compact orientable manifold of dimension $m,$
$r $ be an  even integer,
$v \in H^{r}(M,\rat),$ and define
$$
L : H^{\bullet}(M, \rat) \lorw H^{\bullet + r} (M, \rat) ,\qquad 
u \lorw u\cup v.
$$
Show that,  in a suitable sense, $L$ coincides with the dual 
map $ L^{\vee}.$ 

\n 
State and prove the analogous statement for any $r.$
}
\end{exe}

\newpage
\section{Lecture 3: Complex manifolds}
\label{l3}
We discuss various notions of conjugation on complex vector spaces,
tangent and cotangent bundles on  a complex manifold, the standard 
orientation of a complex manifold, the quasi complex structure, 
$(p,q)-$forms, $d'$ and $d'',$ Dolbeault and Bott-Chern
cohomology.

\subsection{Conjugations}
\label{conjugations}
Let us recall the following distinct notions of conjugation.

\begin{itemize}

\item 
There is of course the usual conjugation in $\comp:$
$\gamma \lorw \ov{\gamma}.$

\item  
Let $V$ be a real vector space and $V_{\comp}:= V \otimes_{\real}
\comp$ be its complexification. 
There is the natural $\real-$linear isomorphism
given by  what we call {\em c-conjugation}
(i.e. coming from the complexification)
$$
c \, :\,  V_{\comp} \lorw V_{\comp}, \qquad u\otimes \gamma \lorw u \otimes
\overline{\gamma}.
$$
We may denote $c(v)$ by ${\ov{v}}^{c}$

\item
Let $V$ and $V'$ be real vector spaces, 
$P: V_{\comp} \lorw V'_{\comp}$ 
be complex linear.  Define the {\em operational conjugation}
the complex linear map
\begin{equation}
    {\ov{P}}^{o}: V_{\comp} \lorw V'_{\comp}, 
    \qquad
    {\ov{P}}^{o} (v) \, = \, \ov{  P( {\ov{v}}^{c} )}.
    \end{equation}
    Note that the assignment $P \to {\ov{P}}^{o}$
    defines  a real linear isomorphism
    $$
    Hom_{\comp}(V_{\comp} , V'_{\comp} ) 
    \simeq_{\real} Hom_{\comp}(V_{\comp} , V'_{\comp} ) 
    $$
    which is  not complex linear. 
    As the reader should check, $\ov{P}^{c} = \ov{P}^{o}$
    in the space with conjugation
    $Hom_{\comp}(V_{\comp} , V'_{\comp} ) = Hom_{\real}(V, 
    V' ) \otimes_{\real} \comp.$

\item
Let $W$ be a complex vector space. The 
{\em conjugate $\overline{W}$ of $W$} 
is the complex vector space such that
 $W = \overline{W}$ as  real vector 
spaces, but such that 
the scalar multiplication in $\overline{W}$ is defined as
$$
\gamma  \cdot w : = \overline{\gamma} w
$$
where, on the right-hand-side
 we are using  the given complex scalar multiplication
on $W.$
The identity map $Id_{W}: W\lorw {\overline{W}}$ is a real linear 
isomorphism, but is not complex linear.

\item
Let $f: S \lorw \comp$ be a map of sets. Define
${\ov{f}}^{r} : S \lorw \comp$ as ${\ov{f}}^{r}(s) = \ov{ f(s) }.$

\end{itemize}

\bigskip
As the reader should check, we have, for every 
$P \in Hom_{\comp}(V_{\comp} , V'_{\comp} )$
\begin{equation}
\label{ahfc}
\fbox{$
{\ov{P}}^{o}(v)  \, = \, {\ov{P}}^{c} (v)  \, = \, 
{ \ov{P}}^{r} ( {\ov{v}}^{c} ) \, = \, 
\ov{ P( {\ov{v}}^{c} )}.$}
\end{equation}

\subsection{Tangent bundles on a  complex manifold}
\label{tbocm}
Let $X$ be a complex manifold of dimension $n,$ 
$x \in X$ and  
$$
(U\, ; \, z_{1}= x_{1} + i\,y_{1}, \ldots , z_{n}=
x_{n}+ i\,y_{n} )
$$
be a holomorphic chart for $X$ around $x.$

\n
We have the following notions of ``tangent bundles'' for $X.$
Depending on the context, the symbols
$\partial_{x_{j}},$ $dx_{j},$ etc.  will denote either vectors in 
the fibers at $x$  of the corresponding   bundles,
 or the corresponding  local frames for those bundles.

 \begin{itemize}

     \item
     $T_{X}(\real),$
the {\em real tangent bundle}. The fiber $T_{X,x}(\real)$  
has real rank $2n$ and it
is the real span 
$$
\fbox{$ \real \langle \,  \partial_{x_{1}}, \ldots , \partial_{x_{n}},
 \partial_{y_{1}}, \ldots , \partial_{y_{n}} \, \rangle .$}
 $$

 \item $T_{X}(\comp):= T_{X}(\real)\otimes_{\real} \comp,$
the {\em complex tangent bundle}. The fiber $T_{X,x}(\comp):=
T_{X,x}(\real)\otimes_{\real} \comp$  
has complex rank $2n$ and 
is the complex span 
$$
\fbox{$ \comp \langle \,  \partial_{x_{1}}, \ldots , \partial_{x_{n}},
 \partial_{y_{1}}, \ldots , \partial_{y_{n}} \,  \rangle .$}
 $$
Note that $T_{X} (\real) \subseteq T_{X}(\comp)$
via the natural real linear map $v \to v \otimes 1.$
Set 
\begin{equation}
    \label{natb}
    \fbox{$
\partial_{z_{j}}:= \frac{1}{2}( \partial_{x_{j}} - i\, \partial_{y_{j}}),
$}\qquad
\qquad 
\fbox{$\partial_{\overline{z}_{k}} :=  \frac{1}{2}(
\partial_{x_{j}} + i\, \partial_{y_{j}}).$}
\end{equation}
Clearly, we have
\begin{equation}
    \label{natbx}
    \fbox{$
\partial_{x_{j}} =\partial_{z_{j}} +  \partial_{\ov{z}_{j} },
$}\qquad
\qquad 
\fbox{$\partial_{\overline{y}_{j}} =  i\, (
\partial_{z_{j}} -  \partial_{\ov{z}_{j}}).$}
\end{equation}

We have $T_{X,x}(\comp)= \comp \langle \partial_{z_{1}}, \ldots ,
\partial_{z_{n}}, \partial_{\overline{z}_{1}}, \ldots,
\partial_{\overline{z}_{n}} \rangle .$

\end{itemize}

\bigskip
In general, a smooth change of 
coordinates $(x,y) \lorw (x',y')$
 does not leave invariant
 the two subspaces
$\real \langle \{\partial_{x_{j}} \}\rangle $ and 
$\real \langle \{\partial_{y_{j}} \}\rangle $
of $T_{X,x}(\real).$

 \bigskip
\n
    A local linear change of coordinates
$z'_j = \sum_{k}{A_{jk}} z_k$
(this is the key case to check when dealing
with this  kind of questions)
produces a change of basis
in $T_{X,x}(\comp):$
$$
\partial_{z'_j} =\sum_j{ (A^{-1})_{kj} \partial_{z_k} }, \qquad
\partial_{\ov{z}'_j} =\sum_j{ \ov{(A^{-1})_{kj}} \partial_{\ov{z}_k} }.
$$
The simple yet important
consequence of this observation 
 is that  a holomorphic change of coordinates
$z \lorw z'$ fixes  the two subspaces  
 $\comp \langle \{ \partial_{z_{j}} \} \rangle \subseteq T_{X,x}(\comp)$ 
 and
$\comp \langle \{ \partial_{\overline{z}_{j}} \} \rangle
\subseteq T_{X,x}(\comp).$  It follows that we may define
the following complex smooth vector bundles:

\begin{itemize}
    
\item $T'_{X}$ the {\em holomorphic tangent bundle}. The fiber
$$
\fbox{$
T'_{X,x}\, =\,  \comp \langle\, \partial_{z_{1}}, \ldots, 
\partial_{z_{n}}\,  
\rangle $} 
$$ 
has complex rank $n.$ $T'_{X}$ is a holomorphic vector bundle. 

\item $T''_{X}$ the {\em anti-holomorphic tangent bundle}. The fiber
$$
\fbox{$
T''_{X,x}= \comp \langle \, \partial_{\overline{z}_{1}}, \ldots, \partial_{
\overline{z}_{n}} \,   \rangle  $}
$$ 
has complex rank $n.$ It is an anti-holomorphic vector bundle,
i.e.  it admits  transition functions which  are
 conjugate-holomorphic. 
 
\end{itemize}
 
\bigskip
We have  a canonical injection and
a canonical  internal direct sum decomposition into complex sub-bundles
\begin{equation}
\label{ttt}
\fbox{$ T_{X } (\real) \, \subseteq  \,
T_{X}(\comp) = T'_{X} \oplus T''_{X}.$}
\end{equation}

\bigskip
Composing the injection with the projections we get
canonical real isomorphisms
\begin{equation}
\label{tttiso}
\fbox{$ T'_{X} \, \simeq_{\real} \, T_{X}(\real) \, \simeq_{\real} \,
 T''_{X}, $}
 \end{equation}
 \begin{equation}
\label{tttisox}
\fbox{$ \partial_{z_{j}} \leftarrow \partial_{x_{j}} \to 
\partial_{\ov{z}_{j}}, $} \qquad \qquad  
\fbox{$i\, \partial_{z_{j}} \leftarrow \partial_{ y_{j} } \to 
 -i\, \partial_{\ov{z}_{j}}$}
 \end{equation}

\bigskip
The conjugation map $c: T_{X}(\comp) \lorw T_{X}(\comp),$
$c(v) = { \ov{v} }^{c},$ is a real linear
isomorphism
which is not complex linear and it
induces a real linear isomorphism 
\begin{equation}
\label{rciso}
\fbox{$ c\, : \, T'_{X}\, \simeq_{\real} \,T''_{X},$} \quad 
\fbox{$ c(\partial_{z_{j}} )\,  = \,
\partial_{\overline{z}_{j}},$}  \qquad 
\fbox{$ c\, : \, T''_{X}\, \simeq _{\real}\, T'_{X}, $} \quad 
\fbox{$ c(\partial_{\overline{z}_{j}} ) \,  = \, 
\partial_{z_{j}},$} 
\end{equation} 
and complex linear isomorphisms 
\begin{equation}
\label{cciso}
\fbox{$ c\, : \, T'_{X}\, \simeq _{\comp} \, \overline{T''_X}, $} \qquad 
\qquad
\fbox{$ c\, :\,   T''_{X}\,  \simeq_{\comp}\,
\overline{T'_X}.$}
\end{equation} 

\bigskip

\begin{exe}
\label{proj}
{\rm ({\bf Standard vector bundles on $\pn{n}$})
Verify the assertions that follow. 

\n
Let $\pn{n}$ be the complex projective space of dimension $n.$
It is the set of  equivalences classes of 
$(n+1)-$tuples of not-all-zero complex numbers
$ [x]= [x_{0}: \ldots : x_{n}],$ where we identify
two such tuples if they differ by a non-zero multiplicative constant.
One endows $\pn{n}$ with the quotient
topology stemming from the identification
$\pn{n}=(\comp^{n+1}\setminus \{0 \})/ \comp^{*}.$
 $\pn{n}$ is a  compact topological space.
Consider the following open covering consisting 
of $n+1$ open subsets of $\pn{n}:$
$U^{j}= \{\, [x] \, |\;
x_{j} \neq 0\, \},$ $ 0 \leq j \leq n.$ We have homeomorphisms
\[
\tau_{j}: \,  U^{j} \lorw \comp^{n}, \qquad 
[ x_{0}: \ldots : x_{n}] \lorw \left( z_{0}^{j}:=\frac{x_{0}}{x_{j}} ,
\;
\stackrel{\hat{j}}\ldots  \; ,\, z_{n}^{j}: = 
\frac{x_{n}}{x_{j}} \right).
\]
There are  biholomorphisms
$$
\tau_{k} \circ \tau_{j}^{-1}: \comp^{n}\setminus \{ z_{k}^{j} = 0 \}
\lorw
\comp^{n}\setminus \{ z_{j}^{k} = 0 \}
$$
so that $\pn{n}$ is a complex manifold of dimension $n.$ 
The transition functions are
\[ \left\{ \begin{array}{cclc}
z_{l}^{k} &=& (z^{j}_{k})^{-1} \, z_{l}^{j} & l \neq j, \\
z_{j}^{k} &=& (z_{k}^{j})^{-1}. &
\end{array}
\right.
\]
Let $a \in \zed.$ Consider the nowhere vanishing holomorphic functions
$(z_{k}^{j})^{-a} $ on $U^{j} \cap U^{k}.$
They define a holomorphic line bundle $L_{a}$ on $\pn{n}.$
For $a \geq 0$ 
the holomorphic sections of $L_{a}$ are in natural bijection
with the degree $a$  polynomials in $\comp[x_{0}, \ldots, x_{n}]$
so that the zero locus of a section coincides with the zero set
of the corresponding polynomial (counting multiplicities). 

\n
The line bundle 
$L_{1}$ is called the {\em hyperplane bundle}.

\n The line bundles
$L_{a}$, $a <0,$ have no holomorphic sections.

\n
The line bundle $L_{-1}$ is the tautological line bundle, i.e.
the one that has as fiber over a point  $[x] \in \pn{n}$
the line $\lambda (x_{0}, \ldots , x_{n}) \subseteq \comp^{n+1}.$

\n
One has
$L_{a+b} \simeq L_{a}\otimes L_{b},$ and
$L_{a} \simeq L_{b}$ iff $a=b$, even as topological complex vector 
bundles.

\n
These are the only (topological complex) holomorphic line bundles on
$\pn{n}.$

\n
Use the atlas given above above to compute the transition functions of
$T'_{\pn{n}}$ and check that
the anti-canonical line bundle
$K_{\pn{n}}^{*} := \det{  T'_{\pn{n}} } \simeq L_{n+1}.$

\n
There is a non-splitting 
exact sequence (called the Euler sequence) of holomorphic vector bundles
$$
0 \lorw L_{-1} \lorw L_{0}^{\oplus n+1} \lorw T'_{\pn{n}}\otimes L_{-1} 
\lorw 0.
$$
The total Chern class of $T'_{\pn{n}}$ is $c (T'_{\pn{n}})=
(1 + c_1(L_1))^{n+1}.$ See \ci{b-t} for the definition of Chern
classes and their basic properties.
}
\end{exe}

\subsection{Cotangent bundles on complex manifolds}
\label{cbocm}
Denote by  $\{ dx_{1}, \ldots , dx_{n}, dy_{1}, \ldots,
dy_{n} \}$ the basis dual to $\{ \partial_{x_{1}}, \ldots , \partial_{x_{n}},
 \partial_{y_{1}}, \ldots , \partial_{y_{n}} \},$
 by
$\{ dz_{1}, \ldots , dz_{n} \} $ the basis dual to $\{ \partial_{z_{1}}, 
 \ldots , \partial_{z_{n}} \},$ and
 by
 $\{ d\overline{z}_{1}, \ldots , d\overline{z}_{n} \} $ 
 the basis dual to $\{ \partial_{\overline{z}_{1}}, 
 \ldots , \partial_{\overline{z}_{n}} \}.$

 \begin{exe}
     \label{sugus}
     {\rm 
Verify the following identities.
\begin{equation}
\label{dzdzbar}
    \fbox{$ dz_{j}\,  = \, dx_{j} + i\,dy_{j},$} \qquad \qquad 
\fbox{$ d \overline{z}_{j}\,  =\,  dx_{j}- i\,  dy_{j}.$}
\end{equation}
\begin{equation}
\label{dzdzbarx}
\fbox{$
dx_{j} \, = \, \frac{1}{2} \, ( dz_{j} \, + \, d\ov{z}_{j} ),
$}
\qquad \qquad 
\fbox{$ d y_{j}\,  =\,  
\frac{1}{2i} \, ( dz_{j}\, -   d\ov{z}_{j}) .$}
\end{equation}
 }
 \end{exe}

 \bigskip
We have the following vector bundles on $X$.

 \begin{itemize}
     
\item $T^{*}_{X}(\real),$ the {\em real cotangent bundle}, with fiber
$$
\fbox{$
T^{*}_{X,x}(\real) \, =\,  \real \langle \, 
dx_{1}, \ldots , dx_{n}, dy_{1}, \ldots,
dy_{n} \, \rangle . $}
$$

\item 
$T^{*}_{X}(\comp):= (T_{X}(\comp))^{*} = T_{X}^{*}(\real)
\otimes \comp,$ the {\em complex cotangent bundle}, with fiber
$$
\fbox{$
T^{*}_{X,x}(\comp) \,=\, \comp \langle\,  
dx_{1}, \ldots , dx_{n}, dy_{1}, \ldots,
dy_{n} \, \rangle . $}
$$

\item  $T'^{*}_{X}$ the {\em holomorphic cotangent bundle}, with fiber
$$
\fbox{$
T'^{*}_{X,x} \, = \,  \comp \langle \,  dz_{1}, \ldots  , dz_{n} \,
\rangle . $}
$$ 
It is a holomorphic
vector bundle, dual to $T'_{X}.$

\item $T''^{*}_{X}$ the {\em anti-holomorphic cotangent bundle}, with fiber
$$
\fbox{$
T''^{*}_{X,x} \, = \, 
\comp \langle \,  d\overline{z}_{1}, \ldots  , d\overline{z}_{n} 
\,
\rangle . $}
$$
It is an anti-holomorphic
vector bundle, dual to $T''_{X}.$
\end{itemize}

\bigskip
We have a canonical injection
and a canonical internal direct sum decomposition 
\begin{equation}
\label{uup}
\fbox{$
T^{*}_{X} (\real) 
\, \subseteq  \,
T^{*}_{X} (\comp) \, =\,  T'^{*}_{X} \, \oplus \,  T''^{*}_{X}.$}
\end{equation}

\bigskip
Composing the injection with the projections we get
canonical real isomorphisms
\begin{equation}
\label{tttisod}
\fbox{$ T'^{*}_{X} \, \simeq_{\real} \, T^{*}_{X}(\real) \, \simeq_{\real} \,
 T''^{*}_{X}, $}
 \end{equation}
 \begin{equation}
\label{tttisoxd}
\fbox{$ \frac{1}{2} \, d z_{j}
\leftarrow dx_{j}
\to  \frac{1}{2} \,  d \ov{z}_{j} , $} 
\qquad \qquad  
\fbox{$
\frac{1}{2i} \, d z_{j}
\leftarrow dy_{j}
\to  -\frac{1}{2i} \,  d \ov{z}_{j} ,
$}
 \end{equation}

 \bigskip

The conjugation map $c: T^{*}_{X}(\comp) \lorw T^{*}_{X}(\comp)$ 
is a real linear
isomorphism,
which is not complex linear. It  
induces a real linear isomorphism 
\begin{equation}
\label{strciso}
\fbox{$ c\, : \,  T'^{*}_{X}\, \simeq _{\real}\, T''^{*}_{X}, $}
\quad 
\fbox{$ c(d{z_{j}} ) \, =\, 
d\overline{z}_{j}, $} \qquad 
\fbox{$ c\, :\,  T''^{*}_{X}\, \simeq _{\real}\, T'^{*}_{X},$} \quad 
\fbox{$ c(d\overline{z}_{j} )\,  =\, 
d{z_{j}},$} 
\end{equation} 
and complex linear isomorphisms 
\begin{equation}
\label{stcciso}
\fbox{$ c\, :\, T'^{*}_{X}\, \simeq _{\comp}\, \overline{T''^{*}_{X}},$} 
\qquad  \qquad 
\fbox{$ c\, : \,   T''^{*}_{X} \, \simeq_{\comp}\, 
\overline{T'^{*}_{X}}.$}
\end{equation}

\bigskip
Let $f(x_1,y_1, \ldots, x_n, y_n) = u (x_1,y_1, \ldots, x_n, y_n) 
+ i\, v(x_1,y_1, \ldots, x_n, y_n)$
be a smooth complex-valued function in a neighborhood of $x.$
We have 
\medskip

\begin{equation}
    \label{diff}
\fbox{$ df \,
= \,
du + i\, dv \,  = \,
\sum_j{ \frac{\partial f}{\partial z_j} dz_j } +
\sum_j{ \frac{\partial f}{\partial \ov{z}_j }  d\ov{z}_j } . $}
\end{equation}

\bigskip
\begin{exe}
    \label{exeverdiff}
    {\rm 
Verify the equality (\ref{diff}).    
    }
    \end{exe}

\subsection{The standard orientation of a complex manifold}
\label{tcooacm}
\begin{pr}
\label{cmao}
A  complex manifold $X$ admits   a canonical orientation.
\end{pr}
{Proof.} By Exercise \ref{saor} it is enough to exhibit
an atlas where the transition functions have positive
determinant Jacobian.
The holomorphic atlas is one such. In fact,   the
determinant Jacobian
with respect to the change of 
coordinates $(x_{1}, y_{1}, \ldots ,
x_{n}, y_{n}) \lorw 
(x'_{1}, y'_{1}, \ldots ,
x'_{n}, y'_{n})$  is positive
being the square of the absolute value of the  
determinant Jacobian of the holomorphic 
change of coordinates $(z_{1}, \ldots , z_{n})
\lorw (z'_{1}, \ldots , z'_{n}).$
\blacksquare

\bigskip
Here is a more direct proof of Proposition \ref{cmao}.

\medskip
\n
Let $(U;z)$ be a chart in the holomorphic atlas of $X,$
$x \in U.$

\medskip
\n
The  real $2n-$form
\begin{equation}
\label{deforincn}
\fbox{$  { o}_{U } \, := \, 
dx_1 \wedge dy_1 \wedge \ldots \wedge dx_n \wedge dy_n \, $}
\end{equation}
is nowhere vanishing on $U$ and 
defines an orientation for $U.$

\n
One checks that 
$$
{ o}_{U}\,  = \,  (i/2)^{n} \, dz_{1}\wedge d\overline{z}_{1}
\wedge \ldots \wedge
dz_{n}\wedge d\overline{z}_{n}
$$
 and that 
$$
{ o}_{U} = (i/2)^{n}\, (-1)^{  \frac{(n-1)n}{2} } \,
dz_{1} \wedge \ldots  \wedge dz_{n} \wedge 
d \overline{z}_{1} \wedge  \ldots \wedge d\overline{z}_{n}.
$$
Let $(U';z')$ be another chart around  $x.$

\n
We have
$$
J(z(x)) := \det{ || \frac{\partial z'_{j} }{\partial 
z_{k} } (z(x)) ||} > 0.
$$

\n
Since
$$
dz'_{1} \wedge \ldots \wedge dz'_{n} \wedge 
d \overline{z}'_{1} \wedge \ldots \wedge d \overline{z}'_{n}=
|J(z(x))|^{2} \,  dz_{1} \wedge \ldots \wedge dz_{n} \wedge 
d \overline{z}_{1} \wedge \ldots \wedge d \overline{z}_{n},
$$
we have that
$$
{ o}_{U'} \, = \,  |J(z(x))|^{2} \, { o}_{U}.
$$
It follows that if we use a covering
of $X$ by means of holomorphic charts, we can glue
the forms $o_{U}$ using a partition of unity subordinate to the 
covering and obtain a nowhere vanishing real  $2n-$form
$o$ which orients $X$ independently of the covering chosen
within the holomorphic atlas.

\bigskip
 This is the so-called {\em standard orientation}
of $X;$
at every point $x \in X$ it is determined  by the vector (\ref{deforincn}).

\subsection{The quasi complex structure}
\label{cotacs}
The holomorphic tangent bundle $T'_{X}$
of a complex manifold $X$ admits  the complex linear 
automorphism given by multiplication by $i.$

\n
Via the isomorphism (\ref{tttiso}) we get an automorphism
$J$
of the real tangent bundle $T_{X}(\real)$ such that
$J^{2} = -Id_{T_{X}(\real)}.$

\n
The same is true for $T'^{*}_{X}$ using the dual map
$J^{*}.$

\bigskip
\begin{exe}
\label{qcexp}
{\rm
Show that, using a local chart  $(U;z):$
$$
J( \partial_{x_{j}} ) \, = \, \partial_{y_{j}} , \qquad \qquad
J( \partial_{y_{j}} ) \, = \, - \partial_{x_{j}},
$$
$$
J^{*}( dx_{j} ) \, = \, - d y_{j} , \qquad \qquad
J^{*}( d y_{j} ) \, = \, dx_{j}.
$$
}
\end{exe}

\bigskip
The various properties of tangent and cotangent bundles, e.g. 
(\ref{ttt}), (\ref{uup}), (\ref{tttiso}), (\ref{tttisod})
etc.
can be seen 
using $J$ via the eigenspace decomposition of 
$T_{X}(\comp)$ with respect to $J \otimes Id_{\comp}.$

\bigskip
The formalism of a quasi complex structures
allows to start with a complex vector space
$V$ and  end up with a display like (\ref{ttt}):
$$
V_{\real} \subseteq V_{\real}\otimes_{\real} \comp \, = \, V' \,  \oplus 
\, V''
$$
where $V_{\real}$ is the real vector space underlying
$V,$ $V' \simeq_{\comp} V$ and $V'' \simeq_{\comp} \ov{V'}$
and $V_{\real}\otimes_{\real} \comp$ has the conjugation operation.

\bigskip
A second equivalent point of view is detailed in Exercise
\ref{gphv} and starts with a $V$ as above and considers
$V \oplus \ov{V}$ instead.

\bigskip
We mention  this  equivalence in view 
of the use of Hermitean metrics on $T'_{X}.$

\n
These metrics are defined as special  tensors 
in $T'^{*}_{X} \otimes_{\comp} \ov{T'^{*}_{X}}.$

\n
However, in view of the fact that it is convenient
to conjugate them, it may be preferable to view them
as tensors in
the space with conjugation 
$T^{*}_{X}(\comp) \otimes_{\comp} T^{*}_{X}(\comp):$
using the canonical isomorphism (\ref{stcciso}) we can 
view the tensor $h$ as an element
of $T'^{*}_{X} \otimes_{\comp} T''^{*}_{X} \subseteq 
T^{*}_{X}(\comp) \otimes_{\comp} T^{*}_{X}(\comp).$

\n
This is also convenient in view of the use of the 
real alternating
form associated with a Hermitean metric which can then be viewed
as a real  element of $\Lambda^{2}_{\comp}(T^{*}_{X}(\comp)).$

\bigskip

\begin{defi}
\label{qcs}
{\rm 
({\bf Quasi complex structure}) 
A {\em quasi complex structure} on a real vector space $V_{\real}$  of 
finite even dimension $2n$ is a $\real-$linear automorphism 
$$
J_{\real}\, :\,  V_{\real} \,  \simeq_{\real}\,  V_{\real},
\qquad \qquad  J^{2} \, = \,  - 
Id_{V_{\real} }.
$$
}
\end{defi}

\begin{exe}
\label{descent}
{\rm
Show that giving a quasi complex structure
$(V_{\real}, J_{\real} )$ as in Definition \ref{qcs}  is equivalent to endowing 
$V_{\real}$
with a structure of complex vector space of dimension $n.$
(Hints: in one direction define $i \, v := J_{\real} (v);$ in the other
define $J_{\real}$ as multiplication by $i.$)
}
\end{exe}

\bigskip
Let $(V_{\real},J_{\real})$ be a quasi complex structure.

\medskip
\n
Let $V_{\comp}:= V_{\real} \otimes_{\real} \comp$ and
 $J_{\comp}= J_{\real} \otimes Id_{\comp} 
: V_{\comp} \simeq_{\comp} V_{\comp}$ be the 
complexification of $J_{\real}.$

\medskip
\n
The automorphism $J_{\comp}$ of $V_{\comp}$ has eigenvalues
$i$ and $-i.$

\medskip
\n
There are a natural  inclusion and a natural internal direct 
sum decomposition
$$
\fbox{$
V_{\real } \, \subseteq \, V_{\comp} \, = \, V' \, \oplus \, V''$}
$$
where

\begin{itemize}
    
\item
the subspace
$
V_{\real} \, \subseteq  \, V_{\comp}
$
is the fixed locus of the conjugation map
 associated with the complexification,

\item
$V'$ and $V''$ are the  $J_{\comp}-$eigenspaces
corresponding to the eigenvalues $i$ and $-i,$  respectively,

\item
since $J_{\comp}$ is real, i.e. it fixes $V_{\real} \subseteq V_{\comp},$ 
$J_{\comp}$  
commutes 
with the natural conjugation map and 
$V'$ and $V''$
are exchanged by this conjugation map,

\item
there are natural $\real-$linear isomorphisms coming from
the inclusion and the projections to the direct summands
$$
\fbox{$
V' \, \simeq_{\real}\,  V_{\real} \, \simeq_{\real} \, V''$}
$$
and complex linear isomorphisms
$$
\fbox{$V' \, \simeq_{\comp} \, \ov{V''},$} 
\qquad \qquad \fbox{$V'' \, \simeq_{\comp} 
\, \ov{V'},$}
$$
\item
the complex vector space defined by the complex structure (see
Exercise \ref{descent}) is $\comp-$linearly isomorphic to
$V'.$

\end{itemize}

\bigskip
\begin{exe}
    \label{vvv}
    {\rm
    Verify all the assertions made above.
    }
    \end{exe}

\bigskip
The same considerations are true for the quasi complex structure
$(V_{\real}^{*},  J_{\real}^{*}).$ We have 
$$
\fbox{$
V_{\real}^{*} \, \subseteq  \, 
V_{\comp}^{*} \, = \, V'^{*} \, \oplus \, V''^{*}, $}
$$
$$
\fbox{$
V'^{*} \, \simeq_{\real}\,  V_{\real}^{*}  \, \simeq_{\real} \, 
V''^{*}, $}
$$
$$
\fbox{$V'^{*} \, \simeq_{\comp} \, \ov{V''^{*}},$} 
\qquad \qquad \fbox{$V''^{*} \, \simeq_{\comp} 
\, \ov{V'^{*}}.$}
$$
In addition to the natural conjugation map stemming from the 
complexification, $V_{\comp}^{*}$ comes with the operational 
conjugation
map, for its elements are in 
$Hom_{\comp}(V_{\comp}, \real \otimes_{\real} \comp).$
However, these two operations coincide; see (\ref{ahfc}).

\bigskip
\begin{exe}
    \label{vvvd}
    {\rm
    Verify all the assertions made above.
    }
    \end{exe}

\begin{exe}
\label{gphv}
{\rm 
Let $\widetilde{W}:= W \oplus \ov{W}$ and $\iota: \widetilde{W} 
\lorw \widetilde{W}$ be the involution exchanging the summands.
Let $W_{0} \subseteq \widetilde{W}$ be the real subspace fixed by 
$\iota.$

\n
Show that there is a natural quasi complex structure on $W_{0}$ and 
that we get the following structure
$$
W_{0} \,  \subseteq \, {W_{0}}_{\comp}\, = \, 
W' \oplus W''
$$
endowed with the natural conjugation $c$  coming from the complexification.

\n
Show that there is a natural isomorphism
$l: {W_{0}}_{\comp} \simeq_{\comp} \widetilde{W},$ mapping
$W' \simeq_{\comp} W$ and $W'' \simeq_{\comp} \ov{W},$
which is 
compatible with the conjugation  and involution, i.e. $\iota\circ l = 
l \circ c.$

\n
Verify that
$$
W' \, = \, l^{-1}(W) \, = \, \{\, (w,w)\otimes 1 -
(i\,w, i\,w)\otimes i \, | \; w \in W \, \}
$$
and
$$
W'' \, = \, l^{-1} (\ov{W} ) \, = \, 
 \{\, (w,w)\otimes 1 +
(i\,w, i\,w)\otimes i \, | \; w \in W \, \}.
$$

\n
The complex vector space $F = Hom_{\real}(W, \comp)= W^{*} \oplus \ov{W}^{*}$
admits a conjugation-type map
$\ov{ (f,\ov{g})} := (g, \ov{f})$ (here $
\ov{g}(w) := \ov{ g(w)}$) with fixed locus
$$
F_{0} \, = \, \{ \, (f, \ov{f})\, | \; f \in W^{*} \, \}.
$$
Show that the map $m: F_{0} \otimes_{\real} \comp \lorw F$,
$m( (f,\ov{f})\otimes i ) := ( i\, f, \ov{ i\, f} )$
is an isomorphism compatible with the conjugations.
Verify that
$$
W'^{*} \, = \, m^{-1}(W^{*}) \, = \, \{\, (f,\ov{f})\otimes 1 -
(i\,f, \ov{i\,f})\otimes i \, | \; f  \in W^{*} \, \}
$$
and
$$
W''^{*} \, = \, m^{-1}(\ov{W}^{*}) \, = \, \{\, (f,\ov{f})\otimes 1 +
(i\,f, \ov{i\,f})\otimes i \, | \; f  \in W^{*} \, \}.
$$
}\end{exe}

\subsection{Complex-valued forms}
\label{cvf}
\begin{defi}
\label{defofcvf}
{\rm 
Let $M$ be a smooth manifold.
Define the {\em complex-valued smooth $p-$for\-ms}
as 
$$
\fbox{$
A^{p}(M)\, : = \,  E^{p}(M) 
\otimes_{\real} \comp \, \simeq \, C^{\infty}(M, T_{M}(\comp) ).
$}
$$
}
\end{defi}

\bigskip
The notion of exterior
differentiation  extends  
to complex-valued differential forms:
$$
\fbox{$ d\,:\,  A^{p}(M) \lorw A^{p+1}(M).$}
$$

\bigskip
\begin{rmk}
\label{eaht}
{\rm
If $M$ is a smooth manifold,
then we have the complex-valued version of the  de Rham
Isomorphism.
If $M$ is a smooth, oriented, compact manifold,
then  
we have 
the complex-valued versions of the 
Hodge theory statements 
of $\S$\ref{hfhit}. 
No new idea is necessary for these
purposes.
In the remaining part of these lectures, instead,
we are going to 
discuss  the aspects of Hodge theory which are specific
to complex,  K\"ahler and projective  manifolds.
}
\end{rmk}

\bigskip
Let $X$ be a complex manifold of dimension
$n,$ $x\in X,$ $(p,q)$ be 
a pair of non-negative integers and
define  complex vector spaces 
\begin{equation}
\label{tpq}
\fbox{$ 
\Lambda^{p,q}(T^{*}_{X,x}) \, : =  \,
\Lambda^{p}( T'^{*}_{X,x} ) \otimes \Lambda^{q}( T''^{*}_{X,x} ) 
\, \subseteq \,  \Lambda^{p+q}_{\comp}( T^{*}_{X,x} (\comp )   ).
$}
\end{equation}
 There is a canonical internal direct sum decomposition
of complex vector spaces
$$
\fbox{$
\Lambda^{l}_{\comp}( T^{*}_{X,x}(\comp) ) \; = \; \bigoplus_{p+q=l}{
\Lambda^{p,q}( T^{*}_{X,x} )}.$}
$$

\bigskip

\begin{exe}
\label{pqinv}
{\rm 
Verify that a holomorphic change of coordinates
leaves this decomposition invariant. In particular,
we can define smooth complex  vector bundles
$\Lambda^{p,q}(T^{*}_{X})$ and obtain an internal direct sum 
decomposition of smooth complex vector bundles
\begin{equation}
\label{pqdeco}
\fbox{$
\Lambda^{l}_{\comp}( T^{*}_{X}(\comp)) \, =\,  \bigoplus_{p+q=l}{
\Lambda^{p,q}( T^{*}_{X} )}. $}
\end{equation}
}
\end{exe}

\bigskip
\begin{defi}
\label{pqforms}
{\rm ({\bf $(p,q)-$forms})
The {\em space of $(p,q)-$forms} on $X$ 
$$
\fbox{$
A^{p,q}(X) \, :=  \, C^{\infty} (\, X\, ,\,  
\Lambda^{p,q}(T^{*}_{X}\, )\, )
$}
$$
is the  complex vector space 
of smooth sections of the smooth complex vector bundle
$\Lambda^{p,q}(T^{*}_{X}).$
}
\end{defi}

\bigskip

\begin{exe}
\label{locformpq}
{\rm
Verify the following statements.

\n
$A^{p,q}(X)$ is the complex vector space
of  smooth $(p+q)-$forms $u$ which can be written, locally on $U,$
with respect to any holomorphic chart $(U;z),$
as
$$
u = \sum_{|I|=p,|J|=q}{ u_{IJ}(z)\, dz_{I} \wedge d\overline{z}_{J} }.
$$
There is  a canonical direct sum decomposition,
$$
\fbox{$
A^{l}(X)\;  =\;  \bigoplus_{p+q=l}{A^{p,q} (X)} $}
$$
and
$$
\fbox{$ d\,(A^{p,q})\,  \subseteq\,  A^{p+1,q}(X) \, \oplus\,  A^{p, q+1}(X). $}
$$
}
\end{exe}

\bigskip
Let $l =p+q$ and consider the natural
projections  
$$
\pi^{p,q}\, :\,  A^{l}(X) \lorw A^{p,q}(X) \subseteq A^{l}(X).
$$

\bigskip

\begin{defi}
\label{dpdpp}
Define operators
$$
d': A^{p,q}(X) \lorw  A^{p+1,q}(X), \qquad\qquad 
d'': A^{p,q}(X) \lorw  A^{p,q+1}(X)
$$
$$
\fbox{$ d' \,  = \,   \pi^{p+1,q}\circ d,$}
 \qquad  \qquad
\fbox{$ d''\,  = \,  \pi^{p,q+1}\circ d.$}
$$
\end{defi}

\bigskip
Note that
\begin{equation}
    \label{dddsq}
\fbox{$ d =d' +d'',$}  \qquad 
\fbox{$ d''^{2}\, =\, 0\,  =\,  d'^{2},$}  \qquad 
\fbox{$d'd''\, =\,  - d''d'.$}
\end{equation}

\bigskip
\begin{exe}
\label{d'd''loc}
{\rm Verify that if,
in local coordinates,
$$
u\, =\,  \sum_{|I|=p, |J|=q}{ u_{IJ}\, dz_{I} \wedge d \overline{z}_{J}}
\in A^{p,q}(X),
$$
then
$$
d' u\, =  \,    \sum_{I,J}\sum_{j}{ \frac{\partial u_{IJ} }{ 
\partial_{z_{j}} }  \, dz_{j} \wedge dz_{I} \wedge
d\overline{z}_{J} },\qquad d''u \, =  \,
\sum_{I,J}\sum_{j}{ \frac{\partial u_{IJ} }{
\partial_{\overline{z}_{j}} }  \,  d\overline{z}_{j} \wedge dz_{I} \wedge
d\overline{z}_{J} }.
$$
}
\end{exe}

\bigskip

\begin{exe}
\label{conjdd}
{\rm Show that 
$$
\fbox{$
\ov{d'} \, = \, d'',$} \qquad 
\fbox{$\ov{d''} \, = \, d'$},
$$
where the conjugation symbol denotes  the operational conjugation  defined in
$\S$\ref{conjugations}.
}
\end{exe}

\subsection{Dolbeault  and Bott-Chern cohomology}
\label{dcg}
In this section, we look at new cohomology groups stemming
from the complex structure, the Dolbeault cohomology groups
$H^{p,q}_{d''}(X)$ and the Bott-Chern cohomology groups
$H^{p,q}_{BC}(X).$

\bigskip
The former ones appear in the \underline{complex-manifold-analogues}
of the results of $\S$\ref{hfhit}, which should not to be confused
with the \underline{complex-valued-analogues} of Remark \ref{eaht}.

\n
The latter ones will be used to show the important fact
that on a compact
K\"ahler manifold the Hodge Decomposition is independent
of the K\"ahler metric used to obtain it.

\bigskip
\begin{defi}
\label{tdcdef}
{\rm
({\bf The Dolbeault complex})
Fix $p$ and $q.$
The {\em Dolbeault complex} of $X$ is 
the complex of vector spaces
\begin{equation}
\label{dc}
\fbox{$ 0 \lorw A^{p,0}(X) \stackrel{d''}\lorw A^{p,1}(X) \stackrel{d''}\lorw 
\ldots \stackrel{d''}\lorw
A^{p,n-1}(X) \stackrel{d''}\lorw A^{p,n}(X) \lorw 0.$} 
\end{equation}
We have the analogous complex
\begin{equation}
\label{dcd}
\fbox{$ 0 \lorw  A^{0,q}(X) \stackrel{d'}\lorw A^{1,q}(X) \stackrel{d'}\lorw 
\ldots \stackrel{d'}\lorw
A^{n-1,q}(X) \stackrel{d''}\lorw A^{n,q}(X) \lorw 0.$}
\end{equation}
}
\end{defi}

\bigskip
\begin{defi}
 \label{dpqc}
{\rm  ({\bf Dolbeault cohomology})
The {\em Doulbeault cohomology groups} are the cohomology
 groups  of the complex (\ref{dc})
$$
\fbox{$
 H^{p,q}_{d''}(X)\,: =\, \frac{\ke{\,  d''\,:\, A^{p,q}(X)\lorw 
 A^{p,q+1}(X)  } }{
\im{ d''\,:\, A^{p,q-1}(X)\lorw A^{p,q}(X)} }
$}
$$
Define also
$$
\fbox{$
 H^{p,q}_{d'}(X) \, : = \,
 \frac{\ke{ \,  d' \, : \,  A^{p,q}(X)\lorw A^{p+1,q}(X)  } }{
\im{ \,  d'\, :\,  A^{p-1,q}(X)\lorw A^{p,q}(X)}}. $}
$$
}
\end{defi}

\bigskip
\begin{exe}
\label{isodcdd}
{\rm
Show that conjugation induces
canonical complex linear isomorphisms
$$
H^{p,q}_{d''}(X) \; \simeq_{\comp} \;
\ov{ H^{q,p}_{d'}(X) }.
$$
}
\end{exe}

\bigskip
\begin{tm}
    \label{gdl}
    ({\bf Grothendieck-Dolbeault Lemma})
    Let $q>0.$ 
    Let $X$ be a complex manifold and  $u \in A^{p,q}(X)$
    be such that $d''u=0.$ Then, for every point $x \in X,$
    there is an open neighborhood $U$ of $x$ in $ X$ and a 
    form $v \in A^{p, q-1}(U)$ such that 
    $$ u_{|U} \, = \, d''v .$$
    \end{tm}
    {\em Proof.} See \ci{gh} page 25.
    \blacksquare

    \bigskip
    
    \begin{exe}
\label{gdd'}
{\rm 
({\bf Grothendieck-Dolbeault Lemma for $d'$.)}
State and prove the result analogous to
Theorem \ref{gdl} for $d'.$
}
\end{exe}

\bigskip
\begin{rmk}
\label{weildol}
{\rm This remark is the Dolbeault counterpart of the
Weil-de Rham isomorphism Theorem  of Remark 
\ref{weil}. We have the  fine sheaves ${\cal A}^{p,q}_{X}$ 
of germs of local $(p,q)-$forms,
the sheaf 
$\Omega^{p}_{X}$ of germs
of holomorphic $p-$forms on $X,$
i.e. of the form $\sum_{|I|=p}{ u_{I} dz_{I}},$
$u_I$ holomorphic functions.
The Grothendieck-Dolbeault Lemma  \ref{gdl}
implies that
$\Omega^{p}_{X} \lorw ({\cal A}^{p,\bullet}_{X}, d'')$
is a resolution of $\Omega^p_X$ by fine sheaves. 

\smallskip
\n
We get the canonical {\em Dolbeault Isomorphisms}
$$
\fbox{$
H^{q}(X, \Omega^{p}_{X}) \; \simeq\;  H^{p,q}_{d''}(X).$}
$$
The $d'-$version of the Grothendieck-Dolbeault
Lemma  gives isomorphisms
$$
\fbox{$ H^p (X , \ov{ \Omega^q_X    } ) \; \simeq \;  
H^{p,q}_{d'}(X) $}
$$
where $\ov{ \Omega^q_X    }$ is the sheaf of germs
of anti-holomorphic $q-$forms on $X,$
i.e. of the form
$\sum_{|I|=q}{ u_{I} \, d\ov{z}_{I}},$
$u_I$ anti-holomorphic functions.
}
\end{rmk}

\bigskip
\begin{defi} 
\label{bcc}
{\rm
({\bf Bott-Chern cohomology})
The {\em Bott-Chern cohomology groups} of $X$ are defined
as the quotient spaces
$$
\fbox{$
H^{p,q}_{BC}(X)\, := \,
\frac{ A^{p,q}(X) \, \cap \,  \ke{\,d} }{ d'd''(A^{p-1,q-1}(X)) }.
$}
$$
}
\end{defi}

\bigskip
One can prove that if  $X$ is compact, then $\dim{H^{p,q}_{BC}(X)} < \infty.$

\bigskip
\begin{exe}
\label{mapsfbc}
{\rm
Verify that there are  natural maps
$$
H^{p,q}_{BC}(X) \lorw H^{p,q}_{d''}(X),
\qquad
H^{p,q}_{BC}(X) \lorw H^{p,q}_{d'}(X),
 \qquad H^{p,q}_{BC}(X) \lorw
H^{p+q}_{dR}(X, \comp)
$$
which, conveniently assembled, give 
bi-graded algebra homomorphisms. 
}
\end{exe}

\newpage
\section{Lecture 4: Hermitean linear algebra}
\label{l4}
We discuss Hermitean forms on a complex vector space, the associated 
symmetric and alternating forms, 
the induced inner product on the complexified exterior algebra
and we introduce the Weil operator.

\subsection{The exterior algebra on $V^{*}_{\comp}$}
\label{teaov}
Let things be as in $\S$\ref{cotacs}, i.e. we have,
for example,
$$
V^{*}_{\real} \, \subseteq \, V^{*}_{\comp}
\, = \, V'^{*} \, \oplus \, V''^{*}.
$$
Given the natural isomorphism
$\Lambda_{\real}(V^{*}_{\real}) \otimes_{\real} \comp \simeq
\Lambda_{\comp}(V^{*}_{\comp}) ,$
we have the inclusion of exterior algebras
$$
\Lambda_{\real}(V^{*}_{\real}) \, \subseteq \,
\Lambda_{\comp}(V^{*}_{\comp}) 
$$
where the complexified one carries the two identical
$c$-conjugation and operational conjugation
(\ref{ahfc})
which fix precisely the real exterior algebra.

\bigskip
Given any basis $\{ e_{j}^{*}\}$ for $V'^{*},$ the real vector
\begin{equation}
    \label{orh}
\fbox{$
(\frac{i}{2})^{n} \; e_{1}^{*} \wedge \ov{ e_{1}^{*}} \wedge \ldots
\wedge e_{n}^{*} \wedge \ov{ e_{n}^{*}} \; \in \; 
\Lambda^{2n}_{\real}(V^{*}_{\real}) $}
\end{equation}
gives an orientation for $V^{*}_{\real}$ which is independent
of the choice of basis and is therefore considered canonical.
See the calculation following Proposition  \ref{cmao}.

\bigskip
Define
$$
\fbox{$
\Lambda^{p,q}(V^{*}_{\comp} ) \, : = \, \Lambda^{p}_{\comp}(V'^{*})
\, \otimes_{\comp} \,  \, \Lambda^{q}_{\comp}(V''^{*}) 
$}
$$
which can be viewed as sitting naturally in 
$\Lambda^{p+q}_{\comp}(V^{*}_{\comp}).$

\bigskip
\n
Its elements can be written as complex linear combinations
of vectors of the form
$$
f_{1} \wedge \ldots \wedge f_{p} \wedge \ov{g_{1}} \wedge \ldots
\wedge \ov{ g_{q} }, \qquad  f_{j}, \, g_{k} \, \in V'^{*}.
$$

\bigskip
We have real linear isomorphisms induced by conjugation
$$
\Lambda^{p,q}(V^{*}_{\comp} ) \, \simeq_{\real} \, \Lambda^{q,p}(V^{*}_{\comp} ).
$$

\bigskip
We have a natural internal direct sum decomposition
$$
\fbox{$ 
\Lambda^{l}_{\comp}(V^{*}_{\comp})\, = \, 
\bigoplus_{p+q =l} \,
\Lambda^{p,q}(V^{*}_{\comp} ).
$}
$$

\bigskip
If $V^{*}_{\real}$ has a metric $g^{*}$, then we get 
the induced metric, still denoted by $g^{*},$
and the star operator
on $\Lambda_{\real}(V^{*}_{\real})$
 for the given metric  and the  orientation (\ref{orh}).

\subsection{Bases}
\label{bas}
Let $ \{ e_{1}, \ldots, e_{n} \}$ be a basis for $V'$ and
$\{ e_{1}^{*}, \ldots, e_{n}^{*} \}$ be the dual basis for $V'^{*}.$

\bigskip
\begin{exe}
\label{ebdb}
{\rm
Show that
$$
\{ \, \ov{e_{1}}, \ldots, \ov{e_{n}} \, \}   
$$
is a basis for  $V''.$

\n
Show that
$$
\{\, \ov{e_{1}^{*}}, \ldots, \ov{ e_{n}^{*}} \, \}
$$
is the corresponding dual basis for $V''^{*}$ where the conjugation 
is the one associated with the complexification, or equivalently,
the operational one.

\n
Show the analogues of (\ref{natb}), (\ref{natbx}), (\ref{dzdzbar})
and (\ref{dzdzbarx}) i.e.:
show that 
$$
\fbox{$ x_{j} \, := \, e_{j} + \ov{e_{j}}, $} \qquad \qquad
\fbox{$ y_{j} \, := \, 
i\, ( e_{j} - \ov{e_{j}} ) $}
$$
form  a basis for $V_{\real}$
whose dual basis  for $V_{\real}^{*}$ is given by
$$
\fbox{$ x_{j}^{*} \, = \, \frac{1}{2} \, ( e_{j}^{*} + \ov{ e_{j}^{*} } ) ,
$}
\qquad 
\qquad 
\fbox{$ y_j^{*} \, = \, \frac{1}{2i} \, ( e_{j}^{*} - 
\ov{e^{*}_{j} } )$}
$$
and   that
$$
\fbox{$ e_{j}\, = \, \frac{1}{2} ( x_{j} - i\, y_{j}), $} \qquad 
\fbox{$ \ov{e_{j}} \, = \, \frac{1}{2} ( x_{j} + i\, y_{j}), $} 
$$
$$
\fbox{$ e_{j}^{*} \, = \, x_{j}^{*} + i\, y_{j}^{*}, $} \qquad
\qquad 
\fbox{$ \ov{ e_{j}^{*} } \, = \, x_{j} - i\, y_{j}^{*}. $}
$$
}
\end{exe}

\subsection{Hermitean metrics}
\label{hm}
\begin{defi}
\label{defhmn}
{\rm ({\bf Hermitean forms})
A {\em Hermitean form} on a finite dimensional complex vector space
$W$ is a $\comp-$bilinear form
$$
h \, : \, W \times \ov{W} \lorw \comp
$$
such that 
$$
h(v,w) \, = \, \ov{ h(w,v) }, \quad \forall v,\, w\, \in W.
$$
A {\em Hermitean metric} on $W$ is a positive definite
Hermitean form, i.e. one for which
$$
h(v,v) \, > \, 0, \quad \forall  \; 0 \neq v\, \in \,  W.
$$
}
\end{defi}

\bigskip
\begin{exe}
\label{hsap}
{\rm Verify the following assertions.
Let $h$ be a Hermitean form on a complex vector space $W.$
Note that
\ci{we}, $\S I.2$ uses a different convention, where
$f$ is anti-linear in the first variable so that there are sign 
differences in what follows. 

\n
The real bilinear form
on the real vector space $W=\ov{W}$
given by
$S_{h} = Re \,h\, :\,  W \times W \lorw \real $ is 
symmetric.

\n
The form $S_{h}$ is positive definite iff $h$ is a  Hermitean metric.

\n
The real bilinear form $A_{h} = Im \, h \,: W \times W \lorw \real $ is 
anti-symmetric.

\n 
We have, dropping the sub-fix ``h''
$$
S(w,w')=  A(i\,w,w')= - A(w,i\,w'), \qquad A(w,w')=S(w,i\, w') = 
-S(i\,w,w'),
$$
$$
S(i\,w,i\,w') = S(w,w') , \qquad A(i\,w, i\,w') = A(w, w').
$$
Let $S'$ be a symmetric $\real-$bilinear form on the
real vector space $W$ which is invariant under
the $\comp-$linear automorphism of $W$
given by $w \to i\,w.$ Then
$A' (w,w') : = S'(w, i\,w')$ defines a real  alternating form
on the real vector space $W$ and  the form
$S' + i \,A'$ is Hermitean.

\n
Let $A''$ be an alternating $\real$ form on the real vector space
$W$ invariant under $w \to i\, w.$ Then $S''(w,w'):=
A''(i \, w,w')$ defines a symmetric bilinear form
and the form $S''+i\, A''$ is Hermitean.

\n
There are bijections between the following three sets:
the set of Hermitean forms on $W,$
the set of real symmetric  bilinear forms on $W$ invariant under $w \to 
i\, w,$
the set of
real alternating  bilinear forms on $W$ invariant under $W \to i\, W.$ 
}
\end{exe}

\bigskip
\begin{defi}
    \label{twof}
{\rm  ({\bf The alternating bilinear form associated with $h$)}
The {\em alternating bilinear  form} associated with
a  Hermitean form  $h$ is 
$$
\omega_{h} \, : = \, - A_{h} \, = \,   -   \, \mbox{Im} \, h.
$$
}
\end{defi}

\bigskip
A Hermitean form $h$ is a tensor in $W^{*} \otimes_{\comp} \ov{W}^{*}.$

\n
If $\{ e_{j} \}$ is a basis for $W,$ then we get
the dual basis $\{ e_{j}^{*} \}$ for $W^{*}$ and the dual basis
$\{ \e_{j}^{*} \}$ for $\ov{W}^{*}.$

\n
We have that $\e_{j}^{*} = {\ov{e_{j}^{*}}}^{r},$ i.e. 
$\e_{j}^{*}(w)  = \ov{ e_{j}^{*} (w)}$
and we can write
$$
h \, = \, \sum_{j,k}{ h_{jk} \, e_{j}^{*} \otimes {\ov{ \e_{k}^{*}} }}^{r},
\qquad 
h_{jk}:= h(e_{j}, e_{k}), \qquad h_{jk} = \ov{h_{kj}}. 
$$

The slight problem with this set-up is that we wish to perform
conjugation operations in order, for example, to give expressions
for $S_{h}$ and $\omega_{h}$ 
and the expression
$\ov{ e_{j}^{*} \otimes {\ov{ \e_{k}^{*}}}^{r} }^{r} =  
{\ov{ e_{j}^{*}}}^{r} \otimes  
e_{k}^{*}$ does not represent, strictly speaking, an equality
in $W^{*} \otimes_{\comp} \ov{W}^{*}.$

\bigskip
The upshot of  Exercise \ref{gphv} is that we can, 
by setting $V_{\real } := W_{0},$ view $W$ as a $V',$
$\ov{W}$ as a $V''$ etc. and  view $h$ as a $\comp-$bilinear map
$$
\widetilde{ h} \, : \, V'  \times V'' \lorw \comp
$$
with, setting $v': = l^{-1}(v)$ (see Exercise \ref{gphv}) etc.,
$$
h(v,w) = \widetilde{h}(v', \ov{w'}).
$$
The tensor $\widetilde{h} \in V'^{*} \otimes V''^{*} $ can now be viewed in
$V^{*}_{\comp}\otimes_{\comp} V^{*}_{\comp} $
and as such it can be conjugated, 
using the operational conjugation,
in a  way  compatible with all the isomorphisms and conjugations 
considered above.

\bigskip
We are now free to write, with abuse of notation,
\begin{equation}
    \label{efh} \fbox{$
h \, = \, \sum_{j,k}{ h_{jk} \, e_{j}^{*} \otimes {\ov{ e_{k}^{*} } }^{o}} ,$}
\qquad \qquad \fbox{$ h_{jk} \, = \, \ov{h_{kj}} $}
\end{equation}
and we are free to conjugate tensors.

\bigskip
Of course we choose to write $  {\ov{ e_{k}^{*} } }^{o}$ simply
as ${\ov{ e_{k}^{*} } }.$

\begin{exe}
    \label{emoncn}
    {\rm ({\bf The Euclidean metric on $\comp^{n}$})
Let $\{e_{j}\}$ be the standard basis for $\comp^{n}$
and define
$$
h
\, ( \sum_{j}{ a_{j}e_{j}}, \sum_{k}{ b_{k} e_{k} } ) 
\, : = \, \sum_{j}{ a_{j} \ov{b}_{j} }.
$$
Verify, using the bases of Exercise \ref{ebdb}, that 
$$
h\, = \, \sum_{j}{e_{j}^{*} \otimes \ov{e_{j}^{*}} },
$$
$$
S_{h} \, = \, Re\, h \, = \, \sum_{j}{ ( \, x_{j}^{*} \otimes
x_{j}^{*} + y_{j}^{*} \otimes
y_{j}^{*}  \, ) }
$$
$$
\omega_{h} \, = \, - Im\, h  \, = \, - A_{h} \, = \, 
\sum_{j}{ x_{j}^{*} \wedge y_{j}^{*} }
\, = \, 
\frac{i}{2}\,
\sum_{j}{  \, e_{j}^{*} \wedge \ov{ e_{j}^{*}} } .
$$
The ordered basis
$\{ \,  x_{1}^{*}, y_{1}^{*}, \ldots, x_{n}^{*}, y_{n}^{*} \, \}$
is orthonormal with respect to the Euclidean metric
on ${\comp^{n}}^{*} = \real^{2n}$ given by the dual metric
$S_{h}^{*}$. Endow this latter space with the orientation
(\ref{orh}). 
Verify that 
$$
 \frac{1}{n!} \, \omega_{h}^{n} \, = \, dV_{S_{h}^{*}}.
$$
}
\end{exe}

\bigskip
Using the expression (\ref{efh}) for the Hermitean
form $h,$ we deduce that
\begin{equation}
\label{efo}
\fbox{$
\omega_{h} \, = \, \frac{i}{2}\,  \sum_{j,k}{ h_{jk}\, 
e_{j}^{*} \wedge \ov{ e_{k}^{*} } }, $}
\qquad \qquad \fbox{$ h_{jk} \, = \, \ov{h_{kj}} $}
\end{equation}

\bigskip
The $2-$form $\omega_{h}$ is a real $(1,1)-$form.
By Exercise \ref{hsap}, 
giving a real $(1,1)-$form $\omega$
which, as an alternating form on $V',$ is invariant under multiplication
by $i$ on $V',$ is equivalent to giving
a Hermitean form $h_{\omega}.$ In this case,
we have  $\omega_{h_{\omega}} = \omega.$

\bigskip
The Hermitean form is a metric iff $||h_{jk}||$ is positive definite
(with respect to  any basis) which in turn is equivalent to the associated 
$(1,1)-$form $\omega_{h}$ being {\em positive}, i.e. 
being such that
$i\, \omega_{h} (v', \ov{v'}) >0$ for every $0 \neq v' \in 
V'.$

\bigskip
The Graham-Schmidt process ensures that if $h$ is a Hermitean metric,
then we can find a {\em unitary basis} for $h,$ i.e. a basis
$\{ e_{j} \}$ for $V'$ such that
$h(e_{j}, e_{k}) = \delta_{jk}$ so that
$$
h\, = \, \sum_{j}{ e_{j}^{*} \otimes \ov{ e_{j}^{*} } } 
$$
and 
\begin{equation}
    \label{gft}
\omega_{h} \, = \, \frac{i}{2} \,
\sum_{j}{ e_{j}^{*} \wedge \ov{ e_{j}^{*} } } \, = \,
\sum_{j}{ x_{j}^{*} \wedge \ov{ y_{j}^{*} } }
\end{equation}
from which  it is apparent that $\omega$ is a real $(1,1)-$form.

\bigskip
\n
We also have
$$
S_{h} \, = \, Re\, h \, = \, \sum_{j}{ ( \, x_{j}^{*} \otimes
x_{j}^{*} + y_{j}^{*} \otimes
y_{j}^{*}  \, ) }
$$
and  the orientation (\ref{orh}), giving the volume element
$$
dV_{S^{*}_{h}} \, = \, x_{1}^{*} \wedge y_{1}^{*} \wedge \ldots
\wedge x_{n}^{*} \wedge y_{n}^{*}.
$$
A straightforward calculation, i.e. taking
the $n-th$ exterior power of the expression (\ref{efo}),  gives the following
\begin{pr}
    \label{wvo}
    Let $h$ be a Hermitean metric on a complex vector space
    $V'$. Using the standard orientation on the space
    $V_{\real}^{*}$ and the metric $S_{h}^{*}$ on $V_{\real}^{*},$ we have the 
    equality:
    $$
   \frac{1}{n!} \, \omega^{n} \, = \, dV_{S_{h}^{*}}.
    $$
\end{pr}

\bigskip
\begin{rmk}
\label{wineq}
{\rm 
{\bf (Wirtinger Inequality)}
If $h$ is a Hermitean metric on a complex vector space $W,$
then we get two positive top  forms on the real vector space underlying
$W^{*}$, i.e.  $\frac{1}{n!} \,\omega_{h}^{n}$ and $dV_{S^{*}_{h}}.$
By Proposition \ref{wvo} they coincide.

\n
If $U \subseteq W$ is a real vector subspace of dimension
$2k,$ then we have two positive  forms
on $U:$ the restriction  ${S_{h}^{*}}_{|U}$ and the restriction
of $\frac{1}{k!} \, {\omega^{k}_{h}}_{|U}.$

\n
Wirtinger inequality states that
$$
\frac{1}{k!}\, {\omega^{k}_{h}}_{|U} \, \leq \;
{S_{h}^{*}}_{|U}
$$
and equality holds iff $U\subseteq W$ is a complex subspace.
See \ci{mu}, page 88 and the discussion that follows,
which culminates with \ci{mu}, Theorem 5.35, 
concerning volume-minimizing submanifolds of $\pn{n}.$
}
\end{rmk}

\subsection{The inner product  and the $\star$ operator
on the complexified exterior algebra $\Lambda_{\comp}(V^{*}_{\comp})$}
\label{tipsocea}
Let $h$ be a Hermitean metric on $V'.$
The metric $ S_{h}^{*}$ on $V^{*}_{\real}$
induces a metric, $g^{*},$ on $\Lambda_{\real}(V_{\real}^{*})$ 
and, using the orientation (\ref{orh}), the  $\star$ operator
on $\Lambda_{\real}( V^{*}_{\real}).$

\bigskip
The metric $g^{*}$ is a symmetric, positive definite  bilinear form on 
$\Lambda_{\real}(V_{\real}^{*}).$
We consider its complexification, i.e. the $\comp-$bilinear
form $g^{*} \otimes Id_{\comp}$ on $\Lambda_{\comp}(V_{\comp}^{*}),$
and define a Hermitean metric
$\langle \, , \, \rangle$  on 
$\Lambda_{\comp}(V_{\comp}^{*})$ by setting
\begin{equation}
    \label{eahvs}
    \fbox{$
\langle \, u, v \, \rangle \, : = \, (g^{*}\otimes Id_{\comp}) \, 
(u,\ov{v}).$}
\end{equation}

\bigskip
\begin{exe}
    \label{oist}
    {\rm 
Verify that $\langle \, , \, \rangle$  is a Hermitean metric, that
the spaces $\Lambda^{l}_{\comp}( V^{*}_{\comp})$ are mutually 
orthogonal, that the spaces $\Lambda^{p,q}( V^{*}_{\comp})$
are mutually orthogonal and that, given  a  $h-$unitary basis
for $V'$:
\begin{equation}
    \label{totoo}
|| \, e^{*}_{J} \wedge \ov{ e^{*}_{K}} \, || \, = \, 2^{ |J| + |K|   }.
\end{equation}
(See also Exercise \ref{ittwo})
}
\end{exe}

\bigskip
\begin{defi}
    \label{defofccst}
    {\rm
    The {\em $\star$ operator} on $\Lambda_{\comp}( V^{*}_{\comp} )$
    is defined to be the $\comp-$linear extension
    of the $\star$ operator on $\Lambda_{\real}( V^{*}_{\real} ).$
    }
    \end{defi}

    \bigskip
\begin{exe}
    \label{defprop}
    {\rm 
    Show that
\begin{equation}
    \label{sthesat}
    u\wedge \overline{\star \, v} = \langle u,v \rangle \; dV_{S^{*}_{h}},
    \qquad
    \forall \, u,\, v \; \in \Lambda^{l}_{\comp}(V^{*}_{\comp}).
    \end{equation}
Using (\ref{sthesat}), show that 
this new, extended $\star $ operator gives isometries
\begin{equation}
    \label{pqnpq}
\Lambda^{p,q}( V^{*}_{\comp} ) \simeq_{\comp}  \Lambda^{n-q,n-p}
( V^{*}_{\comp} ).
\end{equation}
In particular, 
\begin{equation}
    \label{pistco}
   \star  \; \pi^{p,q} \, = \, \pi^{n-q, n-p}\; \star.
   \end{equation}
    Finally, observe that since the real dimension of $V_{\real}$
   is even, (\ref{stst}) implies that 
   \begin{equation}
       \star \star_{| \Lambda^{p,q}(V^{*}_{\comp})  } \,  = \,
       (-1)^{p+q} Id_{\Lambda^{p,q}(V^{*}_{\comp})}.
       \end{equation}
   }
   \end{exe}

   \bigskip
The explicit form of this new extended $\star$ operator using
the real $S_{h}^{*}-$orthonormal ordered basis $\{ x_{1}^{*}, y_{1}^{*} , 
\ldots , 
x_{n}^{*}, y_{n}^{*} \}$ associated with a $h-$unitary basis $e_{j}$
is identical to the non-complexified one. 

\bigskip
The expression of the $\star$ operator using the  
$ e_{j}^{*},$ $ \ov{e_{j}^{*}}$ basis 
can be found in \ci{we}, pagg 19-20.

\subsection{The Weil operator}
\label{twop}
The Weil operator can be defined for any Hodge structure.
It is very convenient in view of the definition
and the use of polarizations of pure Hodge structures;
see $\S$\ref{hs}.

\bigskip
Here we look at the exterior algebra
 $\Lambda_{\comp}( V^{*}_{\comp})$ which is a direct sum 
 of the   weight $l$ pure Hodge structures 
  $\Lambda^{l}_{\comp}( V^{*}_{\comp});$ see $\S$\ref{teaov}.

\begin{defi}
    \label{cweil}
    {\rm 
    {(\bf The Weil operator)}
   The {\em Weil operator} is the complex linear isomorphism
   $\Lambda_{\comp}( V^{*}_{\comp}) \lorw \Lambda_{\comp}( 
   V^{*}_{\comp})$ induced by 
   multiplication by $i$ on $V',$ i.e.
   $$
   C \, : = \,  \sum_{p,q}{ i^{p-q} \, \pi^{p,q}}.
   $$
   }
   \end{defi}
   
   \bigskip
   The Weil operator is real, i.e. it preserves
   the real subspace $\Lambda_{\real}(V_{\real}).$

   \bigskip
   By (\ref{pistco}),  we have
 \begin{equation}
     \label{cstco}
     C \, \star \;  - \;  \star \, C \, = \, 0
     \end{equation}
     Since $\dim_{\real}{V_{\real}^*} $ is even, 
    $$\star\, \star \, = \,  \sum_{l=0}^{2n}{ (-1)^{l} \pi^{l}} \, = \, w
$$
where $w$ is
     the so-called {\em de Rham operator} $w.$
     
     \bigskip
     The following relations follow:
     $$
     w \, = \, \star\, \star \,  \, = \, C^{2}, 
     $$
$$
\star^{-1} \, = \, w \, \star \, = \, \star \, w, \qquad
C^{-1} \, = \, w \, C \, = \, C \,w.
$$

\newpage
\section{Lecture 5: The Hodge theory of Hermitean  manifolds}
\label{hmocm}
We discuss Hermitean metrics on complex manifolds, the $\Delta'$ 
and $\Delta''$ Laplacians, the $\Delta'$ and $\Delta''$ harmonic 
forms, the corresponding Hodge Theory on a compact complex manifold,
including  Kodaira-Serre Duality.

\subsection{Hermitean metrics on complex manifolds}
\label{hmocmll}
A Hermitean metric  $h$ on a complex manifold  $X$ is the assignment
of  a Hermitean metric
\begin{equation}
\label{hmd}
h(- , -)_{x} \, : \,  T'_{X,x} \otimes_{\comp}  \overline{T'}_{X,x} \lorw \comp  
\end{equation}
for every $x \in X$ varying smoothly  with $x,$ i.e.
such that the functions
$$
h_{jk}(z): =  (\partial_{z_{j}} , \partial_{z_{k}})
$$
are smooth on the open set $U.$

\bigskip
Recalling the discussion culminating with
(\ref{efh}),
 using the local chart $(U;z),$ the Hermitean metric $h$
can be expressed on $U$ in tensor form as
\begin{equation}
\label{h}
h \,  = \,  \sum_{jk}{ h_{jk}\, (z) dz_{j}\otimes d\overline{z}_{k}}.
\end{equation}

\bigskip
The  real $(1,1)-$form
$$
\omega= \omega_{h}= -  Im h \in A^{1,1} (X)
$$
is called the {\em associated $(1,1)-$form}  of the metric $h.$

\bigskip
Using the chart $U,$ it can be written as
$$
\omega= \frac{i}{2} \sum_{jk}{
h_{jk}(z) \,  dz_{j}\wedge d\ov{z}_{k} }, \qquad \qquad
h_{jk} = \ov{h_{kj}}.
$$
The metric $h$ can be recovered from the associated $(1,1)$ form.

\bigskip
\begin{exe}
\label{eohmpof}
{\rm ({\bf Existence of Hermitean metrics on complex manifolds})
Show that any complex manifold admits hermitian metrics on it.
}
\end{exe}

\bigskip
\begin{exe}
\label{funct}
{\rm ({\bf Restriction of Hermitean metrics})
Let $f: X^{n} \lorw Y^{m}$  a holomorphic map of complex 
manifolds of the indicated dimensions $n \leq m$
such that $df: T'_{X} \lorw f^{*}T'_{Y}$
is everywhere of  rank $n,$ i.e. injective.
Let $h_{Y}$ be a Hermitean metric on $Y.$
Show that one can induce a Hermitean metric
$h_{X}$ on $X$ such that $\omega_{h_{X}}= f^{*} \omega_{h_{Y}}.$
In particular, if $f:X \lorw Y$ is an embedding of complex manifolds,
then $h_{X}: ={h_{Y}}_{|X}$ is a Hermitean metric
and $\omega_{h_{X}}= {\omega_{h_{Y}}}_{| X}.$
See \ci{gh}, page 29.
}
\end{exe}

\bigskip
An important related fact is {\em Wirtinger's Theorem}, \ci{gh},
page 31,
\ci{mu}, page 88. See also Proposition \ref{wvo} and Remark
\ref{wineq}.

\bigskip
\begin{exe}
    \label{wt}
    {\rm ({\bf Wirtinger's Theorem})
    Prove Wirtinger's Theorem: let $Y \subseteq X $ be a complex 
    submanifold  of complex dimension $k,$ $h$  be a Hermitean metric on $X$
    and $\omega$ be the associated form; then
    \begin{equation}
	\label{wtgfm}
    vol (Y) \, = \, \frac{1}{k!} \, \int_{Y}{ \omega^{k}_{|Y} }.
    \end{equation}
    }
 \end{exe}

 \bigskip
  Note the following special feature of complex geometry
  expressed by (\ref{wtgfm}): the volume of $Y$
    is expressed as the integral over $Y$ of a globally defined
    differential form on $X.$ This does not occur in general
    in the real case.
    See \ci{gh}, page 31. 

\bigskip
\begin{exe}
\label{fsm}
{\rm ({\bf The Fubini-Study metric on $\pn{n}$})
Verify all the following assertions.
See  \ci{gh} page 30. See also \ci{mu}, pages 86-87.

\n
 Let $ [\, x_{0}: \ldots : x_{n}\, ]$ be homogeneous coordinates
on $\pn{n}.$ The expression $log (|x_0|^{2} +  \ldots  + | 
x_{n}|^{2} )$ is well-defined on $\comp^{n+1}\setminus \{0 \}.$
The differential form 
$$
\omega':= \frac{i}{2\pi} \, d'\,d'' \,\log
(|x_0|^{2} +  \ldots  + | 
x_{n}|^{2} ) \; \in\,  A^{1,1}( \comp^{n+1}\setminus \{0 \} )
$$
is $\comp^{*}-$invariant. Show that it   descends to a $(1,1)-$form
$\omega \in A^{1,1}(\pn{n}).$ 
On the chart $U^{0}$ (and on any chart, keeping track of indices)
we have 
\[
\omega = \frac{i}{2\pi} \, d'\,d'' \,\log
(1  +   |z_{1}|^{2}+\ldots  + | z_{n}|^{2} )=
\frac{i}{2\pi} \,
\left[
\frac{\sum_{j}{ dz_{j}\wedge d\overline{z}_{j}   } }{
1+ \sum_{j}{|z_{j}|^{2}}  } 
- 
\frac{ \left( \sum_{j}{ \overline{z}_{j} \, dz_{j} }     
\right) \wedge 
\left( \sum_{j}{ {z}_{j}\, d\overline{z}_{j} }
\right)  }{  \left( 1+ \sum_{j}{|z_{j}|^{2}}  \right)^{2}     } \right]
\]
At the point $[\,1: 0 : \ldots :  0\,]$
$$
\omega = \frac{i}{2\pi} \sum_{j}{ dz_{j}\wedge d\overline{z}_{j}} >0
$$
so that $\omega$ is the associated $(1,1)-$form for a Hermitean metric
defined on $\pn{n}.$ 
This metric is called the {\em Fubini-Study} metric of $\pn{n}.$

\n
Show that $d\omega =0,$ i.e. the Fubini Study metric is K\"ahler. 

\n
Show that  $\int_{\pn{n}}{\omega^{n}}=1.$

\n
Show that $[\omega] \in H^{2}(X, \real)$
is Poincar\'e dual to the homology class $\{H \} \in
H_{2n-2}(X, \real)$ associated with 
a hyperplane $H \subseteq \pn{n} $ and also that $\omega$
is the curvature form associated with a Hermitean metric
on the hyperplane bundle $L_{1}$ which is therefore
{\em positive} in the complex differential geometric sense
and {\em ample} in the algebraic geometric sense.
Finally, show that $[\omega]$ is the
first Chern class with $\real$ coefficients of
$L_{1}$ (this requires some non-trivial unwinding of the definitions;
see \ci{gh}, page 141).
}
\end{exe}

\subsection{The Hodge theory of a compact Hermitean manifold}
\label{thtchv}
Let   $(X,h),$  be a {\rm Hermitean manifold},
i.e. a complex manifold $X$  endowed with a  Hermitean metric $h.$

\bigskip
The smooth manifold underlying 
 $X$ carries the natural 
orientation given by the complex structure on $X.$

\bigskip
The metric $h$ gives rise to  a Hermitean
metric
$\langle \, , \, \rangle $
on the exterior algebra bundle $\Lambda (T_{X}^{*}(\comp)).$ See
(\ref{eahvs}).

\bigskip
We  extend, as in Definition \ref{defofccst},  the $\star$ operator
to complex-valued forms and
we have 
\begin{equation}
\label{stcvf}
u \wedge \overline{ \star \, v} \, = \,  \langle u,v\rangle \, dV
\end{equation}
where $dV$ is the Hermitean volume element, i.e.
the one associated with the orientation
and the Riemannian metric associated with the Hermitean metric.
See Exercise \ref{defprop}.

\bigskip
One has that
$$
\fbox{$
\star \, : \,   \Lambda^{p,q}_{\comp}( T^*_{X}) \lorw
\Lambda^{n-q,n-p}_{\comp}( T^*_{X})
$}
$$
is a complex linear isometry.

\bigskip
We will consider formal adjoints with respect to the following
metric.
\begin{defi} {\rm 
    Whenever the integral converges, e.g. for compactly supported
    forms,
 define
$$
\fbox{$ \langle \langle u,v \rangle \rangle \, = \,
\int_{X}{ \, \langle u , v \rangle \, dV }
$}
\, = \, \int_{X}{ u \wedge \ov{ \star  \, v} }.
$$
}
\end{defi}
    
\bigskip
\begin{defi}
\label{fohg}
({\bf The $d'$ and $d''$  Laplacians})
$$
\fbox{$ d'^{\star}=  - \star d'' \star,$} \qquad 
\fbox{$ d''^{\star}= (d'')^{\star}=  - \star d' \star $}
$$
$$
\fbox{$ \Delta'= d'd'^{\star} + d'^{\star}d' , $} \qquad 
\fbox{$ \Delta''=
d''d''^{\star} + d''^{\star}d''. $} 
$$
\end{defi}

\bigskip
The definition is motivated by the following

\begin{exe}
\label{tatadjok}
{\rm 
Show, as in Proposition \ref{dadjd*}, that
$d'^{\star}$ is the formal adjoint to
$d'$ and that  $d''^{\star}$ is the formal adjoint to $d''.$
Show that $\Delta'$ and $\Delta''$ are self-adjoint.
}
\end{exe}

\begin{defi}
\label{dbarpqh}
{\rm ({\bf Harmonic $(p,q)-$forms})
The space of  {\em $\Delta'-$harmonic $(p,q)-$forms}
of $X$ is 
$$
\fbox{$
{\cal H}^{p,q}_{\Delta'}(X): = \ke{ \, (\,\Delta' : 
A^{p,q}(X) \lorw A^{p,q}(X)}
\, ) $}
$$
and the
space of  {\em $\Delta''-$harmonic $(p,q)-$forms}
of $X$ is 
$$
\fbox{$ 
{\cal H}^{p,q}_{\Delta''}(X): = \ke{ \, ( 
\Delta'' : A^{p,q}(X) \lorw A^{p,q}(X)} \, ). $}
$$
}
\end{defi}

\bigskip
These spaces  depend on $h.$

\bigskip
The following results  are the
complex analytic versions of the Hodge Orthogonal 
Decomposition Theorem, of the Hodge Isomorphism Theorem
and  of the Poincar\'e Duality Theorem. One replaces 
the de Rham Cohomology with the Doulbeault cohomology 
and $\Delta-$harmonic forms with $\Delta''-$ harmonic forms
etc.

\begin{tm}
\label{hodtcav}
({\bf Hodge theory of compact complex manifolds})
Let $(X,h)$ be a \underline{compact} Hermitean
manifold. Then for every bi-degree
$(p,q):$ 

\bigskip
\n
(a) ({\bf Hodge Orthogonal Decompositions for $d'$ and $d''$})
there are orthogonal direct sum decompositions
$$
\fbox{$
A^{p,q}(X) = {\cal H}^{p,q}_{\Delta'}(X) \stackrel{\perp}\oplus
d'( A^{p-1, q}(X) )  \stackrel{\perp}\oplus
d'^{\star}( A^{p+1, q}(X)) ,
$}
$$
$$
\fbox{$
A^{p,q}(X) = {\cal H}^{p,q}_{\Delta''}(X) \stackrel{\perp}\oplus
d''( A^{p, q-1}(X) )  \stackrel{\perp}\oplus
d''^{\star}( A^{p, q+1}(X)) .$} 
$$

\bigskip
\n
(b) ({\bf Hodge Isomorphisms for $d'$ and $d''$})
there are  isomorphisms of \underline{finite}
dimensional complex vector spaces
$$\fbox{$
{\cal H}^{p,q}_{\Delta'}(X) \simeq H^{p,q}_{d'}(X), $}
$$
$$
\fbox{$
{\cal H}^{p,q}_{\Delta''}(X) \simeq H^{p,q}_{d''} (X).$}
$$

\bigskip
\n
(c) ({\bf Kodaira-Serre Duality}) 
the complex bilinear pairings
$$
H^{p,q}_{d'}(X) \times H^{n-p, n-q}_{d'}(X) \lorw \comp,
\qquad 
(u,v) \lorw \int_{X}{ u \wedge v},
$$
$$
H^{p,q}_{d''}(X) \times H^{n-p, n-q}_{d''}(X) \lorw \comp,
\qquad 
(u,v) \lorw \int_{X}{ u \wedge v},
$$
are   non-degenerate dualities.
In particular, there is a canonical isomorphism
$$
\fbox{$
H^{n-q}(X, \Omega^{n-p}_{X})
\, \simeq \,  
H^{q}(X, \Omega^{p}_{X} )^{\vee}. $}
$$
\end{tm}

\begin{exe}
    \label{ksdexe}
    {\rm 
Exercise: prove (b) and (c) using (a) as in Theorems
\ref{hit} and \ref{pd}.
}
\end{exe}

\newpage
\section{ Lecture 6:  K\"ahler manifolds}
\label{kama}
We discuss the K\"ahler condition, the fundamental identities of 
K\"ahler geometry, the Hodge Decomposition via the $d'd''-$Lemma
and Bott-Chern cohomology (as in \ci{dem1}) and some  topological 
implications of the Hodge Decomposition.

\bigskip
It is on compact K\"ahler manifolds that Hodge theory becomes a 
formidable tool which highligths some of the amazing properties
that these manifolds enjoy: non-vanishing of even Betti numbers,
parity of the odd Betti numbers, Kodaira-Serre symmetry and Hodge
symmetry for the Hodge numbers $h^{p,q}.$ See $\S$\ref{scof}.

\bigskip

It is costumary to denote a Hermitean
manifold $(X,h)$ also by
$(X, \omega),$ where $\omega =\omega_{h}.$ 

\subsection{The K\"ahler condition}
\label{tkc}

\begin{defi}
\label{defk}
{\rm ({\bf K\"ahler metric/manifold})
A  Hermitean metric $h$ on a complex manifold $X$  is called
{\em K\"ahler} if $d\omega=0.$

\n
 A complex manifold $X$ is said to be {\em K\"ahler} if
it admits a K\"ahler metric.

}
\end{defi}

\bigskip

\begin{defi}
\label{defproj}
{\rm ({\bf Projective manifolds)}
A complex manifold is said to be {\em projective}
if it admits a closed holomorphic  embedding 
in some projective space.
}
\end{defi}

\bigskip

\begin{ex}
\label{fsik} 
{\rm 
Any Riemann surface is automatically K\"ahler.
Any compact Riemann surface is projective.

\n
By Exercise  \ref{fsm},
$\pn{n}$ is 
K\"ahler. By Exercise \ref{funct},  any projective manifold is
K\"ahler. 

\n
Any compact complex torus
$\comp^{n}/ \Lambda,$ $\zed^{2n} \simeq \Lambda 
\subseteq \comp^{n}$ a full lattice,
is K\"ahler. However, ``most'' tori of complex dimension
at least two are not projective due to the 
{\em Riemann conditions}. See \ci{gh}, pagg 300-307.

\n
The two examples that follow, the Hopf surface and the Iwasawa
threefold, 
are compact complex manifolds which are not K\"ahler in view of
some of the special properties of compact  K\"ahler 
manifolds that we will establish later in this lecture.

\n
The Hopf surface \ci{dem1}, $\S 5.7$   is a  compact complex surface
which is  not K\"ahler, in fact its  first Betti number
is one, i.e. it is odd, and this is not possible on a compact
K\"ahler manifold; see Theorem  \ref{pqid}.

\n
The Iwasawa manifold \ci{dem1}, $\S 8.10$ is a compact complex 
threefold which is not K\"ahler since it admits  holomorphic
$1-$forms which are not $d-$closed, an impossibility on a compact 
K\"ahler manifold;  see Theorem \ref{holcl}.
}
\end{ex}

\bigskip
Let $(X, \omega)$ be a Hermitean manifold.
Since $\omega$ is a real  $(1,1)-$form, the three conditions
$$
d\omega =0, \qquad d' \omega =0, \qquad  d'' \omega =0,
$$
are equivalent to each other and to the condition
$$
\frac{ \partial h_{jk} }{ \partial z_{l} } \, = \,
\frac{ \partial h_{lk} }{ \partial z_{j} } , \quad  1 \leq j,k,l \leq 
n. 
$$

\bigskip
Let $(X,\omega)$ be  a K\"ahler manifold.
By Proposition
\ref{gphv}
$$
dV \, = \,  \frac{1}{n!}\, \omega^{n},
$$
where $dV$ is the volume element associated with
the Riemannian metric associated with $h$ and the canonical orientation.

\bigskip
There is an important difference
from the Hermitean case:  the right-hand-side
being closed   defines a cohomology class which
cannot be exact on a compact manifold because of Stokes'
Theorem and the fact that the integral of the 
left-hand-side cannot be zero. This implies that
all the relevant powers of $\omega$ define non-zero cohomology classes.
We have the following remarkable consequences for the topology
of a compact K\"ahler manifold.

\begin{tm}
    \label{cckmt}
Let $X$ be a compact K\"ahler manifold.
Then, for every $\, 0 \leq k \leq n,$ 
$$b_{2k}(X)\, = \, \dim_{\real}{H^{2k}(X, \real)} \, > \, 0.
$$

\n
The fundamental class $[V] \in H^{2n-2k}$ of a closed analytic 
 subvariety
of $X$   of dimension $\dim_{\comp}(V)=k$ is non-zero.
\end{tm}
{\em Proof.}
Let $h$ be a K\"ahler metric and $\omega$ the associated $(1,1)-$form.
Since $\omega$ is $d-$closed, so are $\omega^{k},$ for every $k \geq 
0.$
Since  $\frac{1}{n!}\int_{X}{\omega^{n}} = vol (X) \neq 0,$
Stokes' Theorem implies that the cohomology class
$[\omega^{n}]$ is non-zero and so are all the $[\omega^{k}]$,
$ 0 \leq k \leq n,$ due to the obvious relation
$[\omega^{n}] = [\omega^{k}] \cup [\omega^{n-k}].$

\n 
For the second statement, argue as follows. 
By Exercise \ref{wt}, and by the properties 
of the current of integration
along $V$ $\int_{V}$ (see Remark \ref{curr}),  we have that the pairing
$( [ \int_{V}], [\omega^{k}]) = \int_{V_{reg}}{ \omega^{k}_{|V_{reg}} } 
>0,$ so that $[V]= [ \int_{V}] \neq 0.$
\blacksquare

\bigskip
Let $(X, h)$ be a Hermitean manifold.
Let $g$ be the associated Riemannian metric on the underlying
smooth manifold, $ \nabla$ be the Levi-Civita
connection  on the real tangent bundle of $X$ associated with $g;$ in fact
there is one such acting on any tensor bundle. 

\medskip
\n
Let $J: T_X(\real) \lorw 
T_X(\real)$ be the quasi complex structure on $X,$
i.e. the associated  automorphism  of
the real  tangent bundle. Using a holomorphic
local chart $(U;z)$ around $x \in X,$ we have
 $J( \partial_{x_j} ) = \partial_{y_j}$
and $J( \partial_{y_j} ) = - \partial_{x_j}.$

\bigskip
\begin{exe}
\label{hmsasgj}
{\rm
Show that assigning a Hermitean metric $h$ on the complex manifold
$X$ is equivalent to assigning a $J-$invariant 
Riemannian metric on the underlying smooth manifold.
}
\end{exe}

\bigskip
\begin{tm}
\label{allhappytog}
({\bf Characterization of K\"ahler metrics})
Let $(X, h)$ be a Hermitean variety.
The following conditions are equivalent.

\medskip
\n
(a) $d \omega =0,$ i.e. $h$ is K\"ahler;

\smallskip
\n
(b) the types of the complexified tangent vectors are
preserved under parallel transport;

\smallskip
\n
(c) for every real parallel vector field $\eta$ along a smooth curve
$\gamma,$ $J  \eta$ is parallel along $\gamma;$

\smallskip
\n
(d) $\nabla \, J \, = \, 0,$ i.e. the almost complex structure is
parallel;

\smallskip
\n
(e) $ \nabla \, \omega \, = \, 0,$ i.e. the associated $(1,1)-$form 
is parallel;

\smallskip
\n
(f) $h$ admits a  potential locally on $X,$ i.e.
a function $\varphi$ such that, in local holomorphic
coordinates,
$g_{jk} = \frac{ \partial^2 \varphi}{\partial z_j \partial \overline{z}_k};$

\smallskip
\n
(g)  for every $x \in X$ there is a holomorphic
chart $(U;z)$ centered at $x$ such that
$g_{jk}(z) = \delta_{jk} + O (|z|^2);$

\smallskip
\n
(h) the torsion of the Hermitean metric $h$ is zero.
\end{tm}
{\em Proof.} See \ci{mok}, page 18 and \ci{gh} page 107.
See also \ci{dem1}, $\S 5.8.$ See also Exercise \ref{equivok}
 below.
\blacksquare

\bigskip
\begin{exe}
\label{equivok}
{\rm
Show that $(g)$ implies $(a)$ by directly differentiating the expression
$$
\omega \, = \, \frac{i}{2} \, \sum_{j,k}{ ( \delta_{jk} + [2] ) \,dz_{j} 
\wedge d\ov{z}_{k} }
$$
for $\omega$ at a given point.

\n
Show that $(a)$ implies $(g)$ as follows.

\n
Show that the Graham-Schmidt process implies the existence of 
holomorphic coordinates at any given point $x$ such that $h_{jk}(x)=
\delta_{jk}$ so that
$$
\omega \, = \, \frac{i}{2} \, \sum_{j,k,l}{ ( \delta_{jk} +
a_{jkl} \, z_{k} + a_{jk \ov{l} } \, \ov{z}_{k} ) \, dz_{j} \wedge 
d \ov{z}_{k} } \, = \, \frac{i}{2} \, \sum_{j,k}{ h_{jk}(z) \, 
dz_{j} \wedge d\ov{z}_{k} }.
$$
Show, using  that $h_{jk}(x) = \ov{ h_{kj}(x) } ,$ that
$$
\ov{a_{jkl} }\, = \, a_{kj\ov{l}}.
$$
Show, using $d\omega=0,$ that
$$
a_{jkl} \, = \, a_{lkj}.
$$
Use these relations to show that, by setting
$$
b_{juv} \, = \, - a_{vju},
$$
the  change of  coordinates defined by the equality 
$$
z_{j} \, = \,  w_{j} + \frac{1}{2}\, \sum_{u,v}{ b_{juv} \, w_{u} w_{v} }
$$
gives
$$
\omega \, = \,  \frac{i}{2} \, \sum_{j,k}{ ( \delta_{jk} + [2] ) \,
dw_{j} \wedge d\ov{w}_{k} }.
$$
}
\end{exe}

\bigskip

\begin{rmk}
\label{usef}
{\rm 
The  characterization
$(a)\, =\, (g)$  is used in $\S$\ref{tfiokg}
 to prove the fundamental  identities of K\"ahler geometry
by first showing them on $\comp^n$ with the Euclidean  metric
and then by observing that these identities
involve only the metric and its first derivatives
so that proving the Euclidean case is enough. See \ci{gh} page 115.
}
\end{rmk}

\subsection{The fundamental identities of K\"ahler geometry}
\label{tfiokg}
\begin{defi}
\label{deflandl}
{\rm 
({\bf The operator $L$ and its adjoint}).
Let $(X, h)$ be a Hermitean manifold and define
$$
L \, : \, A^{\bullet}(X) \lorw A^{\bullet +2}(X),  \quad
L(u)\, := \, \omega \wedge u,
$$
$$
L^{\star} \, : A^{\bullet}(X) \lorw A^{\bullet -2}(X),
\qquad L^{\star} (u) \, := \, 
\star^{-1} L \star
$$
}
\end{defi}

\bigskip
The operators $d',$ $d'',$
$d'^{\star},$ $d''^{\star}$
 $\Delta',$ $\Delta''$,
$L$ and $L^{\star}$ have bi-degree
$(1,0),$ $(0,1)$, $(-1,0),$ $(0,-1)$, $(0,0),$ $(0,0),$
$(1,1)$ and $(-1,-1).$

\bigskip
The {\em commutator $[A,B]$} of two operators of bi-degree $(a',a'')$ and 
$(b',b'')$ and of total degree $a=a'+a'',$ $b= b'+b''$ is defined as
\begin{equation}
    \label{degop}
    [A,B] \, : = \, A\,B  \, - \,  (-1)^{a\,b} \, B \, A.
\end{equation}

\bigskip

\begin{tm}
\label{bciokg}
({\bf Fundamental identities of K\"ahler geometry})

\n
Let $(X, h)$ be a K\"ahler manifold.
Then 
\begin{equation}
    \label{fikg3}
\fbox{$    [ d''^{\star}, L ] \, = \, i \, d' , $} \quad 
\fbox{$  \; [ d'^{\star} , L ] 
\, = \, -i \, d'', $} \quad
\fbox{$   [ L^{\star}, d'' ] \,  = \, -i \, d'^{\star} ,  $}
\quad
\fbox{$  \;  [ L^{\star}, d' ] \, = \, i \, d''^{\star}.$}
\end{equation}
\begin{equation}
    \label{fikg1}
\fbox{$[ d', d''^{\star} ] \, = \, 0, $}
\qquad  \fbox{$ [d'' ,d'^{\star}] \, = \, 0, $}
\end{equation}
\begin{equation}
\label{fikg2}
\fbox{$ \Delta \, = \, 2 \Delta' \, = \, 2 \Delta'', $}
\end{equation}
In addition,  $\Delta,$ $\Delta'$ and $\Delta''$
commute
with $d$
$\star,$ $ d'$, $d'',$ $d'^{\star},$ $ d''^{\star},$
$L$ and $L^{\star}.$

\n
Finally, $\Delta$ preserves the $(p,q)-$decomposition,
i.e. $\Delta$ commutes with $\pi^{p,q}.$
\end{tm}
{\em Proof.} The first and second  relations (\ref{fikg3})
are conjugate to each other and so are the third and fourth.

\n
The first and third are adjoint to each other.

\n 
It follows that the four relations are equivalent to each other.

\n
We prove the fourth relation in the guided  Exercise 
\ref{commrel}.

\smallskip
\n
We now prove (\ref{fikg1}). Since the two relations are conjugate to 
each other it is enough to show the first one.
We have
$$
i \, ( d' d''^{*}\,  +\,  d''^{*} d')  \, = \, 
d' (i\,d''^{*}) \, + \, (i \, d''^{*}) d' \, =\,
d' ( L^{*} d'  \, - \, d' L^{*}) \, + \, 
( L^{*} d'  \, - \, d' L^{*})  d' \, = 
$$
$$
=\, d' L^{*} d' \, -\, d'^{2} L^{*} \, + \, L^{*} d'^{2} \, - \, 
d' L^{*} d' \, = \, 0.
$$

\smallskip
\n
We now prove (\ref{fikg2}), i.e. that 
that $ \Delta_{d} = \Delta_{d'} + \Delta_{d''}.$
$$
\Delta_{d}\, = \, (d'+d'')(d'^{*} + d''^{*}) + 
(d'^{*} + d''^{*})(d'+d'') \, = \, 
$$
$$
( d' d'^{*} + d'^{*}d') \, + \, ( d''d''^{*} + d''^{*}d'') \, + \, 
( d' d''^{*}\, +\, d''^{*} d' )  \, + \, ( d''d'^{*} \, + \, d'^{*} 
d'')   \, = \, \Delta_{d'} \, + \, \Delta_{d''} \, +\,  0  \, +\,0.
$$
We now prove that $\Delta_{d'} = \Delta_{d''}.$
$$
-i \, \Delta_{d'} \, = \, -i \, d'd'^{*} \, - \, i\, d'^{*}d' \, = \, 
d' (-i \, d'^{*}) \, + \, (-i\, d'^{*} ) d' \, = 
$$
$$
= \, d' ( L^{*} d'' \, - \, d'' L^{*} ) \, + \, 
( L^{*} d'' - d'' L^{*} ) d' \, = \, 
d'L^{*} d'' \, - \, d' d''L^{*} \, + \, L^{*} d''d' \, - \, 
d''L^{*} d' \, = 
$$
$$
d'L^{*} d'' \, + \, d'' d'L^{*} \, - \, L^{*} d'd'' \, - \, 
d''L^{*} d' \, = \, - (  d''L^{*} d' \, - \, d'' d'L^{*} \, + \, L^{*} 
d'd'' \, - \, 
d'L^{*} d'' ) \, =  
$$
$$
= \, - [ \,  d''(L^{*}d' \, - \, d'L^{*} ) \, + \, (L^{*}d' - d'L^{*}) 
d'' \, ] \, = \, - [ \, d'' (i \, d''^{*}) \, + \, ( i\, d''^{*}) d'' 
\, ] \, = \, - i\, \Delta_{d''}.
$$
(\ref{fikg1}) implies that $\Delta'$ commutes with $d''$ and that
$\Delta''$ commutes with $d'.$

\n
(\ref{fikg2}) implies that $\Delta$ commutes with $d',$ $d''$ and 
hence with $d = d'+d''.$

\n 
Since $\Delta$ is self-adjoint and commutes with
$d,$ $d'$ and $d'',$ $\Delta$ commutes with 
$d^{*},$ $d'^*$ and $d''^{*}.$

\n
One verifies  that $\Delta$ commutes with $L$ and $L^{*}$
by using (\ref{fikg3}).

\n
Finally, since $\Delta'$ (and $\Delta''$) are of bi-degree
$(0,0),$ (\ref{fikg2}) implies that $\Delta$ is also of type
$(0,0),$ i.e. that it preserves bi-degrees.
\blacksquare

\bigskip
\begin{exe}
    \label{commrel}
 {\rm ({\bf Proof of the commutation relations})
 Prove that
\begin{equation}
\label{crel}
 [ L^{*} , d' ] \, = \, i \, d''^{*}
 \end{equation}
as follows. See \ci{gh} pages 111-114.

\smallskip
\n
By Remark \ref{usef} it is enough to prove the relation 
for the Euclidean metric on $\comp^{n}.$

\n
Define operators
$$
\wedge_{j} \, : \, A^{p,q}_{c} \lorw  A^{p+1,q}_{c}, \qquad
\qquad \wedge_{j} (u) \, := \, dz_{j} \wedge u,
$$
$$
\ov{\wedge}_{j} \, : \, A^{p,q}_{c} \lorw  A^{p,q+1}_{c}, \qquad
\qquad \ov{\wedge}_{j} (u) \, := \, d\ov{z}_{j} \wedge u,
$$
$$
\partial_{j} \, : \, A^{p,q}_{c} \lorw  A^{p,q}_{c}, \qquad
\qquad \partial_{j} ( u_{JK} \, dz_{J}\wedge d\ov{z}_{K} ) \, := \,
 \frac{\partial \,u_{JK} }{\partial z_{j}} \, dz_{J}\wedge d\ov{z}_{K}
$$
$$
\ov{\partial}_{j} \, : \, A^{p,q}_{c} \lorw  A^{p,q}_{c}, \qquad
\qquad \ov{\partial}_{j} ( u_{JK} \, dz_{J}\wedge d\ov{z}_{K} ) \, := \,
 \frac{\partial \,u_{JK} }{\partial \ov{z}_{j}} \, dz_{J}\wedge d\ov{z}_{K}.
$$
We have 
$$
L \, = \, \frac{i}{2} \sum_{j}{ \wedge_{j} \circ \ov{\wedge}_{j} }
$$
$$
d' \, = \, \sum_{j}{ \partial_{j} \circ \wedge_{j} }, \qquad \qquad
d''\, = \, \sum_{j}{ \ov{\partial}_{j} \circ \ov{\wedge}_{j} }.
$$
Note that, concerning adjoints, we have
$$
(A\circ B)^{*} \, = \, B^{*} \circ A^{*}, \qquad
\qquad (\lambda \, A)^{*}\, = \, \ov{\lambda} \, A^{*}.
$$
Each of these operators has its own adjoint:
$$
\wedge_{j}^{*}, \quad \ov{\wedge}_{j}^{*}, \quad \partial_{j}^{*}, 
\qquad \ov{\partial}_{j}^{*},
$$
$$
L^{*} \, = \, - \frac{i}{2} \,  \sum_{j}{
\ov{\wedge}^{*}_{j} \circ \wedge^{*}_{j} 
},
$$
$$
d'^{*} \, =\,  \sum_{j}{ \wedge_{j}^{*} \circ \partial_{j}^{*} },
\qquad
\qquad
d''^{*} \, =\,  \sum_{j}{ \ov{\wedge}_{j}^{*} \circ \ov{\partial}_{j}^{*} }.
$$
We may omit the ``$\circ.$''

\medskip
\n
Use the definition of
the inner product via integration and
integration by parts to show that

\begin{equation}
\label{111}
\partial_j^* \,  =\,  - \ov{\partial}_j, 
\qquad \qquad 
\ov{\partial}_j^* \, = \,  -\partial_j.
\end{equation}

\medskip
\n
Note, also using (\ref{111}), that 
\begin{equation}
\label{222}
\mbox{ $\; \partial_j$ and $\ov{\partial}_k$ 
commute with each other, with $\wedge_j$ and   
$\ov{\wedge}_k,$ with $\wedge_j^*$ and $\ov{\wedge}_k^*.$ }
\end{equation}

\medskip
\n
The goal of what follows is to show 
that (\ref{ez3}), (\ref{ez33}) and (\ref{333}) below hold.

\medskip
\n
In what follows $f$ and $g$ are compactly supported smooth
functions on $\comp^{n}.$

\medskip
\n
Use the orthogonality properties
of the inner product on the exterior algebra
(see $\S$\ref{tim}) to verify that 
 for every index $j \notin J:$
$$
\langle \langle \, \wedge_j^* ( f \,  dz_J \wedge d\ov{z}_K) \, ,
\,
g\, dz_L \wedge d\ov{z}_M \, \rangle \rangle \, = \, 
\langle \langle \,   f \, dz_J \wedge d\ov{z}_K) \, , \,
 g\, dz_j \wedge dz_L \wedge d\ov{z}_M \, \rangle \rangle \, = \, 0
$$
for every pair of multi-indices $L$ and $M$ and deduce that

\begin{equation}
\label{3A}    
 \wedge_j^* ( f\, dz_J \wedge d\ov{z}_K ) =0, \qquad   j \, \notin  \, J,
 \end{equation}
and, in a similar way, that
\begin{equation}
\label{3AA}    
 \ov{\wedge}_j^* ( f\, dz_J \wedge d\ov{z}_K ) =0, \qquad   j \, \notin  \, K.
 \end{equation}

\medskip
\n
Use the fact that $||dz_j||^2=2$ at every point of 
$\comp^n$ (see Exercise \ref{oist}) and 
 the orthogonality properties of the inner product 
on the exterior algebra to show that
$$
\langle \langle \, \wedge_j^* ( f \, dz_j \wedge dz_J \wedge d\ov{z}_K) \, ,
\,
g\, dz_L \wedge d\ov{z}_M \, \rangle \rangle \, = \, 
\langle \langle \,   f \, dz_j \wedge dz_J \wedge d\ov{z}_K) \, , \,
g \,  dz_j \wedge dz_L \wedge d\ov{z}_M \, \rangle \rangle \, = 
$$
$$
= \, 2\, \langle \langle \,  f \, dz_j \wedge dz_J \wedge d\ov{z}_K \, ,
\,
g\, dz_L \wedge d\ov{z}_M \, \rangle \rangle \,
$$
for every pair of multi-indices $L$ and $M$.
Deduce that 

\begin{equation}
\label{3B}
 \wedge_j^* ( f \, dz_j \wedge dz_J \wedge d\ov{z}_K ) 
\,=\, 2\, f \, dz_J \wedge d\ov{z}_K, 
\qquad
j \, \notin\, J
\end{equation}
and, in a similar way, that
\begin{equation}
\label{3BB}
 \ov{\wedge}_j^* ( f \, d\ov{z}_j \wedge dz_J \wedge d\ov{z}_K ) 
\,=\, 2\, f \, dz_J \wedge d\ov{z}_K, 
\qquad
j \, \notin\, K.
\end{equation}
Clearly, the left-hand-sides of (\ref{3B}) and (\ref{3BB}) are zero
if, respectively,  $j \in J,$ $j \in K.$

\medskip
\n
Note that (\ref{3AA}) and (\ref{3BB}) imply
\begin{equation}
\label{ez3}
\wedge_{j} \ov{\wedge}^{*}_{k} \, + \, 
\ov{\wedge}^{*}_{k} \wedge_{j} \, = \, 0
\end{equation}
and, by conjugation,
\begin{equation}
\label{ez33}
\ov{\wedge}_j \wedge_k^{*} + \wedge_k^{*} \ov{\wedge}_j \, =\, 0.
\end{equation}

\medskip
\n
Use (\ref{3A}) and (\ref{3B}) to verify that
\[
\wedge^{*}_{j} \wedge_{j} ( f \, dz_{J} \wedge 
d\ov{z}_{K} )  \, =\, 
\left \{
\begin{array}{ll}
    0 &  \mbox{ if $j \, \in\,  J.$ } \\
    2 f \, dz_{J} \wedge 
d\ov{z}_{K}
&
\mbox{ if $j\, \notin \,  J$ }
\end{array}
\right.
\]
and that
\[
\wedge_{j} \wedge_{j}^{*}( f \, dz_{J} \wedge 
d\ov{z}_{K} )  \, =\, 
\left \{
\begin{array}{ll}
   2 f \, dz_{J} \wedge 
d\ov{z}_{K}
&
\mbox{ if $j\, \in \,  J$ } \\
0 &  \mbox{ if $j \, \notin\,  J.$ } 
\end{array}
\right.
\]
Deduce that
\begin{equation}
    \label{333jj}
    \wedge_{j} \wedge^{*}_{j} + \wedge^{*}_{j} \wedge_{j}
    \, = \, 
2\, Id.
\end{equation}

\medskip
\n
Let $j \neq k.$ Use (\ref{3B}) to show that
\begin{equation}
    \label{quasar}
 \wedge_{j}^{*} \wedge_{k} ( f \, dz_{j} \wedge
 dz_{J} \wedge d\ov{z}_{K} ) \, = \, 
 - \wedge_{k} \wedge_{j}^{*} ( f \, dz_{j} \wedge
 dz_{J} \wedge d\ov{z}_{K} ).
 \end{equation}

 \medskip
 \n
Let $j \neq k.$ Show  that
 \begin{equation}
     \label{qusar}
 \wedge^{*}_{j} \wedge _{k} ( dz_{J}\wedge d\ov{z}_{K} ) \, = 
 \,
 \wedge_{k} \wedge _{j}^{*} ( dz_{J}\wedge d\ov{z}_{K} )  \, = \,
 0, \qquad \quad j \, \notin \, J.
 \end{equation}

 \bigskip
 \n
Note that  (\ref{333jj}), (\ref{quasar}) and (\ref{qusar})
imply
\begin{equation}
\label{333}
\wedge_j^* \wedge_k + \wedge_k \wedge_j^*  \,
= \, 2\, \delta_{jk}\,Id
\end{equation}

\medskip
\n
Verify (\ref{crel}) using 
(\ref{111}), (\ref{222}), (\ref{ez3})
and (\ref{333}) as follows. 

\medskip
\n
$$
L^* d' \, = \, 
- \frac{i}{2} \sum_{j,k}{ \ov{\wedge}_j^* \wedge_j^* \partial_k \wedge_k}
\, = 
\mbox{  by (\ref{222}) } \, = \, 
\,  - \frac{i}{2} 
\sum_{j,k}{  \partial_k  \ov{\wedge}_j^*\wedge_j^*\wedge_k} \, =
$$
\[
 = \, - \frac{i}{2} \left[ 
\left( \sum_{j=k}{  \partial_j  \ov{\wedge}_j^*\wedge_j^*\wedge_j} 
\right)
\, + \,
\left( \sum_{j\neq k}{  \partial_k  \ov{\wedge}_j^*\wedge_j^*\wedge_k} 
\right) 
\right]  \, = \, \mbox{ by (\ref{333}) }
\]
\[
= \, 
- \frac{i}{2} \left[ \;
\left(
\sum_{j}{ - \partial_j  \ov{\wedge}_j^*\wedge_j\wedge_j^*} \, + \, 
2\, \sum_j{ \partial_j \ov{\wedge}_j^*}
\right)
\, - \, 
\left( 
\sum_{j\neq k}{ \partial_k \ov{\wedge}_j^* \wedge_k \wedge^*_j}
\right) \;
\right] \, = 
\]
$$
= \, \mbox{
by (\ref{ez3})}
\, = \,  - \frac{i}{2} \, \sum_{j}{ \partial_j \wedge_j \ov{\wedge}_j^* 
\wedge^*_j } \, - \, 
i\, \sum_{j=k}{ \partial_j \ov{\wedge}_j^{*}  } \, - \, 
\frac{i}{2} \, \sum_{j\neq k}{
\partial_k \wedge_k \ov{\wedge}_j^* \wedge_j^* } \, =
$$
$$
= \, - \frac{i}{2} \, \sum_{j,k}{ 
\partial_k \wedge_k \ov{\wedge}_j^* \wedge_j^* }
\, - \, i\, \sum_j{\partial_j \ov{\wedge}^*_j } \, = \, 
  d' L^* + i \, \sum_j{ (- \partial_j) \ov{\wedge}_j^* }  \, =
$$
$$
=\, \mbox{
by (\ref{111})
}\, = 
d'L^* + i\, d''^*.
$$
}
\end{exe}
    
\bigskip
\begin{cor}
    \label{corroc}
Let $(X,h)$ be a K\"ahler manifold. Then
$$
\fbox{$
{\cal H}_{\Delta}^{l}(X,\comp) \, = \, 
{\cal H}_{\Delta'}^l(X)  \, = \, {\cal H}_{\Delta''}^l(X) ,
$}
$$
$$
\fbox{$
{\cal H}_{\Delta'}^{p,q}(X)  \, = \, {\cal H}_{\Delta''}^{p,q}(X).
$}
$$
Set,  for $p+q=l$
$$
{\cal H}_{\Delta}^{p,q}(X) \,:= \,
{\cal H}_{\Delta}^{l}(X,\comp) \, \cap \, A^{p,q}(X).
$$
Then
$$
\fbox{$
{\cal H}_{\Delta}^{p,q}(X) \, = \, {\cal H}_{\Delta'}^{p,q}(X) 
\, = \, {\cal H}_{\Delta''}^{p,q}(X), 
$}
$$
and
$$
\fbox{$
{\cal H}_{\Delta}^{l}(X,\comp) \, = \,
\bigoplus_{p+q=l}{ {\cal H}_{\Delta}^{p,q}(X) } \, = \, 
\bigoplus_{p+q=l}{ {\cal H}_{\Delta'}^{p,q}(X) } \, = \, 
\bigoplus_{p+q=l}{ {\cal H}_{\Delta''}^{p,q}(X) }.
$}
$$
\end{cor}
{\em Proof.} 
It is a consequence
of the fact that $\Delta$ commutes with $\pi^{p,q}.$
The first  set of   two equalities follow from
(\ref{fikg2}).
The second set of two equalities
from the fact that if $u$ is $\Delta-$harmonic,
then it is $\Delta'-$harmonic ($\Delta''-$harmonic, reps.) so that the 
$(p,q)-$compo\-nen\-ts of $u$ are $\Delta'-$harmonic 
($\Delta''-$harmonic, resp.),
hence $\Delta-$harmonic.
\blacksquare

\bigskip
At this point we could state and prove the Hodge Decomposition Theorem.

\n
However, without some more preparation, 
it would not be clear from the proof that
the resulting decomposition is independent
of the K\"ahler metric. 
We prefer to prove this most important fact along the way.

\subsection{The Hodge Decomposition for compact K\"ahler manifolds}
\label{thdkm}
The following lemma clarifies the relation among the various notions
of  forms being ``closed'' in the compact 
K\"ahler context. It  is used to show 
that the Hodge decomposition is canonical, i.e. independent of the 
metric.

\begin{lm}
\label{many}
({\bf $d'd''-$Lemma})
Let $(X,h)$  be  \underline{compact} K\"ahler manifold  
and $u \in A^{p,q}(X)\, \cap \, \ke{\,d}.$

\n
The following are equivalent.

\smallskip
\n
$(a)$ $u$ is $d-$exact;

\smallskip
\n
$(b')$ $u$ is $d'-$exact;

\smallskip
\n
$(b'')$ $u$ is $d''-$exact;

\smallskip
\n
$(c)$ $ u$ is $d'd''-$exact;

\smallskip
\n
$(d)$  $u \in ({\cal H}^{p,q}(X))^{\perp}.$
\end{lm}
{\em Proof.}
Since $u=d(d''v)= d'(d''v)=d''(-d'v),$ we see that 
$(c)$ implies $(a),$ $(b')$ and $(b'')$ on any complex manifold.

\n
The rest of the proof needs the compactness as well
as the K\"ahler assumption.

\n
Since $(X,h)$ is K\"ahler, we have $\Delta'=  \Delta''= \frac{1}{2} 
\Delta$  so that 
${\cal H}^{p,q}(X) :=
{\cal H}^{p,q}_{\Delta'}(X) = {\cal H}^{p,q}_{\Delta''}(X)=
{\cal H}^{p,q}_{\Delta} (X).$ In what follows we need the compactness
of $X$ to ensure that the inner products are defined.
If $(a)$ holds, then $u=dv$ and $\langle \langle
dv, w \rangle \rangle= 
\langle \langle
v, d^{\star} w \rangle \rangle= 0$ for every $w  \in {\cal  H}^{p,q}(X),$
since $\Delta-$harmonic forms are characterized by Lemma
\ref{dd*closed}.

\n
This shows  that $(a)$ implies $(d).$

\n
An analogous argument shows that $(b')$ implies $(d)$ and that $(b'')$ 
implies $(d).$

\n
We are left with proving that $(d)$ implies $(c).$

\n
Since  $du=0$ and $u$ is of pure type
$(p,q),$ we have that $d'u=d''u=0.$ 

\n
By Theorem \ref{hodtcav} for $d'',$ $u = d''u_{1} + d''^{\star} u_{2}$
for a $u_{1}$ of type $(p,q-1)$ and  a $u_{2}$
of type $(p, q+1).$ 

\n
Since $d''u=0,$ we have
$d''d''^{\star} u_{2} =0.$ 

\n
On the other hand,
$0 = \langle \langle d'' d''^{\star} u_{2}, u_{2}
\rangle \rangle = || d''^{\star} u_{2} ||^{2},$ i.e. $d''^{\star} u_{2}=0$
and $u = d'' u_{1}.$

\n
Theorem \ref{hodtcav} for $d'$ implies that
$u_{1} = h + d'w_{1} + d'^{\star} w_{2},$
where $h \in {\cal H}^{p,q-1}(X),$
$w_{1}$ of type $(p-1,q-1)$ and $w_{2}$ of type
$(p+1, q-1).$

\n
Using (\ref{fikg1}), we get
$$
u = d''u_{1}\, = \,  d''h+ d''d'w_{1}+ d''d'^{\star}w_{2}
\, = \,
-d'd''w_{1} - d'^{\star}d'' w_{2}.
$$
Since $d'u=0,$ the equality above implies
that $d'd'^{\star} (d''w_{2})=0.$  

\n
The same argument as 
above, proving that $d''^{\star} u_{2}=0,$ shows that $d'^{*}d''w_{2}=0,$
i.e. $u = d'd'' (-w_{1}).$
\blacksquare

\bigskip
\begin{tm}
 \label{thdtwh}
 ({\bf The Hodge Decomposition})
Let $X$ be a compact  K\"ahler manifold.
The natural maps
$$
H^{p,q}_{BC}(X) \lorw H^{p+q}_{dR}(X, \comp)
$$
are injective. We denote their images by $H^{p,q}(X)
\subseteq H^{p+q}_{dR}(X, \comp).$

\n
There is an internal direct sum decomposition
$$
\fbox{$ 
H^{l}_{dR}(X, \comp)\, = \, \bigoplus_{p+q=l}{H^{p,q}(X)}
$}
$$
satisfying the equality
$$
\fbox{$
H^{p,q}(X) \, = \, \ov{ H^{q,p}(X) }.
$}
$$
\end{tm}
{\em Proof.}
Consider the  diagram
$$
\begin{array}{ccc}
 H^{p,q}_{d''}(X)  & \stackrel{g}\lorw & {\cal H}^{p,q}_{\Delta''}(X)   \\
 \uparrow f & & \downarrow h \\
 H^{p,q}_{BC}(X) & \stackrel{l}\lorw & H^{p+q}_{dR}(X, \comp). 
\end{array}
$$
The maps $f$ and $l$ stem from  the definitions. See
Exercise \ref{mapsfbc}.
The map $g$ is the inverse of the Hodge Isomorphism
Theorem \ref{hodtcav}  for $d''$  and it
selects the $\Delta''-$harmonic representative
for a Dolbeault class.
The map $h$ is defined by the fact that
a $\Delta''-$harmonic class is $\Delta-$harmonic, hence
$d-$closed.

\medskip
\n
CLAIM I: {\em $f$ is bijective}.
Injectivity follows from Lemma \ref{many}, $(b'')=(c).$
Surjectivity follows from the fact that every 
element $a_{D} \in H^{p,q}_{d''}(X)$ can be represented
by a $\Delta''-$harmonic  $(p,q)-$form $a$ which is then necessarily 
$d-$closed
so that the associated class $a_{BC} \in H^{p,q}_{BC}(X)$ maps 
to $a_{D}.$

\medskip
\n
CLAIM II: {\em the diagram commutes}.
Let $a \in A^{p,q}(X) \cap \ke{\, d}.$
Let $a_{BC} \in H^{p,q}_{BC}(X),$ $a_{D} \in H^{p,q}_{d''}(X)$
and $a_{dR} \in H^{p+q}_{dR}(X,\comp)$ be the corresponding
elements. Let $a' \in {\cal H}^{p,q}_{\Delta''}(X)$
be the $\Delta''-$harmonic representative
of $a_{D}$ and $a'_{dR}
\in H^{p+q}_{dR}(X, \comp)$ be the corresponding class.
We need to show that $a_{dR} = a'_{dR},$ i.e. that
$a-a'$ is $d-$closed.  Note that
$a' \in A^{p,q}(X) \cap \ke{\, d}.$ We have that
$f( (a-a')_{BC}) =  0.$ By CLAIM I, $(a-a')$ is $d'd''-$closed.
We conclude by Lemma \ref{many}, $(a)=(c).$

\medskip
\n
By 
Corollary \ref{corroc}, $h$ is injective.
It follows that $l$ is injective and
that the image
of  ${\cal H}_{\Delta''}^{p,q}(X)$ in
$H^{p+q}_{dR}(X, \comp),$
being the image of $H^{p,q}_{BC}(X) $ which
depends only on the complex structure of $X,$ is independent
of the K\"ahler metric employed.

\n
By
Corollary \ref{corroc}.
we can assemble, for all pair of indices
$(p,q)$ such that $l= p+q,$ the isomorphisms
$l$
and reach the desired conclusion.
\blacksquare

\bigskip
\begin{rmk}
\label{mohdt}
 {\rm
 ({\bf The meaning of the Hodge Decomposition})
The Hodge Decomposition expresses a basic and important
property of the complex cohomology of
compact K\"ahler manifolds. It does {\em not} state
that any complex-valued de Rham class $a_{dR}\in
H^{p+q}_{dR}(X, \comp )$ admits a canonical
representative $a \in A^{l}(X)$ and associated
$(p,q)-$com\-po\-nen\-ts. It states that
$a_{dR}$ can be canonically decomposed into 
a sum $\sum{ a_{dR}^{p,q}},$
whose terms  one calls the $(p,q)-$components
of $a_{dR}$ and which do not depend
on the choice of a K\"ahler metric.
The choice of a K\"ahler metric allows to
represent these $(p,q)-$components using
$\Delta-$harmonic $(p,q)-$forms which depend on the metric.

\n

\n
This canonical decomposition 
of $H^{l}_{dR}(X, \comp )  $ and the canonical isomorphism
$H^{l}_{dR}(X, \comp) \simeq H^{l}(X, \zed) \otimes_{\zed}
\comp$ lead to the notion of pure Hodge structure
of weight $l.$ 
 }
\end{rmk}

\bigskip
\begin{rmk}
    \label{nonatm}
    {\rm
For a general complex manifold $X$, there is no natural map
$H^{p,q}_{d''}(X) \lorw H^{p+q}_{dR}(X, \comp):$
the map $h\circ g$ depends on the choice of a metric
and $f$ is not, in general, invertible.

\n
In the compact K\"ahler case, we have shown that
$k: =h\circ g= l \circ f^{-1}$ is well-defined and depends only on the 
complex structure on $X$ and we have
$$
H^{p,q}(X) \, := \, l(H^{p,q}_{BC}(X)) \, = \, 
k(H^{p,q}_{d''}(X))\, =\,   
h( {\cal H}^{p,q}_{\Delta''}(X)).
$$   
    }
  \end{rmk}

\subsection{Some consequences}
\label{scof}
The following Theorem is trivially false on a complex manifold and
the example of the  Iwasawa manifold Example \ref{fsik} shows 
that it is false in general on compact complex manifolds.
However, it is true on compact complex surfaces.
\begin{tm}
\label{holcl}
({\bf Holomorphic vs. closed forms})
Let $X$ be a \underline{compact} K\"ahler manifold and
$u \in H^{0}(X, \Omega^{p}_{X})$ 
be a holomorphic $p-$form, i.e. a $(p,0)-$form such that
$d''u=0.$

\n
Then $u$ is closed.
\end{tm}
{\em Proof.} 
Let $h$ be a K\"ahler metric on $X.$ 
Since $d''^{*}u$ has bi-degree
$(p,-1),$ it must be the zero form. 
 It follows that
$u$ is $\Delta''-$harmonic, hence $\Delta-$harmonic, hence
$d-$closed. (See the remark after the Proof of Lemma
\ref{dd*closed})
\blacksquare

\bigskip
Let us assemble some of the properties
we have proved 
concerning the Hodge numbers
$$
h^{p,q} (X) \, := \, \dim_{\comp}{\, H^{q}(X, \Omega^{p}_{X} ) } \, = \,
\dim_{\comp}{ \, H^{p,q}_{d''}(X) } \, = \,
\dim_{\comp}{ \, H^{q,p}_{d'}(X)}.
$$
for compact complex manifolds and for compact K\"ahler manifolds.
Note that we do not have $h^{p,q}(X) = h^{q,p}(X),$ unless
$X$ is K\"ahler, e.g. the Hopf surface for which
$h^{0,1}(X) = h^{1,0}(X) +1 =1.$ Note that in this case the so-called
Hodge-Fr\"olicher spectral sequence collapses, as it does for any
compact complex surface, so that even this last condition does not
imply $h^{p,q}(X)= h^{q,p}(X).$ It is only equivalent to $b_{l}(X)=
\sum_{p+q=l}{\, h^{p,q}(X)}.$

\bigskip
In what follows we omit the K\"unneth-type relations. See
\ci{gh}, page 105.

\begin{tm}
\label{pqid}
Let $X$ be a compact K\"ahler manifold of dimension $n.$
We have
\begin{equation}
\label{pqid1}
h^{p,q}(X) \,   < \, \infty,
\end{equation}
\begin{equation}
\label{pqid2}
h^{0,0} (X) \, = \, h^{n,n}(X) \, = \,  1, 
\end{equation}
\begin{equation}
\label{pqid3}
h^{p,q} \, = \, h^{n-p, n-q}(X).
\end{equation}
Assume in addition that $X$ is K\"ahler. Then we also have
\begin{equation}
\label{pqid4}
h^{p,q}(X) \, = \, h^{q,p}(X),
\end{equation}
\begin{equation}
\label{pqid5}
b_{l}(X) \, = \, \sum_{p+q=l}{ \, h^{p,q}(X) },
\end{equation}
In particular, $b_{odd}(X) \, = \, even.$
\begin{equation}
\label{pqid6}
h^{p,p}(X) \, \geq  \, 1, \qquad     0 \leq  p \leq n.
\end{equation}

\end{tm}
{\em Proof.} Exercise. \blacksquare

\bigskip
These relations can be visualized conveniently on the so-called
{\em Hodge diamond} \ci{gh} page 117.

\bigskip
There are also some inequalities stemming from the Hard Lefschetz 
Theorem \ref{chl}.
See \ci{dem1}, Corollaire 8.17 for the precise statements, which are
immediate  consequences of Theorem \ref{chl}.b.

\bigskip

\begin{exe}
 {\rm
 Compute $h^{p,q}(\pn{n}).$
 }
 \end{exe}

 \bigskip
 
 \begin{exe}
 {\rm
Compute $h^{p,q}(Q^{n}),$ where
 $Q^{n} \subseteq \pn{n+1}$ is the nonsingular
 $n-$dimensional quadric.
 }
 \end{exe}

\newpage

\section{Lecture 7: The Hard Lefschetz Theorem and the Hodge-Riemann
Bilinear Relations}
\label{thlthrbr}
We discuss Hodge structures and polarizations, the operation of cupping
with the fundamental class of a hyperplane section, the Hard Lefschetz 
Theorem,  the Hodge-Riemann Bilinear Relations and the Weak Lefschetz 
Theorem (which is also called 
the Lefschetz theorem on hyperplane sections).

\bigskip
The statement of the Hard Lefschetz Theorem
for a compact K\"ahler manifold $(X,\omega)$ of dimension $n$
is perhaps surprising at first.

\medskip
\n
Poincar\'e Duality asserts that 
$$
H^{n-j}(X,\comp)  \, \simeq \, H^{n+j}(X, 
\comp)^{\vee}
$$
In particular, $b_{n-j}(X) = b_{n+j}(X).$

\medskip
\n
The Hard Lefschetz Theorem states that
$$
L^{j}: H^{n-j}(X, \comp) \lorw
H^{n+j}(X, \comp), \qquad \qquad  u \, \lorw \,  \omega^{j} \wedge  u
$$
is an isomorphism.

\n
This isomorphism, coupled with Poincar\'e Duality,
gives rise to
the non-degenerate  bilinear form 
$$
H^{n-j}(X, \comp) \times H^{n-j}(X,\comp) \lorw 
\comp, \qquad \qquad  (u,v) \lorw \int_{X}{\omega^{j}\wedge u \wedge v}.
$$
The Hodge-Riemann Bilinear Relations express the beautiful
signature properties of this bilinear form.

\bigskip
Though the Hard Lefschetz Theorem and the Hodge Riemann Bilinear relations
hold for any compact K\"ahler manifold $(X, \omega),$
in the last two lectures we discuss the projective case
(a projective manifold is compact and K\"ahler)
which  has some additional geometric flavor
with respect to the compact  K\"ahler case: 
 one can choose  $\omega$  to be the fundamental class
of a hyperplane section on $X.$

\bigskip

 \subsection{Hodge structures}
 \label{hs}
The theory of Hodge structures  and polarizations is a  powerful tool
to investigate the topology of algebraic varieties.
It has been employed fruitfully 
to study families of smooth projective varieties (variations
of polarized pure Hodge structures).

 \n
The fact, discovered by Deligne,  that  the singular cohomology
of  a complex algebraic variety  is endowed with a canonical and 
functorial mixed Hodge structure is one of beautiful depth.

\bigskip
The Hodge Decomposition Theorem and the Hodge-Riemann Bilinear 
Relations can be conveniently re-formulated in terms of
Hodge structures and polarizations.

\bigskip
 
 Let $l \in \zed$, $H$ be a finitely generated abelian group, 
$H_{\rat}:= H 
\otimes_{\zed} \rat$, $H_{\real}= H \otimes_{\zed}\real,$
$H_{\comp}= H \otimes_{\zed} \comp$.

\bigskip
\begin{defi}
    \label{defhs}
    {\rm({\bf Hodge structure})
A {\em pure Hodge structure  of weight $l$ on $H$,} or on $H_{\rat}$, 
$H_{\real},$
is a direct sum decomposition 
$$
\fbox{$
H_{\comp}= \bigoplus_{p+q=l}{ H^{p,q}}
$}
$$
such that 
$$
\fbox{$ H^{p,q} \, = \, \overline{ H^{q, p} }. $}
$$
}
\end{defi}

\bigskip
One can change the labeling of the $H^{p,q}$ spaces and change the weight
$l$ to any other weight $l'$ with $l'\equiv l$ $mod \,2.$

\bigskip
The {\em Hodge  filtration} on such a structure
  is the \underline{decreasing} filtration 
$$
F^{p}(H_{\comp})\, :=
\, \bigoplus_{p'\geq p}{ H^{p',l - p'}.}
$$

\bigskip
A {\em morphism of pure Hodge structures}
$f: H \lorw H'$ of the same weight $l$
 is a group homomorphism such that
$f \otimes Id_{\comp}$ is compatible with the Hodge filtration,
i.e. such that it is a filtered map, i.e. such that
$$
f( F^{p}(H_{\comp})) \subseteq F^{p}(H'_{\comp}).
$$
Such maps are automatically 
{\em strict} in the sense that
$$
\im { (\, F^{p} (H_{\comp}) \, ) } \, = \, \im{ (\, 
H_{\comp} \, ) } \cap F^{p} ( H'_{\comp} ).
$$
The category of  pure Hodge structures of weight $l$ is abelian.
In particular, kernels and cokernels of maps of 
pure Hodge structures of weight $l$ are pure Hodge structures of
weight $l.$

\bigskip
\begin{exe}
\label{fnstr}
{\rm
Give examples of filtered maps of vector spaces  which are not strict.
}
\end{exe}

\bigskip
\begin{ex}
\label{eahs}
{\rm
With reference to $\S$\ref{teaov}, the vector space
$\Lambda^{l}_{\comp}( V^{*}_{\comp})$ is a pure Hodge structure of 
weight $l.$
If $f: X \to Y$ is a holomorphic map of complex manifolds,\
then, for every $x \in X,$ the induced map
$f^{*}: \Lambda^{l}_{\comp}(T^{*}_{Y, f(x)}(\comp))
\lorw  \Lambda^{l}_{\comp}(T^{*}_{X, x}(\comp))$
is a map of Hodge structures.
}
\end{ex}

\bigskip
\begin{ex}
    \label{exfohs}
    {\rm
    The Hodge Decomposition Theorem \ref{thdtwh} 
    can be re-formulated by stating that
    if $X$ is a compact K\"ahler manifold, then 
    $H^{l}(X, \zed)$ admits a natural
structure of  a pure Hodge structure of
    weight $l.$
    The reader should check that this structure is
    functorial  with respect to holomorphic
    maps of compact K\"ahler manifolds.
    }
    \end{ex}

\bigskip
\begin{exe}
\label{cfu}
{\rm 
Verify the statements of Example \ref{eahs} and  \ref{exfohs}.
}
\end{exe}

\bigskip
\begin{exe}
    \label{hsl}
    {\rm
Let $X$ be compact K\"ahler and $L$ be a complex  line bundle
on $X.$ Denote its first Chern class and the operation
of cupping with this first Chern class simply by $L.$

\n
Observe that $L$ has bidegree $(1,1)$ as an operator
on $H^{\bullet}(X, \comp)$ in the sense that it maps
a class of degree $(p,q)$ to one of degree $(p+1, q+1).$

\n
Re-define the Hodge structure of weight $l+2r$ on $H^{l+2r}(X, \zed)$
so that $L^{r}: H^{l}(X, \zed) \lorw H^{l+2r}(X, \zed)$
becomes a morphism of Hodge structures of weight $l$.
Do the analogous operation for weight $l+2$ and for
any weight $l+2m,$ $m \in \zed.$

\n
Show directly that $\ke{\,L^{r}}$, $\im{\,L^{r}}$ and $\mbox{Coker} \,
L^{r}$ are  pure Hodge structures of 
weight $l.$ 

}
\end{exe}

\bigskip
Let $C$ be the Weil operator, i.e. 
$$ 
C\,:\, H_{\comp} \, \simeq \, H_{\comp}, \qquad
C(x)\, =\, i^{p-q} x, \quad  \forall \, x\in H^{p,q}.
$$
See $\S$\ref{twop}.
The Weil operator  is real, i.e. it fixes
$H_{\real}.$  

\medskip
\n
Replacing $i^{p-q}$ by 
$z^{p}{\overline{z}}^{q}$
 we get a real  action $\rho$ of $\comp^{*}$ on $H_{\real}$ and
 $C = \rho(i).$
 
\begin{defi}
\label{defpol}
{\rm 
A {\em polarization} of the real pure Hodge structure
$H_{\real}$ of weight $l$ is a  real bilinear form $\Psi$ on $ H_{\real}$
which is  invariant under the action given by $\rho$ restricted to 
$S^{1} \subseteq \comp^{*}$
and such that the bilinear form
$$
\fbox{$ \widetilde{\Psi} (x,y)\, := \,  \Psi(x, C(y))$}
$$ 
is symmetric and positive definite.
}
\end{defi}

\begin{exe}
    \label{note}
    {\rm
    Verify the  following statements.
If $\Psi$ 
is a polarization
 of the real pure Hodge structure
$H_{\real}$ of weight $l,$ 
then $\Psi$ is symmetric
if $l$ is even, and antisymmetric if $l$ is odd.
(Hint: $C^2 = (-1)^l.$)

\n
In any case, $ \Psi$ is non-degenerate.

\n
In addition, for every $0 \neq x \in H^{p,q}$, 
$(-1)^{l}i^{p-q}\Psi(x,\overline{x}) >0,$
where, by abuse of notation, $\Psi$ also denotes 
the $\comp-$bilinear extension of $\Psi$ to $H_{\comp}.$
(Hint: consider $x+\ov{x}$ and $\frac{x-\ov{x}}{i}$)
}
\end{exe}

\begin{rmk}
 \label{important}
 {\rm
 If $H' \subseteq H$ is a Hodge sub-structure,
then $H'_{\real}$ is fixed by  $C$ so that
 $\Psi_{|H'_{\real}}$ is a polarization.
 
 \n
 In particular, $\Psi_{|H'_{\real}}$  is non-degenerate.
}
\end{rmk}

\subsection{The cup product with the Chern class of a hyperplane
bundle}
\label{tchb}
Let $M$ be a smooth manifold.
An $l-$form $\omega \in E^l(M)$ defines a linear map
$$
L\, : \, E^{\bullet}(M) \lorw E^{\bullet + l}(M)\qquad
u \lorw \omega \wedge u.
$$
If $\omega $ is closed, then $L$ descends to de Rham cohomology:
$$
L\, : \, H^{\bullet}_{dR}(M,\real) \lorw H^{\bullet + l}_{dR}(M,\real),
\qquad
u \lorw \omega \wedge u.
$$

\bigskip
Since the de Rham isomorphism $H^{\bullet}_{dR}(M, \real)
\simeq H^{\bullet}(M, \real)$ is an algebra isomorphism,
we may write $u \cup b$ or $u \wedge v$ depending
on our taste and notational convenience.

\bigskip
If, in addition, $M$ is compact and oriented of dimension $m,$
then we also have the Poincar\'e
Duality isomorphism which we denote
$H^{\bullet}(M, \real) \stackrel{PD}\simeq H_{m -\bullet}(M, \real).$

\n
We denote the inverse map by the same symbol.

\bigskip
We have the following geometric interpretation
of the cap, cup (wedge) products in connection with set-theoretic
intersections.

\medskip
\n
The {\em cap} product
$$
\cap \, : \, H^q(M, \zed) \otimes_\zed H_p(M,  \zed) 
\lorw H_{p-q}(M, \zed)
$$
satisfies, taking $\rat-$coefficients, denoting cohomology classes by
Roman letters and homology classes by Greek ones:
$$
\fbox{$
a \cap \beta \, = \, ( a \cup \beta^{PD} )^{PD} \, = \, 
a^{PD} \cap \beta,$}
$$
where $a^{PD} \cap \beta$ is the result of the following operations:
pass from $a$ to $a^{PD};$ select representatives for $a^{PD}$ and for
$\beta$ so that their supports meet transversally;
the correctly oriented intersection is a well-defined homology class.

\bigskip
In other words:
\bigskip

\n
\centerline{
\fbox{ the cap/cup products corresponds to transverse intersections 
via Poincar\'e Duality.}
}

\bigskip
\begin{exe}
\label{intn}
{\rm
Let $L$ be the $y-$axis in $\real^2$ oriented by $dy$ and $C$ 
be the unit circle centered about the origin
oriented counterclockwise. 

\n
Verify that the intersection numbers, with respect to the 
standard orientation for $\real^2,$ are
$i_{ (0,1)}( L \cap C) =-1$ and $i_{ (0,-1)}( L \cap C) =1$ so that the total
intersection number is zero.

\n
Do the analogous  thing in $\comp^2$ and verify that the total
intersection number is $2.$

\n
Verify that if two complex varieties of complementary
dimensions in $\comp^n$ meet transversally
at a point, then the intersection number at that point is positive.

\n
In fact, this is true even if they do not meet transversally
but still meet at isolated points.

\n
See \ci{gh}, pages 49-65 for a nice discussion on intersecting
cycles.
}
\end{exe}

\bigskip
\begin{rmk}
\label{onint}
{\rm ({\bf Borel-Moore Homology})
Two distinct lines $L_1$ and $L_2$ 
in $ \pn{2}_{\comp}$ meet at one point, say $P,$ with intersection 
index one.

\n
Singular cohomology records this fact as follows.

\n
The lines give rise to cohomology classes
$cl(L_i) \in H^2(\comp \pn{2}, \zed)$
and we get 
$$
\deg{ (cl (L_1) \cup cl(L_2)) }=1
$$
where $\deg : H^4(\comp \pn{2}, \zed ) \simeq \zed$
is the isomorphism stemming from the standard orientation
of $\pn{2}.$

\n
Let $Q \in  \pn{2}$ be a point on neither line.
We have $H^4( \comp \pn{2} \setminus Q, \zed)  =\{0 \}.$
The pairing in cohomology is not the  correct
one to deal with intersections in a non-compact ambient manifold.

\n
Moreover, if we consider non-compact subvarieties
and their intersections, e.g. two lines in $\comp^2$
meeting transversally at a point,
we see that  the non-compactness of the subvarieties
also plays a role in
the loss
of information

\n
It turns out that the formalism of Borel-Moore Homology,
with its natural pairing with singular homology, allows for a 
reasonably general intersection theory.
Given two complex subvarieties $V$ and $W$ of a complex 
manifold $X$ of dimension $n$ we have a paring, depending
on $X,$
$$
H^{BM}_j(V, \rat) \times H_k(W, \rat) \lorw H_{2n -j -k}^{BM}
(V \cap W,\rat)
$$
with a natural degree map giving intersection numbers
if $j+k =2n.$

\n
The cycles  in the Borel-Moore theory
are locally finite combinations of simplices so that it is natural to expect
a pairing with the cycles is homology
which are given by a finite number of simplices.
}
\end{rmk}

\subsection{The Hard Lefschetz Theorem and the Hodge-Riemann 
Bilinear Relations}
Let $X \subseteq  \pn{N}$ 
be a nonsingular complex projective variety not contained in any
hyperplane.

\medskip
\n
The space of holomorphic sections of the line bundle
$$
\fbox{$
L \, : = \, {L_1}_{|X},  $}
$$
where $L_1$ is the hyperplane line bundle of Exercise
\ref{proj}, contains a subspace $\cal L$ naturally identified with the space
$\comp[x_0, \ldots, x_N ]_1$ of homogeneous linear polynomials 
in $N+1$ variables.

\medskip
\n
If $s \in { \cal L} \setminus \{0\}$ and $P_s$ is the corresponding polynomial,
then, counting multiplicities, the zero sets
$$
\{ s=0 \} \, = \, \{ X \cap \{P_s=0 \} \}.
$$

\bigskip
We have the Fubini-Study metric $h_{FS}$
on $\pn{N}$ with associated K\"ahler form  $\omega_{FS}.$

\medskip
Restricting to $X$ we get a K\"ahler metric $h$ with associated 
K\"ahler  form $\omega: = {\omega_{FS}}_{|X}$ and we have,
by slight abuse of notation:
$$
 \omega \, = \, (X \cap {\Bbb P}^{N-1})^{PD} \, = 
 \, c_1 (L) \; \in H^2(X, \rat).
$$ 

\bigskip
It follows, via Poincar\'e Duality, that

\bigskip
\fbox{ cupping with $\omega$
corresponds to intersecting with a hyperplane section.}

\bigskip
\begin{exe}
\label{clcve}
{\rm
Let $C \subseteq \pn{N}$ be a nonsingular complex projective curve.
\n
Show that
$$
L^r \, : \, H^{1-r}(C, \rat)  \, \simeq \, H^{1+r} (C, \rat), \quad
\forall \, r \geq 0. 
$$
Show that
$$
\int_C{ \omega \wedge u \wedge \ov{u} } \, > 0, \quad \forall 
\, u \, \in \, H^{0,0}(C).
$$
$$
i^{p-q} \, \int_C{ u \wedge \ov{u} } \, > 0, \quad \forall 
\, u \, \in \, H^{p,q}(C), \; p+q=1.
$$
}
\end{exe}

\begin{exe}
\label{clsf}
{\rm
Let $S \subseteq \pn{N}$ be a nonsingular complex algebraic surface.

\n
Show that
$$
L^r \, : \, H^{2-r}(C, \rat)  \, \simeq \, H^{2+r} (C, \rat), \quad
r =0, \, 2. 
$$
Deduce that
$$
H^2(S, \rat) \, = \, L ( H^0(S, \rat) ) \,\stackrel{\perp}\oplus\, \ke{L}
$$
where the direct sum is orthogonal with respect to the natural
cup product pairing on $H^2(S, \rat),$ each direct summand
is a pure Hodge structure of weight $2$ and 
$\dim_{\rat}{ \ke{L} } = b_2 - b_0 = b_2 -1.$

\n
Let $S \subseteq \pn{3}$ be a nonsingular quadric. Show that
the restriction of the cup product to $\ke{\,L} \subseteq H^{2}(S, 
\rat)$ is negative definite. Show that the same is true no matter
which ample line bundle $L$  one chooses on $S.$
}
\end{exe}

\bigskip
The methods of Exercise \ref{clsf} do not seem to be helpful
in studying the map
$$
L\, : \, H^1(S, \rat) \lorw H^3(S, \rat).
$$

\bigskip
Perhaps, at this point,
the reader unfamiliar with the Hard Lefschetz Theorem \ref{chl}.a will
appreciate its  beauty.

\bigskip
Let
$  0 \leq r \leq n$ and
define the {\em space $ P^{n-r}_{L}$ of  rational $(n-r)-$primitive classes}
as
$$
P^{n-r} \, : = \; \ke{ \, L^{r+1} } \, \subseteq H^{n-r}(X, 
\rat).
$$
Clearly, this space can be defined also for $\zed,$ 
$\real$ and $\comp$ 
coefficients, in which case we use the notation
$P^{n-r}_L(X, \real),$ for example.

\bigskip
\begin{exe}
\label{parpq}
{\rm Show that  
$$
P_L^{j}(X, \comp)
= \bigoplus_{p+q=j}{\, (\, P^{j}_L(X, \comp) \, \cap H^{p,q}(X)}\,), 
\quad 0 \leq j \leq 
n,
$$
i.e. show that $P^j_L \subseteq H^j(X, \rat)$ is a pure sub-Hodge structure
of weight $j.$
}
\end{exe}

\bigskip

\begin{tm}
\label{chl}
(a) ({\bf The Hard Lefschetz Theorem.})
Let $X\subseteq \pn{N}$ be a projective manifold of dimension $n.$

\n
For every $r \geq 0$, we have isomorphisms
$$
L^{r}\, :\,  H^{n-r}(X, \rat) \lorw H^{n+r}(X, \rat).
$$

\smallskip
\n
(b) ({\bf The primitive Lefschetz decomposition.})
For every $r \geq 0$ there is the direct sum decomposition
$$
H^{n-r}(X, \rat ) \,  = \,   \stackrel{\perp}\bigoplus_{j \geq 0}\,  
L^{j}P^{n-r-2j}_L
$$
where each summand
is a  pure Hodge sub-structure of weight $n-r$   and
all summands are mutually orthogonal with respect to the
bilinear form $\int_{X}{ \omega^{r} \wedge - \wedge -}.$

\smallskip
\n
(c) ({\bf The Hodge-Riemann bilinear relations.})
For every $0\leq l \leq  n,$ the bilinear form
$(-1)^{\frac{l(l+1)}{2}}\int_{X}{\omega^{r} \wedge - \wedge -}$
is a polarization of the real pure weight $l$ Hodge structure
$P^{l}_L(X,\real) \subseteq H^{l}(X, \real).$ In particular,
$$
(-1)^{\frac{l(l-1)}{2}} i^{p-q} \int_{X}{\omega^{n-l}\wedge
\alpha \wedge \overline{\alpha} }>0, \qquad \forall \; 0 \neq \alpha
\in P^{l}_L(X, \comp) \,  \cap\,  H^{p,q}(X, \comp).
$$
\end{tm}

\bigskip
\begin{rmk}
    \label{p2chl}
    {\rm 
    Theorem \ref{chl}.a was first stated by Lefschetz
    in his 1924 book. It is known that his ``proof'' contains a gap.
    The first proof is due to Hodge and uses 
    what today we call Hodge Theory. 
    For a proof based on the commutation relations
    and how they relate to a $sl_{2}(\comp)-$action
    on $H^{\bullet}(X, \comp)$ see \ci{gh}, pages 118-122,
    or \ci{dem1}, pages 35-39 and 44-45. A proof
    of 
   Theorem \ref{chl}.c can be found in \ci{we}, \ci{wl} or \ci{vo}.
    
   \n
    In his second paper on the Weil conjectures, Deligne 
    has given an algebraic proof of Theorem
    \ref{chl}.a. It uses algebraic geometry in characteristic $p.$
    See $\S$\ref{tltas} for a discussion.
    It seems unlikely  that those methods will
    give a proof of 
    Theorem \ref{chl}.c. 
    }
    \end{rmk}

\bigskip
\begin{rmk}
    \label{refhtwp}
    {\rm
  In view of Remark \ref{important},
the bilinear form $\Psi,$ which is non-degenerate
on the primitive spaces, is in fact non-degenerate when restricted
to any sub-Hodge structure
of the primitive spaces.
}
\end{rmk}

\bigskip
\begin{exe}
 \label{linalghl}
 {\rm
Show that  Theorem \ref{chl}.a   is equivalent to
the statement that
$L^{j}: H^{n-r}(X, \rat) \lorw H^{n-r+2j}(X, \rat)$
is injective for $0 \leq j \leq r.$

\n
Show that  Theorem \ref{chl}.a  is equivalent to
the statement that
$L^{j}: H^{n+r-2j}(X, \rat) \lorw H^{n+r}(X, \rat)$
is surjective  for $0 \leq j \leq r.$

\n
Show that 
Theorem \ref{chl}.a  implies Theorem \ref{chl}.b. 
  }
 \end{exe}

\bigskip
\begin{exe}
\label{hleqbnd}
{\rm 
Show that   Theorem \ref{chl}.a  is equivalent
to the bilinear form
$$
\Psi ( u,v) \, : = \, \int_{X}{ \omega^{r} \wedge u \wedge v}
$$
on $H^{n-r}(X, \rat)$ being non-degenerate.
}
\end{exe}

\begin{exe}
    \label{hlpn}
 {\rm Prove the Hard Lefschetz Theorem when $X = \pn{n},$ using 
 only the
 product structure on the cohomology ring.
 }
\end{exe}

\begin{ex}
    \label{resth}
    {\rm Let $X \lorw Y$ be the blowing-up of a point  $y \in Y$ on
    a
     smooth projective three-fold $Y,$ $A$ be an ample line bundle
     on $Y$ and $L:= f^{*}A.$
     Then the statement of the Hard Lefschetz Theorem \ref{chl}.a does 
     {\em not} hold for $L.$ In fact, if $D =f^{-1}(y)$ is the 
     exceptional divisor, then $L\cup [D] =0 $  so that $L: H^{3-1}(X, 
     \rat)
     \lorw H^{3+1}(X, \rat)$ is not an isomorphism.
     }
     \end{ex}
     
     \begin{rmk}
    \label{hlsemi}
    {\rm What follows is proved in \ci{demig1}.
    Let $f: X \lorw Y$ be a map of complex projective varieties, $X$
    nonsingular, $A$ ample on $Y,$ $L:= f^{*}A.$
    
    \n
    Then the statement of the Hard Lefschetz Theorem \ref{chl}.a holds
    for $L$
    {\em if and only if} the map $f$ is {\em semi-small}.
    The notion of semi-smallness was introduced by
    Goresky and MacPherson to generalize the Lefschetz Theorem on 
    Hyperplane Sections \ref{cwlt} to the case of
    not necessarily proper complex analytic maps $f : X \to 
    \pn{N}.$  A map $f: X \lorw Y$ is semi-small iff $\dim{X \times_{Y} X}
    = \dim{X}.$ Note that in general, 
    we have the inequality ``$\geq$.''
    These maps occur frequently in complex geometry and in
    representation theory. If $f$ is semi-small, then one also shows 
    that all three  statements of Theorem \ref{chl} hold for $f.$
    }
    \end{rmk}

\subsection{The Weak Lefschetz Theorem}
\label{twlt}
The following is another fundamental result about the topology of projective 
varieties. Essentially, it states that a great deal of the topology
of a projective manifold ``comes'' from its hyperplane sections.

\bigskip
Let $Y  := X \cap \pn{N-1} \subseteq X \subseteq 
\pn{N} $ be a hyperplane section of
a projective manifold  $X$ of complex dimension $n$
embedded in projective space.

\bigskip
\begin{tm}
\label{cwlt}
({\bf The Weak Lefschetz Theorem.})
The natural restriction map
$$
r^{*}: H^{j}(X, \rat) \lorw H^{j} (Y, \rat) 
$$
is an isomorphism for $j \leq n-2$ and is injective for $j=n-1.$

\n
If, in addition, $Y$ is nonsingular, then the natural Gysin map
(i.e. the Poincar\'e dual to the map in homology)
$$
\tilde{r}_*: H^{n+j-2}(Y) \lorw H^{n+j}(X)
$$
is an isomorphism for $j \geq 2$ and is surjective for $j=1.$
\end{tm}
{\em Sketch of proof.} For a proof and other references, see \ci{gh},
page 156.

\n
There is an anaologus result for the homotopy groups, but we do not
state it here.

\n
The complex manifold $X \setminus Y$ is a closed
complex submanifold of $\comp^N$ of dimension $n.$

\n
Andreotti and Frankel have observed that
such a manifold has the homotopy type of a CW-complex of real dimension 
at most $n$ so that 
$$
H^j(X \setminus Y, \rat) = \{0 \}, \qquad \forall \; j \, >\, n.
$$

\n Equivalently, 
$$
H^j_c (X \setminus Y, \rat)\,  = \, 
\{ 0 \},  \qquad \; \forall \; j\, <\, n.
$$

\n
To conclude one uses  the long exact sequence of
relative cohomology
$$
 \ldots  \lorw H^j_c(X \setminus Y, \rat) \lorw H^j(X, \rat) 
\lorw H^j(Y, \rat)  \lorw H^{j+1}_c (X \setminus Y, \rat) \lorw 
\ldots
$$
\blacksquare

\bigskip
\begin{exe}
\label{rama}
{\rm 
C.P. Ramanujan has elegantly shown that
the Weak Lefschetz Theorem and the Kodaira Vanishing Theorem,
i.e. $H^{j}(X, L^{\vee}) = \{0 \}$ for $j <n,$
are equivalent.

\n
Prove this fact by looking at the long exact sequences
stemming from the following 
two exact sequences:
$$
0 \lorw \Omega^{p}_{X} \otimes L^{\vee}  \lorw
\Omega^{p}_{X} \stackrel{r}\lorw {\Omega^{p}_{X}}_{|Y } \lorw  0,
$$
$$
0 \lorw \Omega^{p-1}_{Y} \otimes L^{\vee}_{|Y}  \lorw 
{\Omega^{p}_{X}}_{|Y}  
\stackrel{i}\lorw
\Omega^{p}_{Y}
\lorw  0.
$$
You will need the  functoriality
of the Hodge Decomposition, the Dolbeault isomorphism,
and the fact
$i \circ r$ induces the restriction of $p-$forms from $X$ to $Y.$
 See \ci{gh}, page 157.
}
\end{exe}

\bigskip
Theorem \ref{cwlt} imposes severe restrictions on the topology of 
algebraic varieties.
For example, it implies that if $X$ is a surface then
the natural map of fundamental groups
$$
\pi_{1}(Y,y) \lorw \pi_{1}(X,y)
$$
is surjective. If $X$ is not simply connected, then 
a smooth hyperplane section $Y$ cannot be 
isomorphic to $\pn{1}.$

\n
If $X$ is simply connected and it has dimension at least
three, then 
a smooth $Y$ must be simply connected and the natural
restriction map
$$
Pic(X) \lorw Pic(Y)
$$
of Picard groups (isomorphisms classes of holomorphic line bundles
with product induced  by the tensor product)
is injective with torsion free-cokernel and if the dimension
of $X$ is at 
least four, then  the restriction map is
an isomoprhism, i.e. every line bundle
on $Y$ comes from one on $X.$

\newpage
\section{Lecture 8: Mixed Hodge structures, the Semi-simplicity
Theorem and  the approximability of primitive  vectors}
\label{l8}
We discuss the mixed Hodge structure on the singular cohomology
of a complex algebraic variety, the Semi-simplicity Theorem, 
the Global Invariant Cycle Theorem, the connection between
the Semi-simplicity Theorem and the Hard and Weak Lefschetz Theorems.
We conclude with illustrating an approximation technique
for primitive vectors.

\subsection{The mixed Hodge structure on the cohomology of
complex algebraic varieties}
\label{mhs}
The paper \ci{du} is a  nice introduction to mixed Hodge structures.

\bigskip
The following is Deligne's fundamental result
which endows  in a functorial way
the singular cohomology of  a complex algebraic variety
with a mixed Hodge structure. See \ci{ho2} and \ci{ho3}.

\bigskip
The implications of the existence of this rich  structure 
 are numerous and far-reaching, for they impose severe
 constraints of varieties and the maps between them.
 We try to get a glimpse of this powerful machinery in this last 
 lecture.
 
\bigskip
\begin{tm}
\label{mhssc}
Let $X$ be an algebraic variety. For each $j$ there is an increasing
weight filtration
$$
\fbox{$
\{0 \} \, =\,  W_{-1} \, \subseteq W_{0} \, \subseteq \, \ldots 
\, \subseteq W_{2j} \, = \, H^{j}(X, \rat)$}
$$
and a decreasing Hodge filtration
$$
\fbox{$
H^{j}(X, \comp) \, = \, F^{0} \, \supseteq \, F^{1} \, \supseteq 
\ldots \, \supseteq\,  F^{m}\, \supseteq \, F^{m+1} \, = \, \{0 \}
$}
$$
such that the filtration induced by $F^{\bullet}$ on the 
complexified graded pieces
of the weight filtration endows 
every graded piece $W_{l} / W_{l-1}$ with a pure Hodge structure of 
weight $l.$

\n
This structure is functorial for maps of algebraic
varieties and the induced maps strictly preserve
both filtrations.
\end{tm}

\bigskip
Here is a very partial list of the  properties
of this structure. See \ci{du} and \ci{ho3}, $\S$8.2 for some more
properties.

\begin{itemize}
    
\item
If $X$ is projective and smooth, then
$$
0 \, = \, W_{j-1} \, \subseteq \, W_{j} \, =\, H^{j}(X, \rat)
$$
and one says that {\em $H^{j}(X, \rat)$ has pure weight $j.$}

\item
If $X$ is projective (but not necessarily smooth), then
$$
W_{j} \, = \, H^{j}(X, \rat)
$$
and one says that  {\em $H^{j}(X, \rat)$ has  weights $\leq j.$}

\item
If $X$ is smooth (but not necessarily projective), then
$$
 \{ 0 \}  \, = \, W_{j-1} \, \subseteq \, H^{j}(X, \rat)
$$
and one says that  {\em $H^{j}(X, \rat)$ has  weights $\geq j.$}

\item
If $i:U \hookrightarrow X$ is the inclusion of a Zariski-dense open subset
of a smooth projective manifold, then
$$
W_{j}(H^{j}(U, \rat)) \, =\, i^{*} H^{j}(X, \rat) \, = \, i^{*} 
\, W_{j}(H^{j}(X, \rat)).
$$
\end{itemize}

\bigskip
A beautiful application of the theory of mixed Hodge structures
is the following theorem of Deligne's whose proof
is a nice and simple exemplification of the ``yoga of weights.''

\bigskip
\begin{tm}
\label{iisa}
Let $\,Y \lorw U \hookrightarrow X$ be maps
of complex algebraic varieties, with $Y$ proper,
$X$ proper and nonsingular,
$U$ a Zariski-dense open subset of $X.$

\n
Then the natural images of $H^j(X,\rat) $ and of $H^j(U,\rat)$ in
$H^j(Y, \rat)$ coincide.
\end{tm}
{\em Proof.} See \ci{ho2}, Corollaire 3.2.18 and  \ci{ho3},
Proposition 8.2.6.

\n
By strictness we have
$$
\im{ ( H^{j}( U, \rat) ) } \, \cap\,  W_{j}( H^{j}(Y, \rat) )\, = \, 
\im { ( W_{j}( H^{j}( U, \rat) )  )}.
$$
Since $H^{j}(Y,\rat)$ has weights $\leq j,$ we have
that
$$
\im{ ( H^{j}( U, \rat)  ) } \, = \, 
\im { ( W_{j}( H^{j}( U, \rat) ) )}.
$$
Since $W_{j}( H^{j}( U, \rat))  = i^{*}H^{j}(X, \rat)$
we get the desired equality.
\blacksquare

\bigskip
\begin{ex}
\label{cexe}
{\rm
Let $S^1 \lorw \real^2 \setminus
\{0 \} \subseteq S^{2}$ be the obvious inclusion
of real algebraic varieties.
The conclusion of Theorem \ref{iisa} does not hold
in this situation.

\n
Theorem \ref{iisa} expresses another property peculiar
to {\em complex} algebraic geometry.
}
\end{ex}

\subsection{The Semi-simplicity Theorem}
\label{sst}
A {\em local system} on a (connected) variety $Y$
is a locally constant sheaf ${\cal L}$
of finite dimensional rational vector spaces on $Y.$

\bigskip
Let ${\cal L}$ be a local system on $Y,$ $y \in Y$
and ${\cal L}_y$ the stalk at $y.$

\medskip
\n
The vector space ${\cal L}_y$ is a representation
of the fundamental group $\pi_1 ( Y, y)$ and
the sheaf ${\cal L}$ is the sheaf of sections 
of the vector bundle
$\widetilde{Y} \times_{\pi_1(Y,y)} {\cal L}_y,$
where $\widetilde{Y}$ is the universal cover of $Y.$

\medskip
\n
From this description one deduces a canonical isomorphism
\begin{equation}
\label{lsin}
H^{0}(Y, {\cal L}) \, \simeq \, {\cal L}_{y}^{\pi_{1}(Y,y)}
\end{equation}
where the right-hand-side denotes the subspace
of ${\cal L}_{y}$ of vectors fixed by $\pi_{1}(Y,y).$

\bigskip
\begin{ex}
\label{ehr}
{\rm Let $f:X\lorw Y$ be a proper,
 surjective and smooth map of smooth manifolds
such that $df$ has everywhere maximal rank.

\n
By Ehreshmann Lemma, \ci{dem1}, Lemme 10.2, $f$ is locally topologically 
trivial
over $Y,$ i.e.
all fibers are diffeomorphic and, locally over $Y,$ $X$ is a product
of $Y$ times the typical  fiber.

\n
Circuiting around  closed paths centered at a point  $y \in Y$
produces  a local system on $Y$ by considering the
action of $\pi_1(Y,y)$ on the cohomology
groups $H^j(f^{-1}(y), \rat)$ of the fiber at $y.$

\n
The sheaf on $Y$ so obtained  is canonically isomorphic
to the $j-$th higher direct image 
$R^j f_* \rat_X$ of the constant sheaf  $\rat_X.$
}
\end{ex}

\bigskip
\begin{defi}
\label{defssple}
{\rm
A local system
$\cal L$ on an algebraic variety $Y$ is said to be
{\em semi-simple} if every local subsystem
${\cal L}'$ of $\cal L$ admits a complement, i.e. a  local 
subsystem 
${\cal L}''$ of $\cal L$ such that 
$$
{\cal L} \simeq {\cal L}' \oplus 
{\cal L}''.
$$
}
\end{defi}

\bigskip
\begin{ex}
\label{enss}
{\rm 
The local system ${\cal L}$ on  the circle $S^{1}$
given by the representation 
\[
\rho\, : \, \pi_{1}(S^{1}, q)\, 
\simeq\, \zed\,  \lorw
GL(2,\rat), \qquad  
\rho (1) := 
 \left(
\begin{array}{cc}
    1&1\\
    0&1
    \end{array}
    \right) \]
is indecomposable, i.e. it cannot be written as a direct sum
of non-trivial local sub-systems, but it is  {\em not} semi-simple
since it admits the local sub-system generated by the 
column vector
$(1,0).$

\n
This local system appears in connection with the standard 
degeneration over the unit disk $\Delta \subseteq  \comp$
of a smooth elliptic curve to a rational nodal curve.
}
\end{ex}

\bigskip
The following result of Deligne's
summarizes some very important properties of smooth
projective maps.

\bigskip
\begin{tm}
\label{dss}
 ({\bf The Semi-simplicity Theorem for smooth maps})
Let $f: X^{n} \lorw Y^{m}$ be a smooth, proper  and surjective  
map of nonsingular complex quasi-projective 
varieties of the indicated dimensions and $L$ be a
an ample
line bundle on $X.$
Then 
$$
L^{i}\, :\,  R^{n-m-i} f_{*} \rat_{X} \, \simeq \,  R^{n-m+i} f_{*}
\rat_{X}, \qquad \quad \forall i \geq 0,
$$
\begin{equation}
 \label{dt}
 Rf_{*}\rat_{X} \, \simeq \, 
\bigoplus_{i\geq 0} R^{j}f_{*}\rat_{X} [-j] 
\end{equation}
and
the local systems $R^{j}f_{*}\rat_{X}$
are semi-simple on $Y.$
\end{tm}
{\em Proof.} See \ci{dess} and  \ci{ho2}, Th\'eor\`eme 4.2.6.
\blacksquare

\bigskip
\begin{rmk}
\label{dssi}
{\rm
(\ref{dt}) implies the $E_2-$degeneration of the Leray spectral 
sequence for maps as in Theorem \ref{dss}.
This is another special feature of complex algebraic geometry.
}
\end{rmk}

 \subsection{The Leray spectral sequence}
 \label{tlss}
 \begin{tm} 
 \label{lip}
 ({\bf The Leray spectral sequence})
 Let $f: X \to Y$ be a continuous map of topological spaces
 and $G$ be a sheaf of abelian groups on $X$.
 
 \n
 There is a first quadrant spectral sequence
 $$
 E_{2}^{p,q} \, = \, H^{p}(Y, R^{q}f_{*}G) \, \Longrightarrow
 H^{p+q}(X, G).
 $$
 \end{tm}
 {\em Proof.} See \ci{gh}, \ci{b-t} or \ci{dem2}.
 \blacksquare
 
 \bigskip
 Let us spell out the content of Theorem \ref{lip}.
 
 \bigskip
For every pair of indices
$(p,q),$ the abelian group $H^{p+q}(X, G)$ admits a decreasing filtration
 associated with the map $f$:
 $$
 H^{p+q}(X, G) \, =  \, F^{0} \, \supseteq \, F^{1} \, 
 \supseteq  \, \ldots  \, \supseteq \, F^{p+q} \, \supseteq \,  \{0 \}.
 $$
 
 \medskip
 There is a collection of abelian groups $E^{p,q}_{r}$
 and  of group homomorphisms $d_{r}^{p,q},$ $p,\, q \, 
 \in \, \zed,$ $r \geq 2$ such that:

 \medskip
 \n
-- $E^{p,q}_{r} = \{0 \}$ as soon as either $p<0$ or $q<0$
(this is what ``first quadrant'' means);
 
\smallskip
\n
-- $E^{p,q}_{2} = H^{p}(Y, R^{q}f_{*}G);$

\smallskip
\n
-- $d_{r}^{p,q}: E^{p,q}_{r} \lorw E^{p+r,q-r+1}_{r};$
note that for any fixed $p+q,$ $d_{r}^{p,q} =0$ for every  $r \gg 0;$
 
\smallskip
\n
-- $d_{r}^{2} =0;$

\smallskip
\n
-- $E^{p,q}_{r+1} = \ke{\, d_{r}^{p,q} }/ \im{ \, d_{r}^{p-r, q+r-1} };$
note that for any fixed $p+q,$  $E_{r}^{p,q} = E_{r+1}^{p,q}$
for  every $r \gg 0;$ we denote these stabilized groups by 
$E_{\infty}^{p,q};$
 
\smallskip
\n
-- 
 $
 F^{p}/F^{p+1} = E_{\infty}^{p,q}.
 $
 
 \bigskip
One says that the spectral sequence degenerates
at $E_r$ if $d_{r'}=0$ for every $r'\geq r.$

\bigskip

 \begin{exe}
     \label{edge}
  {\rm 
   Show that there are canonical injective and surjective maps
 $$
 H^{p+q}(X, G) \surj F^{0}/F^{1}= E^{0,p+q}_{\infty} \inj
 E^{0, p+q}_{r},
 $$
 $$
 E^{p+q,0}_{r} \surj E^{p+q}_{\infty} = F^{p+q} / \{ 0 \} 
 \inj H^{p+q}(X,G).
 $$
 These maps are called the {\em edge homomorphism}.
 
 \n
Show that one always has  exact sequences
 $$
  0 \lorw H^{1}(Y, f_{*}G)  \lorw H^{1}(X, G) \lorw E_{\infty}^{0,1} 
  \lorw 0,
  $$
  $$
  0 \lorw H^{1}(Y, f_{*}G) \lorw H^{1}(X, G) \lorw H^{0}(Y, 
  R^{1}f_{*}G) \lorw  H^{2}(Y, f_{*}G) \lorw H^{2}(X, G).
  $$
  }
  \end{exe}

  \bigskip
  \begin{exe}
  \label{sph}
  {\rm 
  Let $f: S^{n} \lorw S^{m}$ be a  fiber bundle with fiber
  $S^{n-m},$ i.e. a bundle structure on a sphere, over a sphere
  in spheres.
  Use the Leray spectral sequence to determine the possible values
  of $(n,m).$ (Hint: for $m=1$ use the fundamental groups; for
  $m\geq 2$ use the fact that $S^{m}$ being simply connected implies
  that the locals systems $R^{j}f_{*}\zed_{S^{n}}$ are constant
  (cf. \ci{b-t}, Theorems II.13.2 and II.13.4). Find examples
  of such bundles.
  }
\end{exe}

\subsection{The Global Invariant Cycle Theorem}
\label{tfp}
Let $f: U \to S$ be a  proper and smooth morphism of 
smooth complex algebraic varieties.
In particular, the fibers of $f$ are smooth, possibly disconnected,
compact smooth  complex manifolds all diffeomorphic to each other.

\bigskip

As discussed in Example \ref{ehr}, we have the {\em monodromy representation}
$$
\pi_{1}(S, s) \lorw Aut(H^{j}( f^{-1}(s), \rat )).
$$

\bigskip
\n
Every cohomology class on $U$ restricts to a $\pi_{1}(S,s)-$invariant
class in $H^{j}(f^{-1}(y),\rat).$

\begin{tm}
    \label{tfpt}
    ({\bf The Global Invariant Cycle Theorem})
Let  $U\subseteq X $ be a smooth compactification of $U.$

\n
Then the natural restriction map
$$
H^{j}(X, \rat) \lorw  H^{j}(f^{-1}(s) , \rat)^{\pi_{1}(S,s)}
$$
is surjective, i.e the image is given by  the invariants.
\end{tm}
{\em Proof.} By (\ref{dt}), Remark \ref{dssi},
Exercise \ref{edge} and
(\ref{lsin}),  the Leray spectral sequence
for $f$ is $E_{2}-$degenerate and we get a natural surjection
$$
H^{j}(U, \rat) \lorw H^{0}(S, R^{j} f_{*}\rat_{U}) \simeq
H^{j}(f^{-1}(s), \rat)^{\pi_{1}(S,s)}.
$$
The conclusion follows from Theorem \ref{iisa}.
\blacksquare

\bigskip
A typical situation to which this theorem applies is 
when one has a  proper surjective morphism
$g: X \to Y$ with $X$ smooth. Then one sets $S\subseteq Y$
to be the locus over which $g$ is smooth, $U:= g^{-1}(S)$
and $f: = g_{|U}.$

\medskip
\n
Every cohomology class on $X$ restricts to a $\pi_{1}(S,s)-$invariant
class in $H^{j}(f^{-1}(y),\rat).$

\medskip
\n
The Global Invariant Cycle Theorem establishes the highly
non-trivial fact, false in general in 
a non-K\"ahler context, that {\em every}
invariant class on the fiber is the restriction of a global class
on $X.$

\subsection{The Lefschetz Theorems and semi-simplicity}
\label{tltas}
We now show that the Weak Lefschetz Theorem
and the Semi-simplicity Theorem
imply the Hard Lefschetz Theorem.

\bigskip
We need the following  simple topological fact which
is a consequence of the discussion of $\S$\ref{tchb}.

\bigskip
\begin{exe}
\label{csg}
{\rm
Let $X \subseteq \pn{N}$ be a projective manifold, $Y$ be the 
transverse intersection of $X$ 
with a general hyperplane, $r: Y \to X$ be the natural map,  $L:=
{L_{1}}_{|X}$, $L':= L_{|Y}.$

\n
Show that the following diagram is  commutative:
\begin{equation}
\label{gysin}
\begin{array}{ccc}
H^{\bullet}(X, \rat) & \stackrel{L^{j}}\lorw & 
H^{\bullet +2j}(X, \rat)  \\
\downarrow r^{*}  & & \uparrow \widetilde{r}_{*} \\
H^{\bullet} (Y, \rat) & \stackrel{L'^{j-1}}\lorw & 
H^{\bullet + 2j -2} (Y , \rat).
\end{array}
\end{equation}
}
\end{exe}

\bigskip
\begin{exe}
\label{scn-1}
{\rm
Let $f: V \lorw W$ be a linear map of
finite dimensional vector spaces, $f^{*}: W^{*} \lorw V^{*}$
the map dual to $f,$ $\phi'$ a non-degenerate bilinear form
on $W,$  $\phi: W \simeq W^{*}$ the resulting isomorphism
and $l: = f^{*} \circ  \phi \circ f:$
$$
\begin{array}{ccc}
    V & \stackrel{l}\lorw & V^{*} \\
    \downarrow f & & \uparrow f^{*} \\
    W & \stackrel{\phi}\lorw  & W^{*}.
\end{array}
$$
Prove that 
$$
\ke{\,   f^{*} \circ \phi }  \, = \, (\im{f})^{\perp_{\phi'}}.
$$
Deduce that, given the diagram  (\ref{gysin})
for $\bullet = n-1:$ 
\begin{equation}
    \label{cucu}
\begin{array}{rcl}
 H^{n-1}(X, \rat) & \stackrel{L}\lorw &  H^{n+1}(X, \rat) \\
 \searrow r^{*} & & \nearrow \widetilde{r}_{*} \\
\,  & H^{n-1}(Y), & \,
 \end{array}
\end{equation}
we have that
$$
\ke{ \,\widetilde{r}_{*}  } \, = \, (\im{\, r^{*}   })^{\perp}.
$$
Show that $L$ is an isomorphism iff 
$$
\ke{ \,\widetilde{r}_{*}  }\,  \cap\,  \im{\, r^{*}   } 
\, = \, \{0 \}
$$
iff the restriction of the intersection form
on $H^{n-1}(Y, \rat)$ to the injective  image of
$$
H^{n-1}(X, \rat)
 \lorw H^{n-1}(Y,\rat)
$$
is non-degenerate.
}
\end{exe}

\bigskip
We also need the following elementary algebraic fact.

\bigskip
\bigskip
\begin{exe}
\label{grt}
{\rm Let $V$ be a finite dimensional representation of a group
$\pi,$ $\psi$ a $\pi-$invariant  bilinear form on $V,$ 
$V^{\pi} \subseteq V$ be  the  subspace of $\pi-$invariant vectors.

\n
Assume that $V$ is completely reducible, i.e. that every
$\pi-$invariant subspace $V'$  of $V$ admits a complement $V'',$
i.e. a $\pi-$invariant subspace $V''$ such that
$V = V' \oplus V''.$

\n
Show that if $\psi$ is non-degenerate, then $\psi_{|V^{\pi}}$ is 
non-degenerate.
(Hint: consider a complement $W$ of $V^{\pi},$ assume
that $V^{\pi}$ and $W$ are not orthogonal. Deduce that
$W$ admits a $\pi-$invariant subspace of dimension $1$ on which
$\pi$ acts trivially. This would be a contradiction).
}
\end{exe}

\bigskip
Let $s$ be the section of $L$ giving the smooth $Y \subseteq X,$
$s'$ another general section whose smooth zero-locus $Y'$ meets
$Y$ transversally at $Z: = Y \cap Y'.$ 
Let $p: \widetilde{X} \lorw X$ be the blowing up of
$X$ along $Z,$ $q: \widetilde{X} \lorw \pn{1}$
be the resulting map associated with the pencil
of hyperplane sections associated with $s$ and $s'.$
We have a canonical embedding 
$ Y \subseteq \widetilde{X}$ as the fiber of $q$ over the point,
still called $s,$ of the pencil corresponding to $s.$
 
\bigskip

\begin{exe}
 \label{expenc}
{\rm 
Verify that the images of $H^{j}(\widetilde{X}, \rat)$ and of
$H^{j}(X, \rat)$ in $H^{j}(Y, \rat)$ coincide.
}
\end{exe}

\bigskip
Let $U \subseteq \pn{1} $ be the locus over which 
$q$ is smooth. Note that $ s \in U.$

\bigskip
After all these preliminaries we can now draw the bridge between
the Weak and the Hard Lefschetz Theorem.

\bigskip
\begin{pr}
\label{bridge}
Assume that the Hard Lefschetz Theorem has been proved for 
$Y.$ 

\n
Then the Hard Lefschetz Theorem for $X$
is equivalent to the statement
that the intersection form on $H^{n-1}(Y, \rat)$ is non-degenerate
when restricted to the injective image
of $H^{n-1}(X, \rat)$ in $H^{n-1}(Y, \rat).$

\n
This last statement, in turn, 
is implied  either  by (a) the Semi-simplicity Theorem, or
by (b) the Hodge-Riemann Bilinear Relations on $Y.$
\end{pr}
{\em Proof.}
We specialize the diagram (\ref{gysin})
to
$$
\begin{array}{ccc}
H^{n-j}(X, \rat) & \stackrel{L^{j}}\lorw & 
H^{n+j}(X, \rat)  \\
\downarrow r^{*}  & & \uparrow \widetilde{r}_{*} \\
H^{(n-1)-(j-1)} (Y, \rat) & \stackrel{L'^{j-1}}\lorw & 
H^{(n-1)+(j-1) } (Y , \rat).
\end{array}
$$
for $j =0$ and $j \geq 2.$

\n
The Weak Lefschetz Theorem and the hypotheses
imply that $L^{j}$ is an isomorphism for $j =0$ and for $j \geq 2.$

\n
The critical case is $j=1$ which gives the diagram (\ref{cucu}).

\n
By taking a pencil as above,  applying the Global Invariant Cycle 
Theorem \ref{tfpt}
and Exercise
\ref{expenc} we deduce that
$$
H^{n-1}(X, \rat)  \, = \, H^{n-1}(Y, \rat)^{\pi_{1}(U, s)}.
$$
By Exercise \ref{scn-1}, we are left with proving that
the restriction of the intersection form on $H^{n-1}(Y, \rat)$
to (the injective image of) $H^{n-1}(X,\rat)$ is non-degenerate.

\n
The Semi-simplicity Theorem, i.e. condition (a), implies that
$H^{n-1}(Y, \rat)$ is a completely reducible
$\pi_{1}(U, s)-$re\-pre\-sen\-ta\-ti\-on.

\n
The sought-for non-degeneration follows from Exercise \ref{grt}.
So that the Hard Lefschetz Theorem for $X$ follows using (a).

\n
The Primitive Lefschetz Decomposition on $Y$
gives
$$
H^{n-1}( Y, \rat) \, = \, P^{n-1}_{L'} \, \stackrel{\perp}\oplus
L' H^{n-3}(Y, \rat).
$$
Using the Weak Lefschetz Theorem we obtain
\[
H^{n-1}(X) \, = \, \left(  P^{n-1}_{L'} \cap H^{n-1}(X) \right)
\, \stackrel{\perp}\oplus
L'\, H^{n-3}(Y, \rat).
\]
Since the intersection form on $H^{n-1}(Y, \rat)$ is non-degenerate
by Poincar\'e Duality, it is non-degenerate
on the direct summand $L'\, H^{n-3}(Y, \rat).$
By orthogonality the sought-for non-degeneration
is equivalent to the intersection form being non-degenerate
on the space  $P^{n-1}_{L'} \cap H^{n-1}(X).$

\n
This space is a sub-Hodge structure of the pure Hodge structure
$ P^{n-1}_{L'} \subseteq H^{n-1}(Y, \rat)$ which is polarized by
the intersection form by virtue of the Hodge-Riemann Bilinear relations
on $Y.$

\n
By Remark \ref{important} this polarization restricts to a polarization
of the Hodge sub-structure 
$P^{n-1}_{L'} \cap H^{n-1}(X)$ so that the intersection form
is non-degenerate on this space; see Remark \ref{note}.
\blacksquare

\bigskip
Proposition \ref{bridge} does {\em not} yield a proof of the Hard Lefschetz 
Theorem, which is a topological statement concerning the cup product 
operation with the fundamental class of a hyperplane section, 
independently of Hodge Theory.

\medskip
\n
In fact, in both cases (a) and (b), the known transcendental ways to 
prove that these assumptions are met  is through Hodge theory,
in fact through the Hard Lefschetz Theorem!

\bigskip
In his second paper on the Weil
Conjectures, Deligne has given an algebraic
proof of the    Hard Lefschetz Theorem for the $l-$adic 
cohomology of a nonsingular projective varieties over algebraically 
closed fields of \underline{positive} characteristic.

\n
The result implies  the Hard Lefschetz Theorem \ref{chl}.a.
However, presently, there is no ``algebraic'' proof of this result
that does not make use of algebraic geometry in 
positive characteristic.

\subsection{Approximability for the space of primitive vectors}
\label{apv}
In this final section I would like to  discuss, in  a simple case,
two  of the techniques introduced in the papers
 \ci{demig1} and \ci{demig2}.
These papers deal with the new structures
on  the rational singular 
cohomology of a projective manifold
that arise in connection with  a
projective  map $f: X \lorw Y.$

\bigskip
The  set-up of this section  is as follows.

\n
Let $f: X \to Y$ be  a map of projective varieties, $X$ 
be 
nonsingular, $L$ be  a hyperplane bundle on $X$, $H$ be a hyperplane
bundle on $Y,$ $M: = f^{*}H.$ 

\bigskip
The first technique  is used to deal, in some cases,
with the natural class map
$$
H_{\bullet}(f^{-1}(y)) \lorw H^{2n -\bullet}(X), \qquad \qquad y\,  
\in\,  Y.
$$

\bigskip

The second is an approximation  technique 
for   certain cohomology classes
which are ``primitive'' with respect to a 
line bundle  $M$  which is generated by its global sections,
but that is  not necessarily a hyperplane bundle,
on a projective 
manifold using  cohomology classes  which are primitive with respect
to certain related hyperplane  line bundles.

\bigskip
The upshot is that via this approximation it is possible to
study some of the properties of the intersection form
on these ``primitive'' spaces.

\bigskip
In order to exemplify these techniques, we discuss
a very special, yet meaningful case:

\bigskip
\begin{tm}
 \label{grmu}
 ({\bf Contractibility criterion})
 Let $f: X \to Y$ a  map of projective varieties of
 even  complex dimension $2m,$ $X$ be nonsingular,
 $y \in Y$ be such that $f^{-1}(y)$ is a union of finitely many
 algebraic varieties  $E_{j}$ in $X$ of dimension $m$
 and $f$ is an isomorphism over $Y \setminus y.$
 
 \n
 Then the intersection matrix
 $$
 (-1)^{m}|| E_{j} \cdot E_{k}||
 $$
 is symmetric positive  definite.
 \end{tm}

 \bigskip
 \begin{rmk}
     \label{notclass}
     {\rm
 Let us remark that the line bundle 
$M: =f^{*}H$ is trivial on the varieties $E_{j}$ so that it cannot be 
ample. In particular, a priori one cannot apply to it
the conclusions of the Hard Lefschetz Theorem and 
of the classical  Hodge-Riemann 
Bilinear Relations.

\n
There is a  re-formulation of
Theorem \ref{grmu} that holds for {\em every} projective map
with domain an  algebraic manifold.
The precise statement requires substantial preparation and is omitted.
See \ci{demig1} and \ci{demig2}.

\n
The case $m=1$  in Theorem \ref{grmu} 
is classical and is the celebrated Grauert-Mumford 
criterion for the contractibility of curves on surfaces.

\n
For $m\geq 2,$ Theorem \ref{grmu} seems new. See \ci{demig1} and 
\ci{demig2}. For example,  if a morphism $f:X \to Y$ contracts a surface
$E$ in a fourfold $X$ and nothing else, then $[E]^{2}>0.$
}
\end{rmk}

\bigskip
\begin{ex}
    {\rm As a  toy-model, the reader can keep in mind the  special case
    $m=1$
    of Theorem \ref{grmu} in all the considerations that follow.
   }
   \end{ex}

 \bigskip
\begin{pr}
\label{lindep}
({\bf Linear independence of the fundamental classes})
Let things be as in Theorem \ref{grmu}.
Then the fundamental classes $[E_{j}] \in H^{2m}(X, \rat)$
 of the contracted varieties are linearly independent.
 
 \n
 In particular, the natural class map
 $H_{2m}(f^{-1}(y), \rat) \lorw H^{2m}(X, \rat)$
 is an injection of pure Hodge structures of weight $2m.$
 \end{pr}
 {\em Proof }.
 Since $H_{2m}(f^{-1}(y), \rat)$ is a pure Hodge structure
 of weight $2m,$ spanned by the fundamental homology classes
 $\{ E_{j} \} $ (which are of type  $(m,m)$)
 we only need to prove that the class map is injective.
 
 \n
Let $ U$
 be an affine neighborhood of $y$ in $Y,$
e.g.  $U = Y \setminus H,$ where $H$ is a hyperplane section of $Y,$
 relative to some embedding, not containing $y.$ 
 
 \n
Let $U' := f^{-1}(U)$ and $g:= f_{|U'} : U' \lorw U.$ 
We look at the Leray spectral 
sequence
for $g.$

\n
We have
natural isomorphisms
\begin{equation}
\label{sky}
R^{j}g_{*}\rat_{U'} \,  \simeq \,  H^{j}(g^{-1}(y), \rat), \qquad \qquad j 
\, > \, 0.
\end{equation}
where the right-hand-side is viewed as a skyscraper sheaf
at $y.$

\n
The sheaves $R^{j}g_{*} \rat_{U'}$ are what one calls
{\em constructible sheaves} on $U,$ 
for every $j.$ In this context,
it only means that they are locally constant when restricted
to $U \setminus y$ and to $y.$

\n
Since $U$ is affine, the theorem on the cohomological dimension of 
affine sets with respect to constructible sheaves (which we do not 
state in these lectures, but which is the natural generalization
of the Andreotti-Frankel result on the homotopy type
of affine complex manifolds quoted in the proof
of the Weak Lefschetz Theorem \ref{twlt})
implies that
\begin{equation}
    \label{epqz}
    H^{p}(U, R^{0}g_{*} \rat_{U'} ) \, = \, 0, \qquad  p \, > \, 2.
 \end{equation}
 On the other end, 
(\ref{sky}) and (\ref{epqz}) imply that $d_{r}^{0,2m}=0$
for every $r\geq 2:$ these differentials land in zero groups.

\n
It follows that $E_{\infty}^{0,2m} =  E_{2}^{0,2m}  \simeq
H^{2m}(f^{-1}(y), \rat).$

\n
By Exercise \ref{edge}, the surjectivity of the  edge
map implies that 
 the natural restriction map
$$
H^{2m}(U', \rat) \lorw H^{2m}(f^{-1}(y), \rat)
$$
is surjective.

\n
By Theorem \ref{iisa}, the restriction map
$$
H^{2m}(X, \rat) \lorw H^{2m}(f^{-1}(y), \rat)
$$
is surjective as well.

\n
The dual map, i.e. the class map,
$$
H_{2m}(f^{-1}(y), \rat) \lorw H_{2m}(X, \rat)\, \stackrel{PD}\simeq \, 
H^{2m}(X, \rat)
$$
is injective.
\blacksquare

\bigskip
\begin{exe}
    \label{kflb}
    {\rm 
Let $(X, \omega)$ be a K\"ahler  manifold, 
$f: X \to \pn{N}$ be a holomorphic map.
Show that $f^*{\omega_{FS}} + \e \, \omega$ is a  K\"ahler
form 
for every $\e >0.$ 

\n
Deduce that if $f: X \to Y,$  $L,$ $H$ and $M$ 
are as in the set-up,  then one can 
represent the real classes $ M+ \e \, L:= 
c_{1}(M) + \e \, c_{1}(L) \in H^{1,1}(X)$
using K\"ahler classes (i.e. associated with K\"ahler metrics).
In particular, the Hard Lefschetz
Theorem and the Hodge-Riemann Bilinear Relations
hold  for $M+ \e \, L.$
}
\end{exe}

\bigskip
\begin{exe}
\label{hlforl}
{\rm
({\bf The Hard Lefschetz Theorem for $M$}; cf. \ci{demig1})
Show that the Hard Lefschetz  Theorem \ref{chl}.a holds 
for $M,$ i.e. show that
$$
M^{r} \, : \, 
H^{2m-r}(X, \rat)\,  \simeq \,  H^{2m+r}(X,\rat).
$$
(Hint: First prove the Weak Lefschetz
Theorem for the smooth sections $Y$ of $M$ on $X$
as in Theorem \ref{cwlt}; the relevant vanishing can be proved
along the lines of the  proof of Proposition \ref{lindep}.
Note that one can 
always choose $Y$ to avoid $f^{-1}(y).$
Use Proposition \ref{bridge}: the point is to be able to use
condition (b) in that theorem  on $Y$ which can be done since
$M_{|Y}$ is a hyperplane bundle.)

\n
Show that the Primitive Lefschetz Decomposition
Theorem \ref{chl}.b holds for $M,$ in particular, having set
$P^{2m}_{M}:= \ke{\, M} \subseteq H^{2m}(X , \rat):$
\begin{equation}
    \label{hlod}
H^{2m}(X, \rat) \, = \, P^{2m}_{M} \, 
\stackrel{\perp}\oplus \, 
M \, H^{2m-2}(X, \rat)
\end{equation}
and,  recalling Exercise \ref{kflb},
\begin{equation}
    \label{eqbn}
    \dim_{\rat}{  P^{2m}_{M} } \, =\,  b_{2m} - b_{2m-2} \,  = \,
\dim_{\real}{  P^{2m}_{M+ \e\, L} (X, \real)}.
\end{equation}
}
\end{exe}

\bigskip
\begin{pr}
    \label{apprxu}
    ({\bf Approximability of $P^{2m}_{M}$})
Let things be as in Theorem \ref{grmu}.

\n
Then, in the Grassmannian $G(  b_{2m} - b_{2m-2}, H^{2m}(X, \real) )$
we have
$$
P^{2m}_{M}(X, \real) \, = \, \lim_{\e \to 0}{
\;P^{2m}_{M+ \e\, L} (X, \real)}
$$
and  
$$
\Psi: =(-1)^{m} \int_{X}{ - \wedge -}
$$
defines a polarization of $P^{2m}_{M}(X, \real) .$
\end{pr}
{\em Proof.}
We have an elementary inclusion
of vector spaces
$$
P^{2m}_{M}(X, \real)  \, \supseteq \, \lim_{\e \to 0}{
P^{2m}_{M+ \e\, L} (X, \real)}.
$$
By (\ref{eqbn}), both sides have the same finite dimension
so that they coincide.

\n
The Hodge-Riemann Bilinear Relations for
$M+ \e \,L$ imply that $\Psi$ is a polarization
for $P^{2m}_{M+ \e \, L} (X, \real),$ for every $\e >0.$

\n
It follows that $\Psi(- , C(-))$ is positive semi-definite
on the limit
$P^{2m}_{M}(X, \real).$ 

\n
By (\ref{hlod}), since $\Psi$ is non-degenerate on $H^{2m}(X,\rat),$
it remains non-degenerate on $P^{2m}_{M}(X, \real).$

\n
Since $C$ is a real automorphism,  $\Psi(- , C(-)),$ is also  non-degenerate on
$P^{2m}_{M}(X, \real).$

\n
It follows that
$\Psi(- , C(-)),$
must be positive definite, i.e. that  $\Psi$ is a polarization
of $P^{2m}_{M}(X, \real).$
\blacksquare

\bigskip
\n
{\em Proof of Theorem \ref{grmu}.}
Note that since one can choose a section of $H$
avoiding $y,$ the injective image of the class map
lands in the $M$ primitive space
$$
H_{2m}(f^{-1}(y)) \,  \subseteq \, P^{2m}_{M}\,  \subseteq \, H^{2m}(X, \rat)
$$
as a sub-Hodge structure.

\n
By Proposition \ref{apprxu} and Remark \ref{important},
 $\Psi = (-1)^{m} \int_{X}{ -\wedge -}$
is a polarization for this image so that
$\Psi(-, C(-))$ is positive definite on
$H_{2m}(f^{-1}(y)).$ 

\n
The statement follows from the fact that
all classes under consideration are
rational and of type $(p,q)=(m,m)$
so that $C$ acts as the identity on
$H_{2m}(f^{-1}(y)).$ 
\blacksquare

\end{document}